\documentclass[aos,preprint]{imsart}
\usepackage{mathrsfs}
\usepackage{bm}
\usepackage{amsfonts}
\usepackage{natbib}
\usepackage{amsmath}
\allowdisplaybreaks[3]
\usepackage{amsthm}
\usepackage{amssymb}
\usepackage{color}
\usepackage{xcolor}
\usepackage{rotating}
\usepackage{graphicx}
\usepackage{bbm}
\usepackage{enumerate}
\usepackage{appendix}

\RequirePackage[colorlinks,citecolor=blue,urlcolor=blue]{hyperref}
\usepackage{xr}
\startlocaldefs
\newcommand{\IF}{\mathbbm{1}}
\newcommand{\argmin}{\mathrm{argmin}}
\DeclareMathOperator{\cov}{Cov}
\DeclareMathOperator{\var}{Var}

\newcommand{\hatX}{\hat{X}}
\newcommand{\hatY}{\hat{Y}}
\newcommand{\hatB}{\hat{B}}
\newcommand{\hatV}{\hat{V}}
\newcommand{\hatS}{\hat{S}}

\newcommand{\noi}{\noindent}

\newcommand{\la}{\lambda}

\newcommand{\beq}{\begin{eqnarray*}}
	\newcommand{\eeq}{\end{eqnarray*}}
\newcommand{\beqn}{\begin{eqnarray}}
	\newcommand{\eeqn}{\end{eqnarray}}

\newcommand{\mb}{\mbox}
\newcommand{\vp}{\varepsilon}

\newcommand{\ei}{\end{itemize}}

\newcommand{\be}{\begin{equation}}
\newcommand{\ee}{\end{equation}}
\newcommand{\nn}{\nonumber}

\newcommand{\ignore}[1]{}{}
\newcommand{\f}{\frac}
\newcommand{\wt}{\widetilde}
\renewcommand{\(}{\left(}
\renewcommand{\)}{\right)}
\renewcommand{\[}{\left[}
\renewcommand{\]}{\right]}

\newcommand{\de}{\delta}
\newcommand{\De}{\Delta}
\newcommand{\ga}{\gamma}

\newcommand{\si}{\sigma}

\newcommand{\e}{\mathbb{E}}

\newcommand{\WT}{\widetilde}

\renewcommand{\P}{\mathbb{P}}

\newcommand{\veps}{\varepsilon}

\newtheorem{lemma}{Lemma}[section]
\newtheorem{proposition}{Proposition}[section]
\newtheorem{theorem}{Theorem}[section]

\newtheorem{corollary}{Corollary}[section]
\newtheorem{remark}{Remark}[section]
\newtheorem{definition}{Definition}[section]
\newtheorem{condition}{Condition}[section]
\numberwithin{equation}{section}
\endlocaldefs

\begin{document}

\begin{frontmatter}
	
	% "Title of the Paper"
	\title{Refined Cram\'er Type Moderate Deviation Theorems for General Self-normalized Sums with Applications to Dependent Random Variables and Winsorized Mean}
	
	\runtitle{Refined general self-normalized moderate deviation}
	
	%\setattribute{journal}{name}{}  % suppress the "submitted to Annals of statistics"
	
	\begin{aug}
		% indicate corresponding author with \corref{}
		% \author{\fnms{John} \snm{Smith}\thanksref{a}\corref{}\ead[label=e1]{smith@foo.com}\ead[label=e2,url]{www.foo.com}}
		% \address[a]{\printead{e1};\printead{e2}}
		\author[A]{\fnms{Lan} \snm{Gao}\ead[label=e2,mark]{lgao@link.cuhk.edu.hk}},
		\author[A, B]{\fnms{Qi-Man} \snm{Shao}\ead[label=t2, mark]{shaoqm@sustech.edu.cn}}
		\and
		\author[A]{\fnms{Jiasheng} \snm{Shi}\ead[label=e3, mark]{jiashengSHI@link.cuhk.edu.hk}}
		%%%%%%%%%%%%%%%%%%%%%%%%%%%%%%%%%%%%%%%%%%%%%%
		%% Addresses                                %%
		%%%%%%%%%%%%%%%%%%%%%%%%%%%%%%%%%%%%%%%%%%%%%%
		
		\address[A]{Department of Statistics, The Chinese University of Hong Kong, \printead{e2,e3}}
		\address[B]{Department of Statistics and Data Science, Southern University of Science and Technology, \printead{t2}}

		\thankstext{t2}{The research is partially supported by National Nature Science Foundation of China NSFC 12031005, Shenzhen Outstanding Talents Training Fund, Hong Kong Research Grants Council GRF-14302515 and 14304917.}

		\runauthor{L. Gao, Q.-M. Shao and J. Shi}

	\end{aug}
	
	\begin{abstract}
		
		%We establish Cram\'er type moderate deviation theorems with significantly improved range for self-normalized weak dependent variables by study a tailored new approximation statistic, the general self-normalized sum $\sum_{i=1}^n X_i / (\sum_{i=1}^n Y_i^2)^{1/2}$ with $\{(X_i, Y_i)\}_{i=1}^n $ being a sequence of independent bivariate random vectors. 
		
		Let $\{(X_i, Y_i)\}_{i=1}^n $ be a sequence of independent bivariate random vectors. In this paper, we establish a refined Cram\'er type moderate deviation theorem for the general self-normalized sum $\sum_{i=1}^n X_i / (\sum_{i=1}^n Y_i^2)^{1/2}$, which unifies and extends the classical \cite{cramer1938nouveau} theorem and the self-normalized Cram\'er type moderate deviation theorems by \cite{jing2003self} as well as the further refined version by \cite{wang2011refined}. The advantage of our result is evidenced through successful applications to weakly  dependent random variables and self-normalized winsorized mean. Specifically, by applying our new framework on general self-normalized sum, we significantly improve Cram\'er type moderate deviation theorems for one-dependent random variables, geometrically $\beta$-mixing random variables and causal processes under geometrical moment contraction.  As an additional application, we also derive the Cram\'er type moderate deviation theorems for self-normalized winsorized mean. 
		%	simultaneously confidence interval by applying the Winsorized mean in ultra-high mean model. 
		%as an additional example of wider class of statistics and inference problems that relates to the general self-normalized sum and our result.}
		
		%Let $\{(X_i, Y_i)\}_{i=1}^n $ be a sequence of independent random vectors. A Cram\'er type moderate deviation for $\sum_{i=1}^n X_i / (\sum_{i=1}^n Y_i^2)^{1/2}$, named the self-normalized sum, is obtained, which extends earlier results by Jing, Shao and Wang (2003) and Wang (2011). Application to dependent random variables and LASSO$^{1/2}$ (Square-root LASSO) are discussed. Consequently, for dependent random variables, we obtain so far the best possible range in moderate deviation for geometric $\beta$-mixing sequence of random variables and causal process under GMC. For LASSO$^{1/2}$ associated with a high-dimensional regression model $y=X\beta_0 + \sigma_{\epsilon} \epsilon$, we obtain the bound for estimation error under least restriction on dimension $p=p(n)$ and the associate confidence level $\alpha=\alpha(n)$, i.e., $\log(p/\alpha)=o(n^{1/3})$, with the error $\epsilon$ have only finite absolute third moment.

	\end{abstract}
	
	\begin{keyword}[class=MSC]
		\kwd[Primary ]{62E17}
		\kwd{62E20}
		\kwd[; secondary ]{62J07}
	\end{keyword}

	\begin{keyword}
		\kwd{Moderate deviation}
		\kwd{Self-normalized sums}
		\kwd{Dependent random variables}
		\kwd{High-dimensional}
		\kwd{Winsorized mean}
	\end{keyword}

	% history:
	% \received{\smonth{1} \syear{0000}}
	
	%\tableofcontents
	
\end{frontmatter}

\section{Introduction}
\label{intro}
\setcounter{equation}{0}

Let $ X_1,X_2,\cdots,X_n$ be independent random variables with
$ \e X_i  = 0  $ and $  \ \e  X_i^2   < \infty  \ \mbox{for} \  i \geq 1 $. {Set
	$B_n^2 = \sum_{i = 1}^n \e  X_i^2 $},
\begin{gather*}
	S_n=\sum_{i=1}^{n}X_i, \quad V_n^2=\sum_{i=1}^{n}X_i^2,  \quad \bar{X} = \f 1n  \sum_{i = 1}^n X_i   \\
	\mbox{and} \quad \hat{\sigma}_n^2  = \frac { 1}{n - 1 }  \sum_{ i = 1 }^n (X_i - \bar{X})^2 .
\end{gather*}
The self-normalized sum is defined by $ S_n / V_n $ and is closely related to the widely-used Student's $t$ statistic $ t_n  = S_n /( \sqrt{n} \, \hat{\sigma}_n) $ in the sense that
\begin{align*}
	\P ( t_n \geq x  ) = \P  ( S_n / V_n \geq  x [ n / ( n + x^2 - 1  ) ]^{1/2} ) .
\end{align*}
Therefore, to investigate the distribution of Student's {\it t} statistic is equivalent to consider that of the less complex self-normalized statistic.

The past three decades have witnessed the flourishing development of asymptotic theory for self-normalized sums of independent random variables. Regarding the sufficient and necessary conditions for the self-normalized central limit theorem, we refer to \cite{gine1997student} and \cite{shao2018necessary} for independent and identically distributed (i.i.d.) random variables and general non-i.i.d. random variables, respectively. Specifically, for the i.i.d. case, the former paper showed that $ S_n / V_n $ is asymptotically standard normal if and only if the common distribution is in the domain of  attraction of the normal law. Basically, there are two ways to measure the accuracy in normal approximation. One method is the absolute error, which concerns the celebrated Berry-Esseen bound and Edgeworth expansion. See \cite{bentkusGotze1996berry}, \cite{bentkus1996berry} and \cite{wang1999exponential} for Berry-Esseen bounds and \cite{hall1987edgeworth} for Edgeworth expansion. Another method is the relative error, which estimates the ratio of the tail probabilities, typically including the Cram\'er type moderate deviations.  \cite{shao1999cramer} established the following self-normalized Cram\'er type moderate deviation result for i.i.d. random variables. If $ \e | X_1 |^3 < \infty $, then
\begin{align*}
	\frac { \P ( S_n > x V_n  ) } {  1 - \Phi (x) }  \longrightarrow 1
\end{align*}
holds uniformly for $ 0 \leq x \leq o ( n^{ 1 / 6  } ) $, where $ \Phi (x) $ is the standard normal distribution function. \cite{jing2003self} further extended the result to general independent random variables. In particular, they obtained that if $ \e | X_i |^{ 3 }  < \infty$ for $ i \geq 1  $, then there exists an absolute constant $ A > 0  $, such that
\begin{equation}
	\bigg| \f{\P(S_n > x V_n  )}{1-\Phi(x)} - 1 \bigg| \leq  A (1+x)^3 \sum_{i=1}^{n}\e |X_i|^3/B_n^3  \label{classical}
\end{equation}
holds uniformly for $ 0 \leq x \leq  B_n (\sum_{i=1}^{n}\e |X_i|^3)^{ - 1/3}$.

In addition, \cite{wang2011refined} corrected the skewness in normal approximation and proved that if $ \e | X_i |^3  < \infty $ for $ i \geq 1 $, then there exist positive constants $ A_0 $ and $ C_0   $ such that
\begin{align}
	& \frac{\P\big(S_n>xV_n+c B_n \big)}{ [ 1-\Phi(x+c) ] \Psi_x } \label{refinedcase} \\
	&\hspace{1cm} =   e^{O_1 \Delta_{n,x}}  \bigg[ 1+ O_2 \Big( (1+x)\sum_{i=1}^{n}\e |X_i|^3 /B_n^3\Big) \bigg], \nn
\end{align}
holds uniformly for $ | c | \leq x / 5 $ and for all $  0 < x \leq  \f 1  3  B_n ( \max_i \e | X_i |^3 )^{ -1/3 } $ and $ x \leq C_0  B_n^3 / \sum_{i = 1}^n \e | X_i |^3 $, where $ | O_1 | \leq A_0 $, $ | O_2 | \leq A_0 $ and
\begin{align*}
	\Psi_x&=\exp\bigg[ \gamma^2\Big(\f{4\gamma}{3}-2\Big)x^3 \sum_{i=1}^{n}\e X_i^3 /B_n^3 \bigg],\\
	\Delta_{n,x}&= ( 1 + x )^3 B_n^{-3} \sum_{i=1}^{n}\e\big[  |X_i|^3 \mathbbm{1}( (1 + x) |X_i| \geq B_n) \big] \nn \\
	& \hspace{0.5cm} + (1 + x)^4 B_n^{-4} \sum_{i=1}^{n}\e\big[|X_i|^4 \mathbbm{1}( ( 1 + x ) |X_i|\leq B_n ) \big] ,
\end{align*}
with $\gamma = \frac 1 2  ( 1 + c  / x ) $. Especially, if $ X_1, \ldots, X_n $ are i.i.d. random variables with $ \e X_1 ^4  < \infty $, then \eqref{refinedcase} implies there exist positive constants $A_0$ and $C_0$ depending on $ \e X_1^2 $ and $ \e X_1^4 $ such that  
\begin{align*}
	& \frac { \P ( S_n > x V_n  ) } {  1 - \Phi ( x  )    } =  \exp \left  \{ -   \frac {  x^3 \e [ X_1 ^3 ] } {   3 \sqrt n ( \e [ X_1^2  ] )^{ 3/2 } }  \right\}  \left [   1 +  O_1  \left( \frac {   1 + x     } { \sqrt {n}  }     +  \frac {  ( 1 + x )^4    } { n }  \right) \right],   \label{W2011}
\end{align*}
uniformly in $ 0  < x \leq  C_0  n^{ 1 / 4 }  $, where $ | O_1 | \leq A_0 $.
Observe that in the i.i.d. case, the classical self-normalized Cram\'er type moderate deviation presented in \eqref{classical} gives a convergence rate of $ (1 + x )^3 / \sqrt{n} $ and the corresponding range of convergence $ x = o (n^{ 1/6 }) $. Thus, by specifying the skewness correction term $ \Psi_x $, \eqref{refinedcase} can improve the result of \cite{jing2003self} in terms of both the convergence rate and the range of convergence when the higher fourth moments exist.

It is worth mentioning that the moment conditions for self-normalized Cram\'er type moderate deviation theorems are much weaker than those in the classical theorems for standardized sums. As a result, to account for robustness against heavy-tailed data, the self-normalized sum would be recommended in real-world applications. We refer to \cite{pena2008self} for a systematic introduction to the theory and  statistical applications of self-normalized statistics.

Due to its rigorous control on the ratio of tail probabilities, the self-normalized Cram\'er type moderate deviation has been successfully applied in high-dimensional statistical analysis, including large-scale multiple testing \citep{fan2007many, liu2010cramer, liu2013cramer}, signal detection \citep{delaigle2009higher}, classification \citep{fan2008high} and feature screening \citep{chang2016local} among others.

Most of the existing works have focused on the classical self-normalized sum, that is, $ \sum_{i = 1}^n X_i / (\sum_{i = 1}^n X_i^2 )^{1/2} $ for independent random variables $\{X_i\}_{i = 1}^n $. Yet, in some scenarios, the sequence used for normalizing in the denominator could be different from the numerator, which occurs for a variety of commonly used studentized nonlinear statistics such as the studentized U-statistic and the studentized L-statistics. Therefore, investigations into general self-normalized processes beyond the classical form are imperative. \cite{shao2016cramer} attempted to extend the Cram\'er type moderate deviation theorem to a more general setting, that is, $ ( \sum_{i = 1}^n X_i + D_{1n} ) / ( ( \sum_{i = 1}^n X_i^2 ) (1 + D_{2n})  )^{1/2}$, where the remainders $ D_{1n} $ and $ D_{2n}$ are measurable functions of $ \{X_i\}_{i = 1}^n $ but negligible. Our present work will establish a fundamental framework in Theorem \ref{thm-main} on the Cram\'er type moderate deviation for a more general self-normalized form of $ \sum_{i = 1}^n X_i / (\sum_{i = 1}^n Y_i^2 )^{1/2} $, where $ \{(X_i, Y_i)\}_{i = 1}^n $ is a sequence of independent bivariate random vectors and $X_i$ and $ Y_i $ could be different from each other. It is worthwhile to mention that our Theorem \ref{thm-main} can cover not only the classical Cram\'er type moderate deviation for standardized sums by \cite{cramer1938nouveau}, but also the self-normalized counterparts by \cite{jing2003self} and \cite{wang2011refined}. 

Our investigation into the general self-normalized sum is also motivated by seeking to develop sharper self-normalized Cram\'er type moderate deviation results for weakly dependent random variables. Though Cram\'er type moderate deviation theory has been well studied for independent random variables, the theory for dependent data remains largely underdeveloped. The biggest challenge is that the classical theory for independent random variables cannot be directly applied due to dependence.  \cite{chen2016self} made the first attempt to develop the theory for self-normalized sums of dependent random variables with geometrically decaying dependence. However, their result can be further improved by applying our framework on the general self-normalized sum. 
% Specifically, to deal with dependence, it is common to separate the sum of weakly dependent random variables by big-block sums and small-block sums and then construct independent random variables of block sums by employing the big-block-small-block technique, which will be presented in more details in the proof of Theorems \ref{thm-2}--\ref{thm_GMC}.
The key observation is that after dividing the weakly dependent random variables into consecutive big blocks and small blocks, its self-normalized sum can be approximated by a general self-normalized sum of independent bivariate random vectors. Therefore, Cram\'er type moderate deviation theorems for the self-normalized sums of weakly dependent random variables can be established based on our fundamental theory on the general self-normalized sum. More details will be presented in Section \ref{Sec3}.

% However, with the degree of dependence of the underlying processes affecting the convergence rate and the range of convergence aside, it's also affected from the approach of using self-normalized sum to approximate pivotal quantities of dependent data. As a matter of fact, for pivotal quantity of interests in dependent data case, the normalizing term usually differs to the classical case with the covariance structure barged in.
%   Therefore, to establish Cram\'er type moderate deviation results for weakly dependent random variables, it is essential to investigate the general self-normalized sum.  
%
%it's usefulness beyond dependent data, such as the case in studying the Winsorized estimator. 

The rest of the paper is organized as follows. Our framework on general self-normalized Cram\'er type moderate deviation is presented in Section \ref{main.results.sec}. Section \ref{applications} shows applications to the self-normalized sums of weakly dependent random variables under one-dependence,  geometrically $\beta$-mixing condition, and geometric moment contraction. 
Section \ref{applicationsW} presents an additional application to studentized winsorized mean that naturally takes the form of a general self-normalized sum.
%construct simultaneously confidence interval for Winsorized estimator in ultra-high mean model. 
Section \ref{proof.sec} is devoted to proofs of the theorems in Sections \ref{main.results.sec}--\ref{applicationsW}. Other technical proof details  are included in the Supplementary Material.

%
%\vskip 3em
%
%Let $(X_1, Y_1), (X_2, Y_2), \ldots$, be independent random vectors satisfying $\e X_i=0$, and without loss of generality we assume that $\sum_{i=1}^n \e X_i^2=\sum_{i=1}^n \e {Y_i}^2=1$.
%Write $S_n = \sum_{i=1}^n X_i, V_n^2 = \sum_{i=1}^n Y_i^2 \ \text{and}\  T_n = S_n/V_n$.
%It is wellknown that $T_n$ converges to standard normal distribution, and we try to estimate the relative error of the convergence.
%There are some occasions when the statistics are not exactly the self-normalized form, so we generalize the denominator to some other random variables $Y_i$.
%
%\noi Actually for this general case, we obtain a result of Cram\'er type moderate deviation  under condition $ \e |X_i|^3 <\infty  $, $ \e |Y_i|^3 < \infty $ and $ \e e^{ X_i^2 \over {Y_i^2+c_0 \e Y_i^2 }} < \infty  $ for some constant $ c_0 > 0$.

\section{Main Results}%Moderate Deviation for General Self-normalized Sum}
\label{main.results.sec}
\setcounter{equation}{0}

Let $(X_1, Y_1), (X_2, Y_2), \ldots, (X_n, Y_n)$ be independent bivariate random vectors satisfying
\begin{equation}
	\e X_i=0 \;\; {\rm for}\;\; i \geq 1 \quad   {\rm and}\quad \sum_{i=1}^n \e X_i^2= 1 = \sum_{i=1}^n \e{Y_i}^2. \label{maincon1}
\end{equation}
We remark that for the convenience of presentation, $\{ X_i\} $ and $\{Y_i\}$ are standardized so $\sum_{i=1}^n
\e X_i^2 =1= \sum_{i=1}^n \e{Y_i}^2$. In other words, one should think of $X_i$ as $X_{n,i}$ and similarly $Y_i$ as $ Y_{n, i} $.
Let
\begin{equation}
	S_n = \sum_{i=1}^n X_i ,\quad V_n^2 = \sum_{i=1}^n Y_i^2 \quad {\rm and} \quad T_n=\f{S_n}{V_n}.
	\label{form}
\end{equation}
We first propose an exponential moment condition as follows.
Suppose there exists some constant $ c_0  \geq  0 $ such that for $x > 0 $ satisfying \eqref{G2},
\begin{equation}
	\e e^{ \min\left\{\frac{X_i^2}{ Y_i^2+c_0 \e Y_i^2 }, \, 2 x X_i \right\} } < \infty.  \label{maincon2}
\end{equation}

The above moment condition links $ X_i $ with $ Y_i $ and shows how they interact with each other. In particular, it is automatically satisfied for the classical self-normalized sum with $Y_i=X_i$.

The following notations  will be used throughout the paper. Define
\begin{align}
	L_{3,n} = \sum_{i=1}^n \big( \e |X_i|^3 + \e |Y_i|^3 \big),  \qquad \qquad \qquad \qquad \qquad \qquad  \label{Def}
\end{align}
\vspace{-1em}
\begin{align*}
	\de_{x,i} &= (1 + x )^3 \Big( \e \big[ |X_i  |^3 \mathbbm{1} (  | ( 1 + x ) X_i|> 1 )\big] + \e \big[  |Y_i  |^3 \mathbbm{1} ( | ( 1 + x  ) Y_i|>1 ) \big] \Big) \nn \\
	&~~ + ( 1 + x )^4 \Big( \e \big[ |X_i |^4 \mathbbm{1} (  | ( 1 + x  ) X_i| \leq  1 )\big] + \e \big[  |Y_i  |^4 \mathbbm{1} ( | ( 1 + x  ) Y_i| \leq 1 ) \big] \Big), \nn \\
	r_{x,i} &= \e\left[ \exp \Big\{  \min\Big( \frac{X_i^2}{ Y_i^2+c_0 \e Y_i^2 }, 2xX_i \Big) \Big\} \mathbbm{1} ( | ( 1 + x ) X_i|>1 )\right],  \nn \\
	R_{x,i} & = \de_{x,i} + r_{x,i},  \quad \de_x =\sum_{i=1}^n \de_{x,i},  \quad  r_x=\sum_{i=1}^n r_{x,i}, \quad \mbox{and} \quad  R_x = \delta_x + r_x. \nn
\end{align*}

\begin{theorem} \label{thm-main}
	Assume \eqref{maincon1} and \eqref{maincon2} are satisfied. In addition, $ \e | X_i |^3 < \infty $ and $ \e | Y_i |^3 < \infty $ for $ i \geq 1 $. Then there exist absolute positive constants $ 0 < c_1\leq 1/4 $ and $ A > 0  $ such that
	\begin{equation}
		\P  (S_n \ge x V_n+c  )  =  [ 1 - \Phi (x+c) ] \, \Psi_x^* \, e^{ O_1 R_x } ( 1+O_2 (1+x) L_{3,n} ),
		\label{main-re}
	\end{equation}
	where
	\begin{equation*}
		\Psi_x^*=\exp \Big\{  x^3 \Big(  \f43 \ga^3 \sum_{i=1}^n \e X_i^3 - 2 \ga^2 \sum_{i=1}^n \e [X_i Y_i^2] \Big) \Big\} \quad {\rm and}\quad \ga = \f 12 ( 1 + \f cx ),
	\end{equation*}
	uniformly for $ |c| \le x/5 $ and for all $ x > 0 $ satisfying
	\begin{gather}
		(1 + x ) L_{3,n} \le c_1 , \qquad x^{ -2  } R_x  \le c_1 ,  \label{G1}\\
		\mb{and} \qquad x \le {\f{1}{4} \land \f{1}{2\sqrt{c_0}}\over \big[  \max_i (  \e|X_i|^3 +\e|Y_i|^3 ) \big]^{1/3}}, \quad \max_i r_{x,i} \le c_1 , \label{G2}
	\end{gather}
	\label{thm-1}	
	where $|O_1| \leq A $ and $ | O_2  |  \leq A $.
\end{theorem}

Theorem \ref{thm-main} unifies the classical standard and self-normalized Cram\'er type moderate deviation theorems as well as the refined version by \cite{wang2011refined}. The assumption \eqref{maincon2}
%%is a crurial condition and causes the primary difficulty in applications of Theorem \ref{thm-main}. However, it
is satisfied for a wide class of statistics, including the block sums of weakly dependent random variables and the self-normalized winsorized mean. More details will be provided in the proof of Theorems \ref{thm-2}--\ref{thm_GMC} and \ref{thm-winsor-2}.

%
%\ignore{
%And for an even more special case where $X_1,X_2,\cdots,X_n$ are i.i.d random variables, the error term $O_1R_x=O_1\delta_x$ in (\ref{main-re}) would be dominated by $xL_{3,n}$ when x is in a relatively small range, for instance $0\leq x\leq O(n^{1/6})$. However $O_1\delta_x$ plays a part when $x$ is in a median range, for instance $ O(n^{1/6})\leq x \leq o(n^{1/2})$. It remains interesting whether it's possible to remove the error term $O_1\delta_x$ for standard self-normalized sum.
%}

The following corollary is a straightforward application of Theorem \ref{thm-main} to the classical standardized sum of independent random variables. The proof will be given in Section \ref{pf-cor} in the Supplementary Material. 
\begin{corollary} \label{cor}
	Let $ X_1, \ldots, X_n$ be i.i.d. random variables with $ \e X_i = 0 $ and $   \e [ X_i^2  ]  = \sigma^2 $. Denote $ S_n = \sum_{i = 1}^n X_i $. If there exists a positive constant $t_0  $ such that $ \e e^{ t X_1 } < \infty $ for $ | t  | \leq t_0 $, then we have
	\begin{align}
		\frac { \P ( S_n > x  \sigma \sqrt n  ) } { 1 - \Phi (x) }  = \exp \Big\{ \frac { x^3   \e X_1^3 } { 6 \sigma^3 \sqrt{n}  } \Big \} \Big[   1 + O \Big(\frac { (1 + x)^4   } {  n  } + \frac   { 1 + x  } { \sqrt n  }\Big) \Big] \label{classical-sum}
	\end{align}
	holds uniformly for $ 0 < x \leq   O ( n^{1/4})  $.
\end{corollary}

\section{Applications to Weakly Dependent Random Variables} \label{Sec3}

\label{applications}
\setcounter{equation}{0}

The Cram\'er type moderate deviation theory has been well studied for independent random variables, yet there are few results available for dependent data. The novel work of \cite{chen2016self} made the first attempt to develop the theory for self-normalized sums of weakly dependent random variables satisfying the geometrically $\beta$-mixing condition or geometric moment contraction. In this section, we will further improve their results by applying our fundamental framework on the general self-normalized sum. Before that, we will start with a Cram\'er type moderate deviation theorem for self-normalized sums of one-dependent random variables. The reason to first investigate under one-dependence is two-fold. First, one-dependence is the simplest scenario of dependency and the result for one-dependent random variables can be applied to  $m$-dependence, where $m$ could also go to infinity. Second, many weakly dependent random variables can be approximated by some one-dependent random variables. A typical example includes the block sums of random variables satisfying geometric moment contraction, which will be presented in Section \ref{Sec:GMC}. Therefore, the following Theorem \ref{thm-2} under one-dependence lays the foundation for Theorem \ref{thm_GMC} under geometric moment contraction.

% In this section, we establish Cram\'er type moderate deviation theorems for self-normalized sums of one-dependent random variables and also improve the results of \cite{chen2016self} for $\beta$-mixing and geometrically moment contracting random variables by applying our main result, Theorem \ref{thm-main}.

% The key point is to construct independent structure based on the original weakly dependent random variables.

\subsection{Cram\'er type moderate deviation under one-dependence}
Let $\xi_1, \xi_2, \ldots$ be one-dependent random variables, which means for $ i, j \geq 1 $, $ \xi_i $ is independent of $ \xi_j $ if $ | j - i | \geq 2$. Put
\begin{equation}
	S_n = \sum_{i=1} ^{n} \xi_i, \quad  V_n^2 = \sum_{i=1}^n \xi_i^2, \quad  \rho_n = \frac{\sum_{i=1}^{n-1} \e \xi_i\xi_{i+1} }{ \sum_{i=1}^{n} \e \xi_i^2  }. \label{notationonedependentbasic}
\end{equation}

Note that $  | \rho_n | \leq 1/2 $. Moreover, in many applications where weakly dependent sequence can be approximated by some one-dependent random variables, the covariances $ \e \xi_i \xi_{i + 1} $ are negligible compared to the variables $ \e \xi_i^2$ due to weak dependence, hence  $\rho_n \to 0 $ as $n \to \infty $. Therefore, $ \rho_n $ can be moved to the remainders and the limiting distribution will still be standard normal. One can find more details in the proof of Theorem \ref{thm_GMC}, which obtains Cram\'er type moderate deviation result under geometric moment contraction by applying Theorem \ref{thm-2}.

Under existence of the fourth moment, we have the following theorem for self-normalized sums of one-dependent random variables.
\begin{theorem} 	\label{thm-2}
	Assume that $ \e \xi_i = 0 $, $ \e  \xi_i ^4 \le a_1^4 $ and $ \e  \xi_i^2 \ge a_2^2 $ for $1 \leq i \leq n$ and $    \rho_n \geq \rho $ for some $ \rho > - 1/2  $. Denote $ a= a_1/a_2 $. Then there exist positive numbers  $a_0$ and $A(\rho)$ depending on $ \rho $ such that
	\begin{equation}
		\P (S_n > x V_n) = \Big[ 1 - \Phi \big( \f{x}{\sqrt{1+2\rho_n}} \big) \Big] \Big( 1 + O_1 a^4 \f{(1+x)^2}{n^{1/4}} \Big) 
		\label{re-dep}
	\end{equation}
	holds uniformly for $ x \in ( 0, a_0 a^{-2}  n^{\f18} ) $, where $|O_1|\leq A(\rho)$.
\end{theorem}

The proof of Theorem \ref{thm-2} relies on the big-block-small-block technique and an application of Theorem \ref{thm-1}. The main idea is to approximate the self-normalized sum of one-dependent random variables by a general self-normalized sum of independent random vectors based on the big blocks. In more details, let the length of big blocks be $ l = [n^{\alpha}] $ for $ 0 < \alpha < 1 $, where $ [x] $  denotes the integer part of $ x $ for any $x > 0$, and the length of small blocks be only 1. Denote $  k = [n / (l + 1)] $. For $ 1 \leq j \leq k $, we define the $j$-th big block by
\begin{equation}
	H_j = \{i: (j - 1) (l + 1) + 1 \leq i \leq j (l +1) -1 \} 
\end{equation}
and the sums over $j$-th big block by
\begin{equation}
	X_j = \sum_{i \in H_j} \xi_i ~~\mbox{and}~~ Y_j^2 = \sum_{ i \in H_j } \xi_i^2.
\end{equation}
Observe that by this construction, the big-block sums $ \{ X_j  \}_{j = 1}^{ k} $ and $ \{ Y_j  \}_{j = 1 }^k $ are both sequences of independent random variables, because $ \{\xi_i \}_{  i = 1 }^n  $ are one-dependent and the adjacent big blocks $ H_j $ and $H_{j+1 } $ are separated by a random variable $ \xi_{j (l + 1)}$. Since the big blocks contain $ (1 - n^{ - \alpha}) $ proportion of the random variables in $ \{\xi_i \}_{i = 1}^n   $, we can approximate the self-normalized sum $ S_n / V_n  $ of $ \{ \xi_i  \}_{ i = 1}^n $ by the general self-normalized sum $ \sum_{j = 1}^k X_j / ( \sum_{j = 1}^k Y_j^2 )^{1/2} $.  The crucial quantity $r_{x,j}$ can be separated into two self-normalized sums of independent random variables due to one-dependence and thus can be bounded by using Lemma A.5 in the Supplementary Material.
Therefore, Theorem \ref{thm-2} can be proved by applying Theorem \ref{thm-1} and calculating the error terms involved. More details of the proof will be provided in Section \ref{pf-thm3.1}.

Compared with the classical result \eqref{classical} for independent data, one-dependence results in a narrower zone of convergence and a slower convergence rate. Moreover, \eqref{re-dep} can be easily extended to general $m$-dependent random variables, where $m$ could depend on $n$ and go to infinity. Indeed, if $ Z_1, \ldots, Z_n  $ are $ m  $-dependent and suppose $ b =  n / m $ is an integer for simplicity, we define $ \xi_j = \sum_{i = 1 + ( j - 1 ) m }^{ j m  }  Z_i $ for $ 1 \leq j \leq b $, then $ \xi_1, \ldots, \xi_b $ are one-dependent random variables and Theorem \ref{thm-2} can be applied.

\subsection{Cram\'er type moderate deviation under $ \beta $-mixing} \label{Sec:beta}
In time series, asymptotic independence conditions such as mixing conditions are usually proposed to replace independence, among which $ \beta $-mixing is an important dependent structure and has been connected with a large class of time series models including ARMA models, GARCH models and certain Markov processes. This subsection provides a Cram\'er type moderate deviation theorem for block-normalized sums of geometrically $\beta$-mixing random variables, which improves the result by \cite{chen2016self}.

Let $ \{  X_i  \}_{i = 1}^n $ be a sequence of random variables. Let $ \sigma_{ - \infty }^{t }    $ and $ \sigma_{ t + m }^{ \infty } $ be $ \sigma $-fields generated by $ \{ X_i \}_{ 1 \leq  i \leq t} $ and $ \{ X_i \}_{ i \geq t + m } $, respectively. The $\beta$-mixing coefficient is given by
\begin{align}
	\beta( m ) := \sup \limits_{t} \e \sup \{  | \P ( B | \sigma_{ - \infty }^t  )  - \P (B) | : B \in \sigma_{ t + m   }^{ \infty }      \}.   \label{beta-coeff}
\end{align}
We say $ \{  X_i  \}_{ i \geq 1}  $ is geometrically $\beta$-mixing if $\beta (m) $ admits an exponentially decaying rate, that is, there exist positive numbers $ a_1 $, $ a_2 $ and $ \tau $ such that
\begin{equation}
	\beta ( m ) \leq a_1 e^{ - a_2 m^{ \tau } }.     \label{g-beta}
\end{equation}

To account for dependence, the block technique is naturally used to estimate the variance of sums of dependent random variables (see \cite{chen2016self}). Set $ l = [ n^{\alpha} ] +1 $ for $ 0 < \alpha < 1  $ and $ k = [ n / l ] $ . For $ 1 \leq j \leq k  $, define the $ j $-th block and corresponding $ j $-th block sum by
\begin{align}
	H_j = \{ i: l (j-1) + 1 \leq i \le  l j  \}   \quad  \mbox{and} \quad  Y_j = \sum_{ i \in H_j } X_i, \label{blocksum}
\end{align}
respectively. The block-normalized sum is then defined by
\begin{equation}
	T_k  = \frac {  \sum_{ j = 1 }^k Y_j  }   { \sqrt{ \sum_{ j = 1 }^k Y_j^2 } }.  \label{block-selfsum}
\end{equation}

\begin{theorem} 	\label{thm-blocksum}
	Let $ \{ X_i\}_{i = 1}^n $ be a $\beta$-mixing sequence satisfying \eqref{g-beta}. Assume $ \e  X_i  = 0 $ and there exist positive numbers $ \mu_1  $ and $ \mu_2 $ such that $ \e |X_i|^{r} \le \mu_1^{r} $ for $r > 4$ and $ \e (\sum_{i= s }^{  s + t } X_i )^2 \ge \mu_2^2 t $ for all $i\geq 1, s \geq 0 , t \geq 1 $. Then, for $ 0 < \alpha < 1 $ and $ \tau > 0 $, there exist positive numbers $ A $ and $ d_0 $ depending on $a_1, a_2, \mu_1, \mu_2, \alpha $, $ \tau $ and $r$ such that,
	\begin{equation}
		\Big | \frac {  \P ( T_k \geq x  ) }  {   1- \Phi(x)  } - 1  \Big | \leq   A  \Big(    \f{ ( 1 + x) ^2  }{ n^{\alpha }  } + \f{ ( 1 + x )^2  \log n } {n^{  \min\{(1 - \alpha)/4, \, \alpha \tau/2 \} }}  \Big)  \label{blocksum1}
	\end{equation}
	uniformly in $ 0 \leq x \leq  d_0 \min \{   n^{ \alpha / 2  }, ( \log n  )^{ - 1/2  } n^{  \min\{  ( 1 - \alpha ) / 8,\, \alpha \tau / 4 \} }  \}  $.  
\end{theorem}

%\begin{remark}
As for the choice of $ \alpha $, for any given $ \tau $, we can always choose $ \alpha $ such that $  1 - \alpha  \le  2 \alpha \tau $, that is, $ \alpha \ge  \f{1}{1+2\tau} $. To optimize the convergence rate and the range of $x$ in \eqref{blocksum1},

\noi  (i) when $ \tau \ge 2  $, let $ \alpha = \f15 $, then 
\begin{equation}
	\f{ \P ( T_k \ge x ) }{ 1 - \Phi(x) } = 1 + O \Big( \f{( 1 + x )^2 \log n }{ n^{1/5} }  \Big)  \label{sce1}
\end{equation}
uniformly for $  x \in (0, d_1 (\log n)^{-1/2} n^{1/10} )$;

\noi  (ii) when $ \tau < 2  $, let $ \alpha =\f{1}{1+2\tau} $, then
\begin{equation}
	\f{ \P ( T_k \ge x ) }{ 1 - \Phi(x) } = 1 + O \Big( \f{( 1 + x )^2 \log n }{ n^{\f{\tau}{2(1+2\tau)} }} \Big) \label{sce2}
\end{equation}
uniformly for $  x \in (0, d_1  (\log n)^{-1/2} n^{\f{\tau}{4(1+2\tau)} } )$.

%\end{remark}

\begin{remark}
	We now compare our result with Theorem 4.2 in \cite{chen2016self}. They proved that given $ \e |X_i|^{r} < \infty $ for $r > 3$ and the same assumption \eqref{g-beta},
	\begin{align}
		\Big|\f{ \P ( T_k \ge x ) }{ 1 - \Phi(x) } - 1 \Big| \leq A \Big[ \frac { ( 1  + x)^2  }{ n^{\alpha} } + \frac { (1 + x)^{5/4} } {n^{ (1 - \alpha ) /8  }   } \Big] \label{remark}
	\end{align}
	uniformly in $  0 \leq x \leq  d_0 ( \min \{ ( \log n )^{ - 4/5 } n^{( 1 - \alpha )/10 }, n^{ \alpha \tau  /2 }, n^{\alpha /2 }\} )  $. (We have to mention that the original version of Theorem 4.2 in \cite{chen2016self} missed the error term $ \frac { (1 + x)^2 } { n^{\alpha} } $ and the corresponding condition $ x \leq d_0 n^{\alpha / 2 } $.)
	When $ (1 - \alpha) > 2 \alpha \tau  $, our results might be worse for some choices of $\alpha$. However,  when $ (1 - \alpha ) \leq 2 \alpha \tau $, our results improve theirs in terms of both the convergence rate and the corresponding range of $x$.  
	The improvement is achieved by applying our framework for general self-normalized sum and correcting the bias to normal approximation by specifying the skewness term $ \Psi_x^*  $.
	
	%	The convergence rate and range of convergence obtained above for (\ref{blocksum1}) and (\ref{blocksum2}), are significantly better compare to the result established by \cite{chen2016self} which shares the same condition. The improvement relies heavily on our main result, Theorem \ref{thm-1}, for general self-normalized sums, in which we specified the bias term $ \Psi_x^*  $ and be able to estimate the tail probability more delicately.
\end{remark}

%\begin{remark}
%	If $\beta$-mixing is strengthen to $\phi$-mixing, then we are able to strengthen the result by get rid of the condition $ \f{1-\alpha}{2} \le \alpha \tau $.	
%\end{remark}

The proof of Theorem \ref{thm-blocksum} again builds on the big-block-small-block technique. Recall that for $ 1 \leq j \leq k $, $ Y_j = \sum_{i \in H_j} X_i $ is a block sum defined in \eqref{blocksum}. We first apply big-block-small-block technique to separate the sequence $ \{ Y_j  \}_{j = 1}^k  $ into consecutive big blocks and small blocks. Let the size of big-blocks be $ m_1 = [ n^{\alpha_1} ] $ for $ 0 < \alpha_1 < 1 - \alpha $ and the size of small-blocks be only $1$. Denote $ k_1 = [ k / ( m_1 +1 ) ]   $. For $ 1 \leq u \leq k_1 $, define the $u$-th big block by
\begin{align}
	I_u = \{ j: (m_1 + 1 ) ( u - 1 ) + 1 \leq j \leq (m_1 + 1 ) u - 1  \},
\end{align}
and the sums over $u$-th big block by
\begin{align}
	\zeta_u = \sum_{j \in I_u} Y_j, \quad \eta_u^2 = \sum_{j \in I_u} Y_j^2.
\end{align}
Then the self-normalized sum $ T_k = \sum_{j = 1}^k Y_j / (\sum_{j = 1}^k Y_j^2) ^{1/2}$ can be approximated by the general self-normalized sum $ \sum_{u = 1}^{k_1} \zeta_u / (\sum_{u = 1}^{k_1} \eta_u^2)^{1/2} $ constructed on the big blocks. However, unlike the big-block sums under one-dependence, the big-block sums $ \{ (\zeta_u, \eta_u)   \}_{u = 1}^{k_1} $ are not independent under $ \beta $-mixing assumption in \eqref{g-beta}, but weakly dependent. Therefore, Theorem \ref{thm-1} can not be directly applied. Note that the adjacent random vectors $ ( \zeta_u, \eta_u ) $ and $ (\zeta_{u+1}, \eta_{u+1} )$ depend on $  \{Y_j\}_{  j = (m_1 + 1 ) (u - 1)+ 1 }^{  (m_1 + 1 ) u - 1 } $ and $  \{Y_j \}_{  (m_1 + 1 ) u + 1  }^{(m_1 + 1 ) (u + 1) - 1 }$, respectively. Since $ Y_j $ defined in \eqref{blocksum} is a block sum of $ \{X_i\} $ and by the $ \beta$-mixing assumption on $ \{X_i\} $, we can see the $\beta$-mixing dependence coefficient between $ ( \zeta_u, \eta_u )  $ and $  (\zeta_{u+1}, \eta_{u+1} ) $ is bounded by $ O(  e^{ - a_2 n^{ \alpha \tau} } ) $, which converges to 0 as $ n \to \infty $. According to Lemma \ref{lemma-indep} \citep{berbee1987convergence} presented in Section \ref{pf-thm3.2}, the weakly dependent random vectors $ \{ (\zeta_u, \eta_u)\}_{u = 1}^{k_1} $ can be replaced with independent random vectors $ \{ ( \tilde{\zeta}_u, \tilde{\eta}_u ) \}_{u = 1}^{k_1} $ that have the same marginal distributions, with probability $ 1 - O (k_1 e^{ - a_2 n^{\alpha \tau}}) $. Moreover, the crucial quantity $r_{x,u}$ can be approximated by two self-normalized sums of independent random variables due to the $\beta$-mixing assumption (see Lemma 5.3) and the block technique. Consequently, our main result in Theorem \ref{thm-1} can be applied to the general self-normalized sum $ \sum_{u = 1}^{k_1} \tilde{\zeta}_u / (\sum_{u = 1}^{k_1} \tilde{\eta}_u^2)^{1/2} $. Detailed proof will be given in Section \ref{pf-thm3.2}.

\subsection{Cram\'er type moderate deviation for causal processes under geometric moment contraction (GMC)} \label{Sec:GMC}
The GMC (see \cite{wu2004limit}, \cite{hsing2004weighted} and \cite{wu2005nonlinear, wu2011asymptotic}) is satisfied by many non-linear time series models including various GARCH models that are commonly used in statistics, econometrics and engineering. In this subsection, we present a Cram\'er type moderate deviation theorem for block normalized sums of random variables satisfying GMC.

Let $ \{ \varepsilon_t\}_{ t \in \mathbbm{Z}  } $ be i.i.d. random variables and define $ \sigma $-fields $  \mathscr{F}_t = \sigma ( \ldots, \varepsilon_{t-1}, \varepsilon_t )  $. Suppose that $ \{ X_i = G_i ( \mathscr{F}_i ) \}_{i \geq 1} $ is a causal process with $ G_i ( \cdot )  $ being a measurable function such that $ X_i $ is well-defined. Let $ \{ \varepsilon_{t}^*\}_{ t \in \mathbbm{Z} } $ be an independent copy of $ \{ \varepsilon_{t} \}_{ t \in \mathbbm{Z} } $ and we similarly define $  \mathscr{F}_t^* = \sigma ( \ldots, \varepsilon_{t-1}^*, \varepsilon_t^* )   $.

\begin{definition}
	{ \it (GMC).} Assume that $ \e  | X_i |^r  < \infty $ for all $ i  \geq 1  $ with $ r > 2 $. Define the functional dependence measure by
	\begin{equation}
		\Delta_r ( n  ) =  \sup \limits_i  \Vert X_i - G_i ( \mathscr{F}_{i - n}^*, \varepsilon_{ i - n + 1 }, \ldots, \varepsilon_{ i }  ) \Vert_r, \label{Delta_n}
	\end{equation}
	where $ \Vert \cdot  \Vert_r  = (  \e  | \cdot |^r   )^{ 1 / r } $. We say $ \{ X_i  \}_{i \geq 1} $ satisfies GMC if there exist positive constants $ a_1  $, $ a_2   $ and $ 0<\tau\leq 1  $ such that
	\begin{equation}
		\Delta_r ( n  )  \leq   a_1 e^{ - a_2 n^{\tau} } .    \label{GMC}
	\end{equation}
\end{definition}

Note that the GMC property \eqref{GMC} implies $ \{ X_i  \}_{  i \geq 1 } $ forgets the past $ \mathscr{F}_{0} = \sigma  ( \ldots, \vp_{-1}, \vp_0 ) $ geometrically fast.

\begin{remark}   \label{remark-theta}
	Define another functional dependence measure as
	\begin{equation}
		\theta_r ( n ) =  \sup \limits_i  \Vert  X_ i -  G_i ( \ldots, \veps_{i - n - 2}, \veps_{i - n - 1}, \veps_{i - n }^*, \veps_{i - n + 1} , \ldots, \veps_{ i } )  \Vert_r.   \label{theta_r}
	\end{equation}
	The property \eqref{GMC} is equivalent to $ \theta_r (n) \leq    a_1' e ^{ - a_2' n^{\tau} }   $ for some positive constants $ a_1' $ and $ a_2' $.
	\ignore{  Indeed,
		\begin{align*}
			\Delta_r ( n ) & \leq \sup_i \Vert X_i - G_t(\ldots,\epsilon_{i-n-1},\varepsilon_{i-n}^*,\varepsilon_{i-n+1},\ldots,\varepsilon_i)\parallel_r \nn\\
			& \quad +\sup_i  \Vert G_i (\ldots,\varepsilon_{i-n-1},\varepsilon_{i-n}^*,\varepsilon_{i-n+1},\ldots,\epsilon_i )\nn \\
			&\qquad \qquad \; -G_i (\mathscr{F}_{i-n}^*,\varepsilon_{i-n+1},\ldots,\varepsilon_i )\parallel_r \nn\\
			& = \theta_r ( n ) + \sup_i \parallel X_i - G_i (\mathscr{F}_{i-n-1}^*,\epsilon_{i-n},\ldots,\epsilon_i )\parallel_r \leq \sum_{ m = n }^{\infty} \theta_r (m),
		\end{align*}
		and in the same manner,
		\begin{equation*}
			\theta_r (n) \leq \Delta_r ( n ) + \Delta_r ( n + 1  ).
		\end{equation*}
		%	together with the fact
		%	\begin{align*}
		%	\sum_{k=m}^{+\infty}a_1e^{-a_2k^{\tau}}&=\sum_{k=m}^{2m-1}a_1e^{-a_2k^{\tau}}+\sum_{k=2m}^{+\infty}a_1e^{-a_2k^{\tau}}\\
		%	&\leq a_1me^{-a_2m^{\tau}}+a_1e^{-a_2m^{\tau}}\sum_{k=2m}^{+\infty}e^{-a_2(1-2^{-\tau})k^{\tau}} \\
		%	&\leq a_1e^{-a_2m^{\tau}/2}+a_3e^{-a_2m^{\tau}} \leq a_4e^{-a_5m^{\tau}}
		%	\end{align*}
		%	
		Observing that $ \int_{n }^{\infty} e^{ - a_2 x^{ \tau } }d x  \leq A e^{ - a_2 n^{ \tau } /2  }   $ for some constant $ A $, it is easy to see \eqref{GMC} is equivalent to $ \theta_r (n) \leq a_1' e^{ - a_2' n^{\tau} }   $.
	}
\end{remark}

Now we assume $\{ X_i \}_{i = 1}^n $ is a sequence of random variables with
\begin{equation}
	\e X_i =0\;,\quad \quad \e|X_i|^4 < \infty\;, \label{GMCc1}
\end{equation}
for all $i\geq 1$. Write $S_{k,m}=\sum_{i=k+1}^{k+m}X_i$. Assume there exists a positive number $ \omega_1$ such that for any $ k \geq 0, m \geq 1 $,
\begin{equation}
	\e(S_{k,m}^2)\geq  \omega_1^2 m.    \label{GMCc2}
\end{equation}

As with the procedure for $\beta$-mixing random variables, we construct the block-normalized sum for random variables satisfying GMC. Let the block size $ m= [ n ^{\alpha}  ] $ for $0<\alpha<1$ and $k= [ n / m ] $. For $ 1 \leq j \leq k  $, define the $ j $-th block and the $ j $-th block sum by
\begin{equation*}
	H_j= \{ i : m(j-1) + 1 \leq i \leq m j \} \quad \mb{and} \quad Y_j = \sum_{ i \in H_j } X_i.
\end{equation*}
The block-normalized sum is then given by
\begin{equation}
	T_k = \f { \sum_{ j = 1 }^k Y_j } { \sqrt{ \sum_{ j =1 }^k Y_j^2  } }.
\end{equation}

\begin{theorem}   	\label{thm_GMC}
	Assume $ \{ X_i\}_{  i = 1 }^n $ is a causal process satisfying \eqref{GMC},  \eqref{GMCc1} and \eqref{GMCc2}. Then we have for $ 0 < \alpha < 1 $ and $ \tau > 0   $, there exist positive numbers $ A $ and $ d_0 $ depending on $ a_1, a_2, \omega_1 , \alpha $ and $ \tau $ such that
	\begin{align}
		\P ( T_k \geq  x  ) =   [ 1 - \Phi (x) ] \Big( 1 + O_1 ( \frac { 1 + x^2 } { n^{ \alpha } } + \frac { 1 + x^2  } { n^{ (1 - \alpha ) / 4  } } ) \Big) \label{GMC_result}
	\end{align}	
	uniformly in $ 0 \leq x \leq d_0 \min \{    n^{ \alpha \tau / 2  } ,n^{ \alpha / 2 },  n^{ ( 1 - \alpha  )/ 8 } \} $, where $ | O_1 |  \leq A $.
\end{theorem}

Corollary 4.3 in \cite{chen2016self} stated that \eqref{remark} also holds for the self-normalized block sum of GMC random variables. Compared with their result, our convergence rate and the associated converging range of $x$ significantly improve theirs. 

The main idea of our proof for Theorem \ref{thm_GMC} is to approximate $ \{Y_j \}_{j = 1 }^k $ by one-dependent random variables and then apply Theorem \ref{thm-2}. Define 
\begin{equation}
	\tilde{Y}_j  = \e ( Y_j | \varepsilon_{l}, m(j - 2) + 1 \leq l \leq m j )
\end{equation}
and 
\begin{equation}
	\tilde{T}_k = \frac {  \sum_{j = 1}^k \tilde{Y}_j } { ( \sum_{j = 1}^k \tilde{Y}_j^2 )^{1/2} },
\end{equation}
where $ m = [n^{\alpha}] $. Since $ Y_j  = \sum_{i \in H_j } X_i$ belongs to $ \mathscr{F}_{mj} $ and depends weakly on $ \mathscr{F}_{m (j - 2)} $ by the GMC assumption \eqref{GMC}, it is intuitive that $ \tilde{Y}_j $ is close to $Y_j$. In particular, we can prove $ \Vert \tilde{Y}_j - Y_j \Vert_r  \leq  a_1 m e^{ - a_2 m \tau} $. Therefore, the self-normalized sum $ T_k $ of $ \{Y_j\}_{j = 1}^k  $ can be well approximated by the self-normalized sum $ \tilde{T}_k $ of $ \{\tilde{Y}_j \}_{j = 1}^k $. Moreover, since $ \{\varepsilon_t \}_{t \in \mathbbm{Z}} $ are i.i.d. random variables, it is easy to see $ \{\tilde{Y}_j \}_{j = 1}^k $ are one-dependent. Consequently,  Theorem \ref{thm_GMC} can be proved by applying Theorem \ref{thm-2} and controlling the errors caused by the approximation by one-dependent random variables. We will present the detailed proof in Section \ref{pf-thm3.3}.

\section{Applications to self-normalized winsorized mean}
\label{applicationsW}

Although the sample mean has always been a prominent unbiased estimator for  a location parameter, it has the troubling disadvantage of being heavily influenced by gross outliers. Yet, robustness is often a desirable property, especially in real-world applications. Thus robust alternatives, typically including the trimmed mean (\cite{rothenberg1964note}), the winsorized mean (\cite{dixon1960simplified}, \cite{huber1992Robust}), and the Huber estimator (\cite{huber1992Robust, huber1973Robust}), are imperative to make more reliable statistical inference for unknown parameters. Suppose we have i.i.d. observations $ Y_1, Y_{2}, \ldots, Y_n $ with common distribution $Y$ and 
$$
\mu =  \e [Y] \quad \mbox{and} \quad \sigma^2 = \var (Y).
$$
For a thresholding parameter $ \tau > 0 $ that determines the tradeoff between bias and robustness, the winsorized mean is defined by 
\begin{equation} \label{winsor}
	\hat{\mu}_{W} = n^{-1}  \sum_{i = 1}^n  f  (Y_i),   
\end{equation}
where 
\begin{equation} \label{f}
	f  (x) =   x \mathbbm{1} ( | x | \leq \tau  )  + \tau \mathbbm{1} (  x > \tau ) - \tau \mathbbm{1} ( x < - \tau ). 
\end{equation}
The trimmed mean is defined by
\begin{equation}
	\hat{\mu}_{T} = c_{\tau, n}^{-1}  \sum_{i = 1}^n  Y_i \mathbbm{1} \{ | Y_i | \leq \tau \},    \label{trimmed}
\end{equation}
where $ c_{\tau, n} = \sum_{i=1}^n  \mathbbm{1} \{ | Y_i | \leq \tau \} $. 
%We suppress the dependence of $f (x)$ on $ \tau $ to simplify the notation. 
Moreover, the Huber loss (\cite{huber1992Robust}) is given by 
\begin{equation}
	\ell_{\tau} (u) = \left\{
	\begin{aligned}
		& \frac 12 u^2   \quad & \mathrm{if}~| u | \leq \tau, \\
		& \tau | u | - \frac 12 \tau^2 \quad  & \mathrm{if}~ | u | > \tau,
	\end{aligned}
	\right.
\end{equation}
which is a compromise between square loss and absolute loss. The Huber estimator is then defined as 
\begin{equation}
	\hat{\mu}_{H} = \mathop{\argmin} \limits_{\mu \in \mathbbm{R}} \sum_{i = 1}^n \ell_{\tau} (Y_i - \mu). \label{huber}
\end{equation}

These robust estimators are common in reducing the impact of outliers and are all asymptotically equivalent to the sample mean when the associated tuning parameter $\tau$ tends to infinity. Compared with the Huber estimator, the trimmed mean and the winsorized mean have explicit formulas and therefore are easier to be applied in real-world applications. 
It is well-known that these robust estimators are asymptotically normal under some regularity conditions. 
%({\color{red} find reference, such as Bickel 1965, and others}) 
Recently \cite{Z2018} obtained a Cram\'er-type moderate deviation theorem for the Huber estimator when allowing the tuning parameter $\tau$ to diverge with the sample size $n$ in some regime, and they applied the result to establish theoretical guarantees for the false discovery rate in multiple testing procedure for population means. 
However, the statistic they investigated depends on the unknown variance, which needs to be well estimated in practice.  

\subsection{Cram\'er-type moderate deviation for self-normalized winsorized mean and trimmed mean}
In this section, we will provide Cram\'er-type moderate deviation theorems for the self-normalized winsorized mean defined in \eqref{S} and self-normalized trimmed mean defined in \eqref{self-trimmed}, as an application of our main Theorem \ref{thm-1}. 
{The self-normalized winsorized mean and trimmed mean} are asymptotically pivotal statistics in the sense that their asymptotic distributions do not depend on unknown parameters as $(n, \tau ) \to (\infty, \infty)$, therefore they can be directly used in the multiple testing of population means with theoretical justification.
In addition, we will see our results for self-normalized winsorized mean and trimmed mean outperform that for the Huber estimator established in \cite{Z2018}.

Since the winsorized mean and trimmed mean have explicit expressions as presented in \eqref{winsor} and \eqref{trimmed}, we can easily construct the studentized counterparts by plugging in the sample variance. The studentized winsorized mean is given by 
\begin{equation} \label{S}
	S_{\tau, n} = \frac { \sum_{i=1}^n  (f(Y_i  ) - \mu) } { \sqrt{    \sum_{i = 1}^n ( f (Y_i) - \hat{\mu}_{W} )^2 } }
\end{equation}
and the studentized trimmed mean is given by
\begin{equation}
	U_{\tau, n} = \frac { \sum_{i \in \mathcal{N}} (Y_i \IF \{ |Y_i| \leq \tau \} - \mu ) } { \sqrt{ \sum_{i \in \mathcal{N}}  (Y_i \IF \{ | Y_i | \leq \tau \} - \hat{\mu}_T ) ^2 }  },
\end{equation}
where $ \mathcal{N} = \{1 \leq i \leq n: | Y_i | \leq \tau  \} $.
Observe that 
\begin{align}
	S_{\tau, n} & =  \frac { \sum_{i = 1}^n  ( f (Y_i) -\mu ) } {  \sqrt{ \sum_{i=1}^n (  f(Y_i )  - \mu + \mu -  \hat{\mu}_W  )^2 } } \nn \\
	& =  \frac {  \sum_{i = 1}^n  ( f (Y_i) -\mu )  } {  \sqrt{ \sum_{i=1}^n (  f(Y_i )  - \mu )^2   - \frac 1 n   \big( \sum_{i = 1}^n ( f (Y_i) - \mu ) \big)^2  } }  = \frac { S_{\tau, n} ^* } { \sqrt{ 1 - \frac 1n ( S_{\tau, n}^* )^2 } },  \nn
\end{align}
where $ S_{\tau, n}^* $ is the self-normalized winsorized mean defined as 
\begin{equation} \label{self-winsor}
	S_{\tau, n}^* =    \f { \sum_{i=1}^n ( f(Y_i ) - \mu )  } { \sqrt{ \sum_{i=1}^n ( f (Y_i ) - \mu )^2 } }. 
\end{equation}
Similarly, we have for the studentized trimmed mean that
\begin{equation*}
	  U_{\tau, n} = \frac { U^*_{\tau, n} } { \sqrt{1 - \frac {1} { c_{\tau, n}} (U^*_{\tau, n})^2 } } ,
\end{equation*}
where $U_{\tau, n}^*$ is the self-normalized trimmed mean defined as 
\begin{equation} \label{self-trimmed}
	   U_{\tau, n}^* = \frac { \sum_{i \in \mathcal{N}} ( Y_i \IF \{ | Y_i | \leq \tau \} - \mu)   } { \sqrt{\sum_{i \in \mathcal{N}} ( Y_i \IF \{ | Y_i | \leq \tau \} - \mu)^2 } }
\end{equation}
and $c_{\tau, n} = \sum_{i = 1}^n \IF \{ | Y_i | \leq \tau \}$. Since the function $ x / (1 - \f 1n x^2)^{1/2} $ is an increasing function for $ 0 < x < n^{1/2} $, we have
\begin{equation*}
	\P ( S_{\tau, n}  >  x )   = \P \Big( S_{\tau, n}^* > \frac {  x } { \sqrt{1 + \f { x^2 } { n } } } \Big)
\end{equation*}
and 
\begin{equation*}
	  \P ( U_{\tau, n} > x   ) = \P \Big( U_{\tau, n}^* > \frac {  x } { \sqrt{1 + \f { x^2 } { c_{\tau, n} } } }  \Big)
\end{equation*}
Therefore, to investigate the limiting properties of $ S_{\tau, n} $ and $U_{\tau, n}$ is equivalent to investigate that for the simpler self-normalized statistics $ S_{\tau, n}^* $ and $U_{\tau, n}^*$, respectively. 

Before stating our Cram\'er type moderate deviation results, let us first present how the self-normalized winsorized mean and trimmed mean connect with the general self-normalized sum investigated in our main Theorem \ref{thm-1}. First for the self-normalized winsorized mean, though \eqref{self-winsor} presents the form of a self-normalized sum of independent random variables for $ S_{\tau, n}^* $, the expectation of $ f (Y) - \mu $ is slightly deviated from $0$ and needs to be calibrated. Denote
\begin{equation*}
	\tilde{\mu} =  \e  f (Y) , \quad \sigma_1^2 = \e ( f (Y) - \tilde{\mu} )^2, \quad \sigma_2^2  = \e ( f (Y) - \mu )^2.
\end{equation*}
Then we can write
\begin{align}
	S_{\tau, n}^*&  = \frac {  \frac { \sum_{i = 1}^n ( f (Y_i) - \tilde{\mu} ) } { \sqrt n \sigma_1 }  - \frac { \sqrt{n} (\mu - \tilde{\mu} ) } { \sigma_1  } } { \sqrt {  \frac { \sum_{i = 1} ^n  ( f (Y_i) - \mu )^2 } { n \sigma_{2} ^2 }  }  } \cdot \frac{ \sigma_1 }  { \sigma_2 } \nn \\
	& = \frac { S_n - c } { V_n } \cdot \frac{ \sigma_1 }  { \sigma_2 },  \label{general-self}
\end{align}
where $ S_n  $, $ V_n  $ and $c$ are denoted by
\begin{gather*}
	S_n  =  \sum_{i = 1}^n   \frac {   f (Y_i) - \tilde{\mu}  } { \sqrt n \sigma_1 }  , \quad V_n^2 =   \sum_{i = 1} ^n  \frac { ( f (Y_i) - \mu )^2 } { n \sigma_{2} ^2 }\\
	\mbox{and} \quad c =  \frac { \sqrt{n} (\mu - \tilde{\mu} ) } { \sigma_1  } .
\end{gather*}
Therefore,
\begin{equation*}
	\P ( S_{\tau, n}^* > x  ) = \P \Big( \frac { S_n - c } { V_n } >  \frac { \sigma_2 } { \sigma_1 }  x \Big).
\end{equation*}
Note that the random variables involved in $ S_n $ and $ V_n $ are different, which means the existing results for classical self-normalized sums cannot be directly applied. 

As for the self-normalized trimmed mean, note that in the numerator, $ \sum_{i \in \mathcal{N}} ( Y_i \IF \{ | Y_i | \leq \tau \} - \mu) $ in \eqref{self-trimmed} is equal to $ \sum_{i = 1}^n (Y_i - \mu ) \IF \{ |Y_i | \leq \tau  \} $. Similarly in the denominator, $  \sum_{i \in \mathcal{N}} ( Y_i \IF \{ | Y_i | \leq \tau \} - \mu)^2   $ is equal to $ \sum_{i = 1}^n [ (Y_i - \mu) \IF \{ | Y_i | \leq \tau \} ]^2 $. Thus we have 
\begin{equation*}
	\begin{split}
	   U_{\tau, n}^* & = \frac { \sum_{i = 1}^n (Y_i - \mu ) \IF \{ |Y_i | \leq \tau  \}  } { \sqrt{ \sum_{i = 1}^n [ (Y_i - \mu) \IF \{ | Y_i | \leq \tau \} ]^2 } } 
	\end{split} 
\end{equation*}
Denote
\begin{align*}
	\mu_0 &= \e [ (Y - \mu) \IF\{ | Y_i  | \leq \tau \} ], \quad    \sigma_{3} ^2 = \e [ (Y - \mu) \IF\{ | Y_i  | \leq \tau \} -  \mu_0]^2, \\
	 & \mbox{and}\quad \sigma_{4}^2 = \e [ (Y - \mu)^2 \IF\{ | Y_i  | \leq \tau \} ].
\end{align*}
Similar to the self-normalized winsorized mean, we can obtain
\begin{equation*}
	 \P (U_{\tau, n}^* > x ) = \P \Big( \frac { S_n^{\circ}  - \delta } { V_n^{\circ} } > \frac { \sigma_4 } { \sigma_3 } x \Big),  
\end{equation*}
where
\begin{align*}
	   S_n^{\circ} &= \sum_{i = 1}^n \frac { (Y_i - \mu) \IF \{ | Y_i | \leq \tau \} - \mu_0 } { \sqrt n \sigma_3 }, \quad  (V_n^{\circ} )^2 = \sum_{i = 1}^n \frac { [ (Y_i - \mu) \IF \{ |Y_i| \leq \tau  \} ]^2 } { n \sigma_4^2 } \\
	   &\hspace{2cm} \mbox{and} \quad \delta = \frac { \sqrt n \mu_0 } { \sigma_3 }.
\end{align*}

Consequently, our result for general self-normalized sums in Theorem \ref{thm-1} can be directly applied to $  S_{\tau, n}^*   $ and $ U_{\tau, n}^* $ to derive the following bias-corrected Cram\'er-type moderate deviation theorems for the self-normalized winsorized mean and trimmed mean under the fourth moment.

\begin{theorem} \label{thm-winsor-2}
	Assume $  \e [ Y ^4 ] < \infty $. Then there exist an absolute positive constant $ c_1 $ and positive constants $ c_2 $ and $ A $ depending on $\sigma$, $ \e [ |Y|^3 ]$ and $ \e [Y^4 ] $,  such that for 
	\begin{equation}
		\tau \geq c_1 n^{1/6} \max \{ (  \e  [ Y ^4 ]  )^{1/2} / \sigma, ( \e [ Y ^4 ]   / \sigma )^{1/3} \} ,  \label{cond-2}
	\end{equation}
	it holds that
	\begin{align}
		\P ( S_{\tau, n} ^* > x ) & =  [ 1 - \Phi (x) ] \exp \Big\{ - \f  { x^3 \e (Y - \mu)^3  } { 3  \sqrt n  \sigma^3 }  \Big \}  \nn \\
		& \quad \times \Big[ 1 + O_1  \Big(  \frac { (1 + x ^4 )  } { n  }     +  \frac { (1 + x  ) \sqrt n  } { \tau^3  } + \frac { (1 + x ) } { \sqrt n  } \Big)  \Big]  \label{re-2} 
	\end{align}
	uniformly for $ x \in (0, c_2  \min \{ n^{1/4},   \tau^3 n^{- 1/2 }   \}  ) $,  where $ O_1 $ is a  bounded quantity satisfying $ | O_1 | \leq A $. Similar result holds for $ \P ( S_{\tau, n} ^* < - x ) $. 
\end{theorem}

\begin{theorem} \label{thm-trimmed-2}
	Under the conditions of Theorem \ref{thm-winsor-2}, the same result as \eqref{re-2} holds for $ U_{\tau, n}^* $.
\end{theorem}

Observe that under the fourth moment, the general framework Theorem \ref{thm-main} enables us to pin down the bias-corrected term $ \exp \{ - \f  { x^3 \e (Y - \mu)^3  } { 3  \sqrt n  \sigma^3 }  \} $ which depends on the skewness of the underlying distribution. After correcting this skewness in normal approximation, the convergence rate and the converging range significantly improve that  given in Theorem \ref{thm-winsor}, where only third moment is assumed.

The choice of $ \tau  $ should be determined by taking both convergence rate and robustness of estimator into account. We observe from \eqref{re-2} that the ratio $ \P ( S_{\tau, n}^*  > x ) / [ (  1 - \Phi (x) )   \exp  \{ - \f  { x^3 \e (Y - \mu)^3  } { 3  \sqrt n  \sigma^3 }    \}   ]$ converges to  1  for $x \in (0,  o ( \min \{  n^{1/4} , \tau^3 n^{- 1/2 } \} ) )$. The widest possible range $ x \in (0,  o (n^{1/4}) ) $ can be achieved by choosing $  \tau \geq O ( n^{1/4})  $. When $ \tau \leq O (n^{1/3}) $, the larger $ \tau $ is, the faster rate of convergence and wider range of $x$ can be obtained. Yet, once $ \tau  $ exceeds $ O (n^{1/3}) $, our result reduces to the bias-corrected Cram\'er-type moderate deviation for the classical self-normalized sample mean (see Theorem 1.1 in \cite{wang2011refined}), which is reasonable because the winsorized mean and trimmed mean are asymptotically equivalent to the sample mean as $ \tau \to \infty $. It is worth mentioning that when $ O (n^{1/6}) \leq \tau \leq O (n ^{1/3}) $, the ratio $  {\P ( S_{\tau, n} ^*  > x  ) }  / \big(  [ 1 - \Phi(x)  ] \exp \{ - \f  { x^3 \e (Y - \mu)^3  } { 3  \sqrt n  \sigma^3 }   \}  \big) $ converges to $1$ at the rate of $ O ( (1 + x )^4 n^{ - 1 } + (1 + x) \sqrt n \tau^{ - 3} ) $ uniformly for $x \in o(\min\{ n^{1/4}, \tau^3 n^{-1/2} \})$. In this regime of $\tau$, though the convergence rate of winsorized mean and trimmed mean could be slightly slower than that of the classical self-normalized sample mean and the ranges of $x$ for convergence could be narrower, the winsorized mean and trimmed mean provide  robuster inference. We will provide the proof of Theorem \ref{thm-winsor-2} in Section \ref{pf-thm-4.2} and the proof of Theorem \ref{thm-trimmed-2} in Section \ref{pf-thm-trimmed-1}.

Theorem 2.3 in \cite{Z2018} is closely related to ours. They established a Cram\'er-type moderate deviation result for Huber estimator $ \hat{\mu}_H $ defined in \eqref{huber} by using a Bahadur representation for the Huber estimator. Theorem 2.1 in their paper reveals that $ \hat{\mu}_H - \mu   $ is asymptotically close to $ n ^{- 1} \sum_{i = 1}^n f ( Y_i - \mu ) $, where $ f(\cdot) $ is defined by \eqref{f}. Therefore, it is easy to see that the Huber estimator $ \hat{\mu}_{H} $ is close to the winsorized mean $ \hat{\mu}_{W} = n^{ - 1 } \sum_{i=1}^n f (Y_i) $ as $ \tau \to \infty $. Theorem 2.3 in \cite{Z2018} for the Huber estimator can be restated as follows. The notation $ a_n \ll b_n  $ means $ a_n = o (b_n) $ as $n \to \infty$.
\begin{remark} \label{remark-zhou}
	Assume $ \e | Y |^3 < \infty $. \cite{Z2018} proved for $ n^{1/4} \ll \tau \ll n^{1/2} $ that 
	\begin{align}
		&  \frac { \P (   { \sqrt n \sigma^{-1} | \hat{\mu}_H - \mu | } > x  ) } { 2 (1 - \Phi (x)) }  \nn\\
		& = 1 + O(1) \Big\{ \frac { ( \sqrt {\log n}  + x )^3  } { \sqrt n } + \frac { 1 + x } {n^{3/10}}  + \frac { ( 1 + x ) \sqrt n } { \tau^2  }  + e ^{ - O ( \frac {n} {\tau^2 } ) }\Big \} \label{zhou}
	\end{align}
	uniformly for $ 0 \leq x = o ( \min\{ \sqrt{ n} / \tau, \tau^2 / \sqrt n \} ) $. 
\end{remark}

Compared to their condition $ n^{1/4} \ll \tau \ll n^{1/2} $, our condition $ \tau \geq O ( n^{1/6}) $ is less restrictive. Moreover, when $ \tau \gg n^{1/4} $, both of our convergence rate and the associated   converging range of $x$ improve theirs. 	
% For practical use of \cite{Z2018}'s result, one has to estimate the population variance. While to use the bias-corrected limiting distribution in our Theorem \ref{thm-winsor-2}, one has to estimate the population skewness. 
Our improvement mainly relies on the explicit formula of the self-normalized winsorized mean presented in \eqref{general-self} and our fundamental result for general self-normalized sum established in Theorem \ref{thm-1}. In addition, since the higher moment $ \e [Y^4 ] < \infty $ is assumed, after correcting the bias in normal approximation, the convergence rate and the associated range of $x$ could be significantly improved.

The common downside of our bias-corrected result in Theorem \ref{thm-winsor-2} and the normal approximation for Huber estimator by \cite{Z2018} in Remark \ref{remark-zhou} is that the limiting distributions depend on unknown parameters. In real-world applications, if reliable estimations for the unknown parameters are unavailable, we can directly use normal approximation for the self-normalized winsorized mean presented in the following theorem, where only third moment is required and the limiting distribution does not depend on any unknown parameters. 
\begin{theorem} \label{thm-winsor}
	Assume $  \e | Y |^3 < \infty $. Then there exist absolute positive constants $ c_1, c_2 $ and $ A $ such that for 
	\begin{equation}
		\tau \geq c_1 n^{1/4} \max \{ \e | Y |^3 / \sigma^2, ( \e | Y |^3  / \sigma )^{1/2} \} ,  \label{cond}
	\end{equation}
	it holds that
	\begin{align}
		\frac {\P ( S_{\tau, n} ^*  > x  ) } { 1 - \Phi(x) } = 1 + O_1  \Big( \f { (1 +x^3) \e | Y|^3  } { \sigma^3 \sqrt{n} } + \frac { (1 + x )\sqrt{n} \e | Y|^3  } { \sigma \tau^2 }   \Big), \label{re1}
	\end{align}
	uniformly for $ x \in (0, c_2 \min \{  n^{1/6} \sigma^3 /  \e | Y |^3   , \tau^2 \sigma / (\sqrt{n} \e | Y |^3 ) \} ) $, where $ O_1 $ is a bounded quantity satisfying $ | O_1 | \leq A $. Similar result holds for $ \P ( S_{\tau, n} ^* < - x ) $. 
\end{theorem}

\begin{theorem} \label{thm-trimmed}
	Under the conditions of Theorem \ref{thm-winsor}, the same result as \eqref{re1} holds for $U_{\tau, n}^*$
\end{theorem}
It can be observed that the convergence rate and the converging range of $x$ also outperform the results of \cite{Z2018} shown in Remark \ref{remark-zhou}, and our condition on $\tau$ is less restrictive. We relegate the proof of Theorem \ref{thm-winsor} to Section \ref{pf-thm-4.1} and the proof of Theorem \ref{thm-trimmed} in Section \ref{pf-thm-trimmed-2} in the Supplementary Material.

%\begin{enumerate}[(i)]
% 
%	\item Second, by comparing our result in \eqref{re1} and theirs in \eqref{zhou}, it can be observed that our rate of convergence is an improvement over theirs. 
%	\item Third, for the range of $ x $ for convergence, it is easy to see that when $ n^{ 1/4} \ll \tau \leq O( n^{1/3} ) $, our range and theirs are the same, which are both $ x \in (0,  o (\tau^2 / \sqrt n)) $. However, when $ \tau \gg n^{1/3} $, our result shows the ratio of two tail probabilities on the left of \eqref{re1} converges to 1 for $ x \in (0, o (n^{1/6}) ) $, while the ratio in their result \eqref{zhou} converges to 1 for $ x \in (0, o  ( \sqrt{n} /\tau ) ) $, which is narrower than our range. 
%	\item Finally, our result only requires $ \tau \geq O (n^{1/4}) $, and their result requires a stricter condition $ n^{1/4}  \ll \tau \ll n^{1/2} $. Notice that the winsorized mean and the Huber estimator are both close to the sample mean as $ \tau \to \infty $, and our result in \eqref{re1} reduces to the Cram\'er-type moderate deviation for self-normalized sample mean when $ \tau \geq O (n^{1/2}) $, however, their result cannot reduce to the classical Cram\'er-type moderate deviation for standardized sample mean. 
%\end{enumerate}
%	 

\subsection{Simultaneous confidence intervals}
Cram\'er type moderate deviation results are useful in providing theoretical guarantees for a wide spectrum of statistical applications, including the multiple testing procedure and multiple confidence intervals for ultra-high dimensional parameters. For an illustrative example, we will construct simultaneous confidence intervals for the means under the following ultra-high mean model by using the studentized winsorized mean estimator defined in \eqref{S}. We consider
\begin{equation*}
	{\bf Z}_i =   \bm {\mu}+ \bm{\epsilon}_i, \quad i=1,\cdots,n,
\end{equation*}
where $\{ {\bf Z}_1,\ldots, {\bf Z}_n \}$ are i.i.d. observations, $\bm \mu=( \mu_1,\cdots, \mu_p)^T\in \mathbb{R}^p$, and $\{\bm \epsilon_1, \ldots, \bm \epsilon_i \}$ are i.i.d. errors. Denote $\bm \Sigma =  \cov(\bm \epsilon_i ) := (\Sigma_{i j})_{p \times p}$. Assume there exist constants $C_1,C_2$ such that $ \max_{ 1 \leq j  \leq p }\e |Z_{ij}|^3 \leq C_1$ and $\min_{1\leq j\leq p} \Sigma_{jj}  \geq C_2$.
\begin{theorem} \label{thm4.3}
	Assume the dimensionality $p$, the significance level $ \alpha $ and the thresholding parameter $\tau$ satisfying $\log ( p / \alpha)= o(n^{1/3})$ and $\tau \gg n^{1/3}$. Then for $\alpha\in (0,1)$, and $t_0$ satisfying the equation
	\begin{equation*}
		\f{t_0}{1+t_0^2/ n }= \Phi^{-1} \big(1-\f{\alpha}{2p} \big),
	\end{equation*}
	we have
	\begin{equation*}
		\f{\sum_{ i = 1 }^n f (Z_{ij}) }{n} \pm \f{ t_0 }{n} \sqrt{\sum_{ i = 1 }^n \big( f (Z_{ij})  \big)^2 -\f{1}{n} \Big[ \sum_{ i = 1 }^n f (Z_{ij})   \Big]^2  } \triangleq (L_j, U_j ), \quad 1 \leq j \leq p
	\end{equation*}
	are the $1-\alpha-o(1)$ simultaneous confidence intervals for $\big( \mu_j \big)_{j=1}^p $, where $f(\cdot)$ is the function defined in \eqref{f}.
\end{theorem}
The proof of Theorem \ref{thm4.3} will be provided in Section \ref{pf-thm4.3} in the Suppplementary Material.

\section{Proofs}
\label{proof.sec}
\setcounter{equation}{0}
In this section, we present proofs of Theorem \ref{thm-main}, Theorems \ref{thm-2}--\ref{thm_GMC} and Theorems \ref{thm-winsor-2}. Throughout the rest of this section, $A$ and $C$ denote positive absolute constants that may take different values at each appearance.

\subsection{Proof of Theorem \ref{thm-main}} 
We prove the theorem for the two scenarios $0 < x \leq 3$ and $x > 3$, respectively. First, we prove it for $0 < x \leq 3$. For this range, it is sufficient to prove a Berry-Esseen bound as the following proposition will show. The proof of Proposition \ref{prop-small-x} is postponed to Section \ref{proof-prop1} the Supplementary Material.
\begin{proposition} \label{prop-small-x}
	For $ 0 < x \leq 3 $, there exists an absolute constant $A> 0$ such that
	\begin{align}
		\big| \P (S_n > x V_n + c )  -  [1 - \Phi (x+c)] \big| \leq  A  L_{3,n}.
	\end{align}
\end{proposition}
Note that  $ 1 - \Phi (3.6)  \leq 1 - \Phi (x + c) \leq 1$ for $0 < x \leq 3$ and $ |c | \leq x /5 $. Thus, it follows from Proposition \ref{prop-small-x} that for $0 < x \leq 3$,
\begin{align}
	\P ( S_n > x V_n + c  ) = [1 - \Phi (x+c)]  (1 + O (1 + x) L_{3, n} ). \label{eq1}
\end{align}
Moreover, it holds for $0 < x \leq 3$ satisfying \eqref{G1} that
\begin{align*}
	| (\Psi_x^*) ^{-1} - 1 | \leq  A  x^3 L_{3, n}  \leq  A L_{3, n} ,
\end{align*}
which combining with \eqref{eq1} entails that 
\begin{align*}
	\P ( S_n > x V_n + c  ) = [1 - \Phi (x+c)] \Psi_x^*  (1 + O (1 + x) L_{3, n} ) . 
\end{align*}
Consequently, we have
\begin{align*}
	\P ( S_n > x V_n + c  ) = [1 - \Phi (x+c)] \Psi_x^* e^{O_1 R_x} (1 + O (1 + x) L_{3, n} ),
\end{align*}
where the quantity  $ |O_1| \leq A$ for some absolute constant. This completes the proof for $ 0 < x \leq 3 $. 

Next we deal with the case $ x >  3 $. By applying the elementary inequality
\begin{equation}
	1 + s / 2 - s^2 \leq  ( 1 + s  )^{1/2}  \leq 1 +  s/ 2   , \label{ele-ineq}
\end{equation}
for $ s = V_n^2 - 1 $, we obtain
\begin{align*}
	\frac {1} {2} (V_n^2 + 1) - (V_n^2 - 1)^2   \leq  V_n \leq \frac {1} {2} (V_n^2 + 1).
\end{align*}
Therefore, plugging in the above upper and lower bounds yields
\begin{equation}
	\P \( S_n > x V_n + c \) \ge \P \( 2x S_n - x^2 V_n^2 \ge  x^2 + 2xc \) \label{lower-bdd}
\end{equation}
and  
\begin{align}
	&\P \( S_n > x V_n + c \) \label{upper-bdd}\\
	& \leq \P \(   2x S_n - x^2 V_n^2 \ge x^2 + 2xc - x \De_{n}    \) \nn \\
	& \quad +  \P \(  S_n > x V_n + c, \, |V_n^2-1| >  x^{ -1 } ( 1 \lor 6 R_x^{1/2} ) \), \nn
\end{align}
where $\De_{n} = \min   \{ 2x \( V_n^2-1 \)^2,  x^{-1} ( 2 \lor 72 R_x ) \} $ and the notation $  a \lor b $ means the maximum of $a$ and $b$. The upper bound holds because
\begin{align*}
	&\P \( S_n > x V_n + c,\,|V_n^2-1|\le  x^{ -1 } ( 1 \lor 6 R_x^{1/2} )  \)  \\
	&\leq \P \( 2x S_n - x^2 V_n^2 \ge x^2+ 2xc - 2x^2 \( V_n^2 - 1 \)^2 ,\, |V_n^2-1|\le  x^{ -1 } ( 1 \lor 6 R_x^{1/2} ) \) \nn \\
	& \leq \P \( 2x S_n - x^2 V_n^2 \ge x^2 + 2xc - x \De_{n} \). \nn
\end{align*}

The following Propositions \ref{prop-1}--\ref{prop-3} draw an outline of the proof for the case $x > 3$. Their proofs are relegated to Sections \ref{pf_prop1}--\ref{proof-prop4} in the Supplementary Material. 

\begin{proposition}
	There exists an absolute constant $A$ such that
	\begin{align}
		\P \big( 2x S_n - &x^2V_n^2 \ge x^2 + 2xc \big) \label{prop1-re} \\
		=& \[1-\Phi (x+c)\]  \Psi_{x}^*   e^{O_1 R_x} \{ 1+O_2(1+x)L_{3,n}  \},
		\nn
	\end{align}
	for $ x > 3$ satisfying \eqref{G1} and \eqref{G2} and $|c| < x/5$, where $ |O_1| \leq A$ and $ |O_2|  \leq A$.
	\label{prop-1}
\end{proposition}

\begin{proposition}
	There exist absolute constants $A_1$ and $ A_2 $ such that
	\begin{align}
		\P \big( 2x S_n - &x^2V_n^2 \ge x^2 + 2xc - x\Delta_{n} \big) 	\label{prop2-re} \\
		\le& \[1-\Phi (x+c)\]  \Psi_{x}^*  e^{A_1 R_x} \{ 1+A_2(1+x)L_{3,n} \},
		\nn
	\end{align}
	for $ x > 3$ satisfying \eqref{G1} and \eqref{G2}, and $|c| < x/5$.
	\label{prop-2}
	
\end{proposition}

\begin{proposition}
	There exist absolute constants $A_1$ and $ A_2 $ such that
	\begin{align}
		& \P \(S_n  \ge x V_n+c, ~ |V_n^2-1|>  x^{ - 1} ( 1 \lor 6 R_x^{1/2} )  \) \label{prop3-re} \\
		& \hspace{4cm} \le A_1 R_x    [ 1 - \Phi( x + c) ]  \Psi_x^* e^{ A_2 R_x } \nn ,	
	\end{align}
	for $ x > 3$ satisfying \eqref{G1} and \eqref{G2}, and $|c| < x/5 $.
	\label{prop-3}
\end{proposition}

We obtain by substituting the results in Propositions \ref{prop-2}--\ref{prop-3} into \eqref{upper-bdd} that
\begin{align*}
	&\P \( S_n > x V_n + c \) \\
	& \leq \[1-\Phi (x+c)\]  \Psi_{x}^*  e^{A_1 R_x} \{ 1+A_1 R_x + A_2(1+x)L_{3,n}  \} \\
	& \leq \[1-\Phi (x+c)\] \Psi_{x}^*  e^{A R_x} \{ 1  + A (1+x)L_{3,n}  \} 
\end{align*}
which together with the result in Proposition \ref{prop-1} yields the desired result \eqref{main-re} for $x > 3$. The proof is completed.

\subsection{Proof of Theorem \ref{thm-2}} \label{pf-thm3.1}
The main idea is to apply the big-block-small-block technique to construct a general self-normalized sum based on an independent sequence to which our main result Theorem \ref{thm-1} can be applied. Denote $  B_n^2 = \sum_{ i = 1 }^n \e \xi_i^2  $. We first apply Berry-Esseen bound for sum of one-dependent random variables to cope with the case $ 0 \leq  x \leq O  ( \sqrt{ \log n } ) $. Note that
\begin{align*}
	&  \left| \P ( S_n \geq x V_n  )  - \Big[  1 - \Phi \Big( \frac  { x }  { \sqrt{ 1 + 2 \rho_n } }  \Big)  \Big]  \right| \nn \\
	& \hspace{0.5cm} \leq \left |  \P \Big( S_n \geq x B_n ( 1 - n^{ - 1/3 } )^{ 1/2 }  \Big)  - \Big[  1 - \Phi \Big( \frac  { x }  { \sqrt{ 1 + 2 \rho_n } }  \Big)  \Big]  \right|  \nn \\
	&   \hspace{0.5cm}  \quad + \left |  \P \Big( S_n \geq x B_n ( 1 +  n^{ - 1/3 } )^{ 1/2 }  \Big)  - \Big[  1 - \Phi \Big( \frac  { x }  { \sqrt{ 1 + 2 \rho_n } }  \Big)  \Big]  \right|  \nn \\
	&   \hspace{0.5cm} \quad +  \P \Big( | V_n^2 - B_n^2 | > n^{ - 1/3 }  B_n^{ 2 } \Big)  \\
	&   \hspace{0.5cm} : = E_1 + E_2 + E_3.
\end{align*}
Recalling the definition of $\rho_n$ in \eqref{notationonedependentbasic}, we have $ \var (S_n) = (  1 + 2 \rho_n ) B_n^2  $. By noticing the assumptions
\begin{align*}
	\e \xi_i^4 \leq a_1^4, \quad \e \xi_2^2 \geq a_2^2, \quad  a = a_1/a_2, \quad   - 1/2 < \rho \leq  \rho_n \leq 1/2
\end{align*}
and applying the Berry-Esseen bound for sums of one-dependent random variables (see \cite{shergin1980convergence}), we obtain
\begin{align*}
	E_1 & \leq   \frac  { A(\rho)  a^3 } { \sqrt{n}  } +  \Big|  \Phi \Big( \frac  { x }  { \sqrt{ 1 + 2 \rho_n } }  \Big)   -  \Phi \Big( \frac  { x ( 1 - n^{ -  1 / 3  } )^{ 1/2 } }  { \sqrt{ 1 + 2 \rho_n } }  \Big)  \Big| \nn \\
	& \leq A ( \rho ) a^3  ( n^{ - 1/2 } + x^2n^{ - 1/3 } ) \leq  A ( \rho ) a^3 (1+x)^2 n^{ - 1/3 },
\end{align*}
where $ A (\rho) $ is a positive constant depending on $ \rho $ and may take different values at each appearance. In the same manner, this above bound applies to $E_2$ as well.
As for $E_3$, it follows by Chebyshev's inequality that
\begin{align*}
	E_3 \leq ( n^{ 1/3  } B_n^{ - 2 } )^2 \e [ ( V_n^2 - B_n^2  )^2 ]  \leq  A a^4  n^{ - 1/3 }.
\end{align*}
Therefore,
\begin{align}
	\left| \P ( S_n \geq x V_n  )  - \Big[  1 - \Phi \Big( \frac  { x }  { \sqrt{ 1 + 2 \rho_n } }  \Big)  \Big]  \right|  \leq A (\rho) a^4 (1+x)^2 n^{ - 1/3 }, \label{BEboundMdependent}
\end{align}
Moreover, observe that $ 1 - \Phi ( \frac { x } {\sqrt{1 +  2 \rho_n } } ) \geq 1 - \Phi ( \frac { 2 \sqrt 3  } {\sqrt {1 + 2 \rho} } ) $ for $ 0 \leq x  \leq 2 \sqrt 3 $ and 
\begin{equation}
	1 - \Phi (  \frac { x } {\sqrt{1 +  2 \rho_n } }  )  \geq  A(\rho) x^{-1}  \exp \Big\{ -   \frac { x^2 } {  2(1 + 2 \rho) } \Big\}
\end{equation}
for $x > 2 \sqrt 3 $. As a consequence, there exists a constant $ c_{\rho} $ depending on $ \rho $ such that \eqref{re-dep} holds for $ 0 \leq x \leq  c_{\rho}  \sqrt{\log n } $.

Next we turn to the proof for $ x >   c_{\rho} \sqrt{ \log n  } $. Let the length of big blocks be $ l =[n^\alpha] $, and each small block contains only one random variable. Denote $ k = [n / (l + 1 )] $. Without loss of generality, let us assume $ n / (l + 1) $ to be an integer, then the sequence of $\{\xi_i\}_{1 \leq i \leq n } $ can be divided into $ k  $ big blocks and $ k $ small blocks. Observe that  if $n / (l + 1 )$ is not an integer, the sequence of $\{\xi_i\}_{1 \leq i \leq n } $ will be divided into $k + 1$ big blocks and $k$ small blocks,  in which the first $ k  $ big blocks are of size $ l $ and the last one is of size $ n - (l +1 ) [ n / (l + 1 )] $. Although the size of the last big block might be different, our analysis also applies under this scenario. 
For $ 1 \leq j \leq k  $, the $j$-th big block and the corresponding block sums are given by 
\begin{gather*}
	H_j = \{ i: (j-1)( l +1) + 1 \le i \le j( l + 1 ) - 1 \}  ~ ~ \mbox{and} ~~ X_j = \sum_{i \in H_j} \xi_i,   \;\;	Y_j^2 = \sum_{i \in H_j } \xi_i^2,
\end{gather*}
Moreover,  we denote
\begin{gather}
	S_{n1} = \sum_{j=1}^{k} X_j, \quad V_{n1}^2 = \sum_{j=1}^{k} Y_j^2, \quad B_{n1}^2 = \e V_{n1}^2,	\label{notationonedependent}\\
	S_{n2} = \sum_{j=1}^{k} \xi_{j(l+1)},  \quad  V_{n2}^2 = \sum_{j=1}^{k} \xi_{j(l+1)}^2, \quad  B_{n2}^2 = \e V_{n2}^2.   \label{smallb-def}
\end{gather}

Observe that $ \{ X_j \}_{1 \leq j \leq k} $ and $ \{ \xi_{j (l +1 )} \}_{1 \leq j \leq k } $ both consists of independent random variables. As $ S_{n1} $ is the sum of $ \xi_i 's $ in big blocks and it is the main part of $ \sum_{i = 1}^n \xi_i  $, while $ S_{n2} $ corresponds to the small blocks. The big-block-small block technique splits the sum $ \sum_{i = 1}^n \xi_i  $ into two parts $S_{n1}$ and $ S_{n2} $, each of which is a sum of independent random variables.  Let $ \tau = B_{n2}/x $ and we do truncations $ \hat{\xi_i} = \xi_i \IF ( |\xi_i| \le \tau ) $ only for the $ \xi_i $'s in small blocks, that is, $ i = j ( l + 1 ) $ for $ 1 \leq j \leq k $, so
\begin{align}
	\P \( S_n \geq x  V_n  \)
	&\le \P (  \hatS_n  \geq x  \hatV_n   ) + \P ( S_n \geq x V_n, \max_{1 \le j \le k} |\xi_{j(l+1)}| > \tau ) , \label{up-bdd}\\
	\P (S_n \geq x V_n)
	&\ge \P (  \hatS_n  \geq x  \hatV_n    ) - \P (   \hatS_n  \geq x  \hatV_n   , \max_{1 \le j \le k} |\xi_{j(l+1)}| > \tau ),  \label{low-bdd}
\end{align}
where $ \hatS_n = S_{n1} + \hatS_{n2} $, $ \hatV_n^2 = V_{n1}^2 + \hatV_{n2}^2  $ with $ \hatS_{n2} = \sum_{j = 1}^k \hat{\xi}_{ j ( \ell + 1 ) } $ and $ \hat{V}_{n2}^2 =  \sum_{j = 1}^k \hat{\xi}_{ j ( \ell + 1 ) }^2 $.
For a positive number $ d_1 > 0 $, we have the upper bound
\begin{align}
	\P ( \hatS_n \geq x \hatV_n ) &\le \P \big(S_{n1} \geq x V_{n1} - d_1 n^{ -\f{ \alpha}{2}  } x B_n\big) \label{uu-bb} \\
	& \quad + \P \big(\hatS_{n2} \geq  d_1 n^{ -\f{ \alpha}{2}  }  x B_n  \big) \nn
\end{align}
and the lower bound
\begin{align}
	\P ( \hatS_n \geq x \hatV_n ) &\ge \P \big(S_{n1} \geq x V_{n1} +   d_1 n^{ -\f{ \alpha}{2}  }  x B_n\big) \label{ll-bb}\\
	&\quad  - \P \big( S_{n1} \geq x V_{n1}, ~V_{n1}^2 < B_n^2 / 4\big) \nn \\
	&  \quad - \P \big( \hatS_{n2} < - d_1 n^{ -\f{ \alpha}{2}  } x B_n + x (\hatV_n - V_{n1}), ~V_{n1}^2 > B_n^2 / 4 \big) . \nn
\end{align}
\noi We can obtain the following bounds for the terms involved in \eqref{uu-bb} and \eqref{ll-bb}. The proofs of Propositions \ref{lemma-bigb} and \ref{lemma-smallblock} are given in Sections \ref{pf_prop6.1} and \ref{pf_prop6.2} in the Supplementary Material.

\begin{proposition} 	\label{lemma-bigb}
	There exists an absolute positive constant $ d_0 $ and a constant $A(\rho, d_0)$ depending on $ \rho   $ and $ d_0  $ such that, for $ d_1  =  \kappa_{\rho} a^2 $ with some sufficiently large constant $ \kappa_{\rho} $,
	\begin{align}
		\P \Big( S_{n 1 }  &\geq x V_{n 1 } + d_1 n^{ - \alpha / 2  } x B_n  \Big) \label{big-block} 	 \\
		&  =  \Big[ 1 - \Phi( \f{x}{ \sqrt{ 1 + 2 \rho_n } } ) \Big]  \Big(  1 + O_1 \Big(  \frac { a^4 x^4  } { n^{ 1- \alpha } }+  \frac {  a^2  x^2  } { n^{ \alpha/ 2 } } +\f{a^3x}{n^{\f{1-\alpha}{2}}} \Big)  \Big),  \nn	
	\end{align}
	uniformly for $ x \in ( 2 , d_0 a^{ -1  } \min \{  n^{ \alpha/4 }, n^{ ( 1 - \alpha ) / 4  } \}  ) $, where $ | O_1 | \leq  A(\rho, d_0) $.
	A similar result holds for $ \P ( S_{n1} \geq x V_{n1} - d_1 n^{ - \alpha/2  } x B_n ) $. Moreover,
	\begin{align}
		\P \Big( S_{n1} &\geq x V_{n1}, ~V_{n1}^2 < \f{B_n^2}{4} \Big) \label{error_low1} \\
		&\le  \f{ A a^4 x^4} {n^{1-\alpha}} \Big[ 1 - \Phi\Big( \f{x}{ \sqrt{ 1 + 2 \rho_n } } \Big) \Big] \exp {\Big\{ A  a^4 \f{x^4} {n^{1-\alpha}} + A  a^3 \f  {x^3} {\sqrt{n}} \Big\} } , \nn
	\end{align}
	uniformly for $ x \in ( 2 , d_0 a^{ - 1} n^{\f12 (1-\alpha)}) $.	
\end{proposition}

\begin{proposition} 	\label{lemma-smallb} \label{lemma-smallblock}
	For $ d_1 =\kappa_{\rho}  a^2   $ with $\kappa_{\rho}\geq 10$, we have
	\begin{align}
		\P \big(\hatS_{n2} > d_1 n^{ - \alpha/2  }  x B_n  \big)  &  \le \exp  \{  -     \kappa_{\rho} x^2 / 6   \}, \label{error_up}  \\
		\P \big( \hatS_{n2} < - d_1 n^{ -  \alpha /2 } x B_n  &  +  x (\hatV_n - V_{n1}), ~V_{n1}^2 >  B_n^2 /4 \big)  \label{error_low2}  \\
		&   \le \exp \{  -  \kappa_{\rho}  x^2 /  14    \} . \nn
	\end{align}
\end{proposition}

Note that for $ x > c_{\rho} \sqrt{ \log n  } $ and sufficiently large $ \kappa_{\rho} $,
\begin{align}
	\frac { \exp\{  - \kappa_{\rho} x^2 / 6   \}  } { 1 - \Phi \big( \frac { x } { \sqrt {1 + 2 \rho_n} }  \big) }  \leq \exp \{ - \kappa_{\rho} x^2 / 12  \} \leq n^{ - 1/4 }.   \label{error}
\end{align}
By choosing $ \alpha = 1/2 $ and combining \eqref{big-block}-\eqref{error}, we obtain
\begin{align}
	\P \Big(\hatS_n \geq x \hatV_n \Big)
	=  \Big[ 1 - \Phi\Big( \f{x}{ \sqrt{ 1 + 2 \rho_n } } \Big) \Big]   \Big( 1 + O \Big( \frac {a ^4  x^2    } { n^{ 1/4 } }  \Big) \Big) \label{trun-re}
\end{align}
uniformly for $ x \in ( c_{\rho} \sqrt{ \log n  } , d_0 a^{ -1  }  n^{ 1/ 8 }  ) $.

In addition, we can bound the error terms in \eqref{up-bdd} and \eqref{low-bdd} as follows. The proof of Proposition \ref{prop-error} is shown in Section \ref{pf_prop7} in the Supplementary Material. 
\begin{proposition} \label{prop-error}
	Under the conditions in Theorem \ref{thm-2}, we have  for $  x \in (c_{\rho} \sqrt{\log n}, d_0 a^{-1} n^{1/8} )$ that 
	\begin{equation}
		\begin{split}
			& \P \Big( \hatS_n \geq x \hatV_n, \max_{1 \le j \le k} |\xi_{j(l+1)}| > \tau \Big) \\
			& \quad \leq A(\rho, d_0) \f{ a^ 4 (1+x)^4} {n^{1/2}}\Bigg[ 1 - \Phi\Big( \f{x}{ \sqrt{ 1 + 2 \rho_n } } \Big) \Bigg]   
		\end{split}
		\label{error-bound1}
	\end{equation}
	and
	\begin{equation}
		\begin{split}
			&\P \Big( S_n \geq x V_n, \max_{1 \le j \le k} |\xi_{j(l+1)}| > \tau \Big) \\
			& \quad \leq A(\rho, d_0)  \f{a^4 (1 + x)^4}{n^{1/2}}\bigg[ 1-\Phi\Big(  \f{x}{\sqrt{1+2\rho_{ n }}}  \Big) \bigg] .
		\end{split}
		\label{error-bound2}
	\end{equation}
\end{proposition}
Consequently, substituting \eqref{trun-re}--\eqref{error-bound2} into \eqref{up-bdd} and \eqref{low-bdd} yields the desired result \eqref{re-dep}. This completes the proof of Theorem \ref{thm-2}.

\subsection{Proof of Theorem \ref{thm-blocksum}} \label{pf-thm3.2}
The proof for Theorem \ref{thm-blocksum} again builds on the big-block-small-block technique, and also exploits a lemma in \cite{shao1996weak} to replace the weakly dependent big blocks and small blocks by independent random variables, respectively. 
We begin the proof by introducing three essential lemmas from the literature. Lemma \ref{lemma-mo1} (Theorem 4.1 of \cite{shao1996weak}) and Lemma \ref{lemma-mo2} (Theorem 10.1.b of \cite{lin2011probability}) concern the bound of moments under weak dependence while Lemma \ref{lemma-indep} (Lemma 2.1 of \cite{berbee1987convergence}) shows that a $ \beta $-mixing sequence of random variables can be replaced by an independent sequence of random variables in a domain whose measure is at least $ 1 - \sum_{i = 1}^n  \beta^{(i)}$.
\begin{lemma} \label{le-dep-moment}
	(Theorem 4.1 in \cite{shao1996weak}.)
	Let $\{X_i, i \ge 1\} $ be a sequence of zero-mean random variables with $ \e | X_i |^r  \le \mu^r$ for $ r > 2 $ and $ \mu >0  $. Assume that mixing condition \eqref{g-beta} holds, then
	\begin{equation*}
		\e \Big[ \Big| \sum_{i=k}^{i=k+m} X_i \Big|^{r'} \Big] \le C m^{r'/2} \mu^{r'}\;,
	\end{equation*}
	for any $ 2 \le r' < r, m \ge 1$ and $ k\ge 0 $, where $C$ is a constant that depends on $  r', r, a_1, a_2 $ and $ \tau $.
	\label{lemma-mo1}
\end{lemma}

\begin{lemma} \label{lemma-mo2}
	(Theorem 10.1.b of \cite{lin2011probability}.)
	Assume $ \{ X_i \}_{ i \geq 1 } $ is a sequence of random variables and $ \beta(n) $ is the $ \beta $-mixing coeffient defined in \eqref{beta-coeff}. Denote by $ \sigma_{1}^k $ and $ \sigma_{k + n}^{\infty}$ the $\sigma$-fields generated by $ \{X_i\}_{1 \leq i \leq k} $ and $ \{ X_{i} \}_{i \geq k + n} $, respectively.
	For $ X \in L_p ( \sigma_1^k ) $ and $Y \in L_q (\sigma_{k+n}^{\infty}) $ with $ p, q, r \ge 1 $ and $ \f1p +\f1q + \f1r = 1 $, we have
	\begin{equation}
		| \e X Y - \e X \e Y | \le 8 \beta(n)^{1/r} \| X \|_p \| Y \|_q \;.  \label{esti-cov}
	\end{equation}
\end{lemma}

For two random variables (or vectors) $X$ and $Y$, define
\begin{equation*}
	\beta( X, Y) = \frac 12 \sup_A  \Big(  \big( \P_{X, Y} -  \P_X \times \P_Y  \big) (A) -  \big( \P_{X, Y} -  \P_X \times \P_Y  \big) (A^{c}) \Big).
\end{equation*}

\begin{lemma} \label{lemma-indep}
	(Lemma 2.1 of \cite{berbee1987convergence}.)
	Let $ \{  \xi_i, 1 \le i \le n \} $ be a sequence of random variables on the same probability space and define $ \beta^{(i)} = \beta ( \xi_i,\, (\xi_{i+1}, \dots, \xi_n) ) $. Then the probability space can be extended with random variables $ \tilde{\xi_i} $ distributed as $ \xi_i$ such that $ \{ \tilde{\xi_i} \}_{ 1 \le i \le n}  $ are independent and
	\begin{equation}
		\P \Big(\xi_i \neq  \tilde{\xi_i}, ~\mbox{for some} ~1 \le i \le n \Big)  \le \beta^{(1)} + \dots + \beta^{(n-1)}.  \label{betaindepapproerror}
	\end{equation}
\end{lemma}

Recall the definition of block sums $\{Y_j\}_{1 \leq j \leq k }$ in \eqref{blocksum}. We set the size of big blocks as $ m_1 = [n^{\alpha_1}] $ for some $ 0 < \alpha_1 < 1 - \alpha $ and the size of small blocks as $1$. Denote $ k_1 = k / (m_1 + 1 )$, where $k = [n / l]$.  For simplicity of presentation, we assume $k / (m_1 +1 )$ to be an integer, as explained in the proof of Theorem \ref{thm-2}. The $ u $-th big block is given by 
\begin{equation*}
	I_u = \Big\{ j : ( m_1 + 1 )( u -1) + 1 \le j \le (m_1 + 1 )u - 1  \Big\}  , \quad {\rm for}\;  1\le u \le k_1 .
\end{equation*}
Define
\begin{gather*}
	S_k = \sum_{j = 1}^k Y_j, \quad V_k^2 = \sum_{j = 1}^k Y_j^2, \quad  \xi_u  = \sum_{j \in I_u } Y_j ,\quad  \eta_u^2 = \sum_{j \in I_u } Y_j^2, \\
	S_{k,1} = \sum_{ u =1}^{k_1} \xi_u, \quad V_{k,1}^2 = \sum_{ u =1}^{k_1} \eta_u^2 ,  \quad S_{k,2} = \sum_{u=1}^{k_1} Y_{ u (m_1 + 1)} , \quad  V_{k,2}^2 = \sum_{ u =1}^{k_1} Y_{ u (m_1 + 1)}^2,\\
	B_n^2 = \sum_{j=1}^k \e Y_j^2,\quad  B_{n,1}^2 = \sum_{ u =1}^{k_1} \sum_{j \in I_u } \e  Y_j^2,\quad  B_{n,2}^2 = \sum_{ u =1}^{k_1} \e Y_{ u (m_1 + 1)}^2 .
\end{gather*}
\noi By Lemma \ref{lemma-mo2}, it is easy to see that under the mixing condition \eqref{g-beta},
\begin{equation}
	\f { B_n^2 } { B_{n,1}^2 } = 1 + O\big(\f{k_1}{k} \big) = 1 + O(n^{- \alpha_1}) , \quad   \f { \e S_{k,1}^2 } { B_{n,1}^2 } = 1 + O( n^{- \alpha} ) + O( n^{- \alpha_1} ) . \label{mo-ratio}
\end{equation}

\noi Let $ \hatY_j = Y_j \IF ( |Y_j| \le b ) $, where $ b = B_{n,2}/(1 + x) $. Parallel to one-dependent case, we separate the big blocks and small blocks right after truncating the terms inside the small blocks. Denote
\begin{align*}
	&\hatS_{k, 2}  = \sum_{u = 1}^{k_1} \hatY_{ u ( m_1 + 1  ) }, \quad \hatS_{k} = S_{k, 1} + \hatS_{k, 2}\\
	{\rm and}\quad  &\hatV_{k,2}^2 =  \sum_{u = 1}^{k_1} \hatY_{ u ( m_1 + 1  ) }^2 , \quad \hatV_{k}^2 = V_{k, 1}^2 + \hatV_{k, 2}^2 .
\end{align*}
It is straightforward that
\begin{align}
	\P \big(  S_k \geq x V_k  \big) &\le \P \big( \hatS_k \ge x \hatV_k \big) + \P \Big( S_k \ge x V_k , \max_{1 \le u \le k_1} |Y_{  u (m_1 + 1) }| > b \Big) , \label{trun1}   \\
	\P \big(  S_k \geq x V_k \big) & \ge \P \big( \hatS_k \ge x \hatV_k \big)	- \P\Big( \hatS_k \ge x \hatV_k, \max_{1 \le u  \le k_1} |Y_{ u (m_1 + 1) }| > b \Big) . \label{trun2}
\end{align}
Further, for the main term $ \P ( \hatS_k \ge x \hatV_k  ) $, we choose $\vp = d_1 n ^{ - \alpha_1 / 2 } \log n  $ with  a positive number $d_1 > 0$ and obtain
\begin{align}
	\P \Big( \hatS_k \ge x \hatV_k  \Big) &\le \P \big( S_{k,1} \ge x V_{k,1} - \vp x B_n \big) + \P\Big( \hatS_{k,2} > \vp x B_n \Big), \label{upboundfortruncatedinbetamixing} \\
	\P \Big( \hatS_k \ge x \hatV_k  \Big) & \ge \P \big( S_{k,1} \ge x V_{k,1} + \vp x B_n \big) \label{lowerbound}\\
	&\qquad  - \P\Big( S_{k,1} \ge x V_{k,1} , V_{k,1}^2  \le \f14 B_n^2 \Big)  \nn \\
	&\qquad - \P\Bigg( \hatS_{k,2} - \f{x \hatV_{k,2}^2 }{ \hatV_k + V_{k,1} } < -\vp x B_n, V_{k,1}^2 > \f14 B_n^2 \Bigg). \nn
\end{align}
The estimate for dominated terms $ \P (  S_{k,1} \ge x V_{k,1} - \vp x B_n ) $ and $ \P (  S_{k,1} \ge x V_{k,1} + \vp x B_n )  $  is presented in the following lemma and the proof will be shown in Section \ref{pf-thm3.2-prop}. 
\begin{proposition} \label{prop-thm3.2}
	Assume $ \vp = d_1 n^{\alpha_1 / 2 } \log n  $ for a positive number $d_1 > 0$ and $ \alpha_1 \leq \alpha \tau $.  Under the conditions of Theorem \ref{thm-blocksum}, there exist a positive constant $ c_0 $ depending on $ d_1, \mu_1 / \mu_2, a_1, a_2 , \alpha$ and $\tau$ such that 
	\begin{align}
		&\P  ( S_{k,1} \ge x V_{k,1} \pm \vp x B_n ) \big/ [1- \Phi(x)  ]  \label{main1} \\
		=& 1 + O\bigg(\f{(1 + x)^4}{ n^{1- \alpha- \alpha_1} }  + \f{1 + x }{n^{(1 - \alpha - \alpha_1 )/2 }} + \f{( 1 + x )^2 }{ n^{\alpha }  } + \f{ ( 1 + x )^2 \log n } {n^{ \alpha_1/2 }} \bigg)  \nn
		%& \P( S_{k,1} \ge x V_{k,1} + \vp x B_n) \nn \\
		%&=  [1- \Phi(x)] \big(1 + O(1) (\f{(1 + x)^4}{ n^{1- \alpha- \alpha_1} } + \f{(1 + x )^3}{\sqrt{n}} + \f{1 + x }{n^{(1 - \alpha - \alpha_1 )/2 }} + \f{( 1 + x )^2 }{ n^{\alpha }  } + \f{ ( 1 + x )^2 } {n^{ \alpha_1/2 }} ) \big)   \label{main2}
	\end{align}
	uniformly for $ x \in (0, c_0 \min \{ n^{( 1 - \alpha - \alpha_1 )/4 }, n^{\alpha \tau/2}, n^{\alpha/2}, (\log n)^{-1/2}n^{\alpha_1/4}\} )$, 
\end{proposition}

For small-block-related error terms, $ S_{k,2} $ and $ V_{k,2}^2 $ can be replaced with the sum of independent random variables by Lemma \ref{lemma-indep}. In addition, following a similar proof to Proposition \ref{lemma-smallb}, we obtain that under $ \vp = d_1 n^{- \alpha_1/2} \log n $ for $ d_1 $ being some positive number depending on $ \mu_1/\mu_2 $, there exist positive numbers $ C_1 $ and $ C_2 $ depending on $a_1$, $a_2$, $ \mu_1 $, $ \mu_2 $, $\alpha$ and $ \tau $ such that	
\begin{align}
	& \P\Big( \hatS_{k,2} > \vp x B_n \Big)  \label{er1} \\
	& \qquad \le \exp  \{- C_1 ( 1 + x  )^2  d_1 \log n   \} +  C_2 \exp  \{ - a_2 n^{\alpha \tau  } / 2   \},    \nn \\
	& \P\Big( \hatS_{k,2} - \f{x \hatV_{k,2}^2 }{ \hatV_k + V_{k,1} } < -\vp x B_n, V_{k,1}^2 > \f14 B_n^2 \Big) \label{er2} \\
	& \qquad \leq  \P \Big(   \hatS_{k,2} - \f{x \hatV_{k,2}^2 }{ B_n /2 } < -\vp x B_n \Big)\nn \\
	& \qquad \le \exp \{- C_1 ( 1 + x  )^2  d_1 \log n   \} +  C_2 \exp \{ - a_2 n^{\alpha \tau  } / 2 \},    \nn
\end{align}
where the error term $C_2 \exp \{ - a_2 n^{\alpha \tau  } / 2 \}  $ is obtained by applying Lemma \ref{lemma-indep} to replace $S_{k, 2}$ and $V^2_{k, 2}$ with sum of independent random variables and the fact that (see \cite{berbee1987convergence})
\begin{equation}
	\beta \big( \{ X_i\}_{i\in \mathcal{J}'}, \{ X_i\}_{ i\in\mathcal{J}'' } \big) \leq \beta \big( \{ X_i\}_{i\leq k}, \{ X_i\}_{ i\geq n + k } \big) \leq \beta( n ), \label{betamixingproperty}
\end{equation}
for any sets $ \mathcal{J}'\subset \{ i\leq k \}$, $\mathcal{J}''\subset \{ i\geq n+k \}$.

Regarding big-block-related error terms, by using Chebyshev's inequality and Taylor expansion, we control the error term in the same manner as the proof of \eqref{Omega1}. We can obtain, 
\begin{align}
	\P \Big( S_{k,1} &\geq x V_{k,1}, V_{k,1}^2  \le \f14 B_n^2 \Big)  \label{er3} \\
	&\le A_5 \f{(1 + x)^4}{n^{1 - \alpha - \alpha_1 }} \Big( 1 - \Phi \big[x \big(1 + O \big(n^{-\alpha} + n^{- \alpha_1 } \big) \big)\big] \Big) \nn
\end{align}
uniformly for  $ x \in (3, c_0 \min \{ n^{( 1 - \alpha - \alpha_1 )/4 }, n^{\alpha \tau/2} \} )$.
When $ 0 < x \leq 3 $, it follows from Lemmas \ref{le-dep-moment} and \ref{lemma-indep} and Chebyshev inequality that under condition \eqref{g-beta}, 
\begin{align*}
	\P \Big( S_{k,1} \geq x V_{k,1}, V_{k,1}^2  \le \f14 B_n^2 \Big) & \leq \P (  V_{k,1}^2  \le \f14 B_n^2 ) \\
	&  \leq   a_1 n^{ \alpha_1 }e^{ - a_2 n^{\alpha \tau }  } + \exp \Big \{  - \frac { B_{n,1}^2 - B_n^2 / 4 } { 2 \sum_{u = 1}^{k_1}  \e \eta_u^4 }  \Big \} \\
	& \leq a_1 n^{ \alpha_1 }e^{ - a_2 n^{\alpha \tau }  }  + \exp \{ - A_1  n^{ 1 - \alpha - \alpha_1 } \mu_2 / \mu_1  \} \\
	& \leq  \frac { A_2 (1 + x)^4 } { n^{ 1 - \alpha - \alpha_1  } } \big[ 1 - \Phi (x) \big].
\end{align*}

Substituting \eqref{main-p}, \eqref{er1}, \eqref{er2}, and \eqref{er3} into \eqref{upboundfortruncatedinbetamixing} and \eqref{lowerbound} yields
\begin{align}
	&\P (\hatS_k \ge x \hatV_k )/ [1- \Phi(x)  ] \label{truncatedresultinbetamixing} \\
	=& 1 + O \bigg(\f{(1 + x)^4}{ n^{1- \alpha- \alpha_1} }  + \f{1 + x }{n^{(1 - \alpha - \alpha_1 )/2 }} + \f{( 1 + x )^2 }{ n^{\alpha }  } + \f{ ( 1 + x )^2 \log n } {n^{ \alpha_1/2 }} \bigg)  \nn
\end{align}
uniformly for $ x \in (0, c_0 \min \{ n^{( 1 - \alpha - \alpha_1 )/4 }, n^{\alpha \tau/2}, n^{\alpha/2}, (\log n)^{-1/2}n^{\alpha_1/4}\} )$.

To avoid redundance, we omit the analysis of the truncation errors in \eqref{trun1} and \eqref{trun2} as their proofs share the same fashion with \eqref{trun-error1} and \eqref{trun-error2}. Consequently, (\ref{truncatedresultinbetamixing}) also holds for $\P(S_k \ge x V_k)$. Finally, we need to balance the error terms by choosing $ \alpha_1 $ and seeking the best convergence rate or largest range for convergence. As a result, we choose $ \alpha_1 = (1 - \alpha)/2 $ when $ (1 - \alpha)/2 \le \alpha \tau$, and choose $ \alpha_1 = \alpha \tau $ when $ (1 - \alpha)/2 > \alpha \tau$, and then the desired result follows. This completes the proof for Theorem \ref{thm-blocksum}.

\subsection{Proof of Theorem \ref{thm_GMC}} \label{pf-thm3.3}

The main idea is to use one-dependent random variables to approximate $ \{Y_j \}_{1 \leq j \leq k } $ and then apply Theorem \ref{thm-2}. Recall that $ m = [n^{\alpha}] $ and $ k  =[n /m ] $. Let
\begin{equation*}
	\tilde{Y}_j=\e(Y_j|\varepsilon_l, mj-2m+1\leq l\leq mj)
\end{equation*}
and
\begin{equation*}
	\tilde{T}_k=\frac{\sum_{j=1}^{k}\tilde{Y}_j}{ ( \sum_{ j = 1 }^k \tilde{Y}_j^2 )^{ 1/2 } }.
\end{equation*}
%	\begin{equation*}
%	\tilde{T}_k=\frac{\sum_{j=1}^{k}\tilde{Y}_j}{\tilde{V}_k},\quad \;\; {\rm where} \; \; \tilde{V}_k^2=\sum_{j=1}^{k}\tilde{Y}_j^2.
%	\end{equation*}
As $ \{\varepsilon_t \}_{t \in \mathbbm{Z}} $ are i.i.d. random variables, $\{ \tilde{Y}_j \}_{j\geq 1}$ are one-dependent. Note that by conditional Jensen's inequality, for $ 2 \leq r \leq 4  $,
\begin{align*}
	&\Big\Vert X_i  - \e  \Big( X_i \Big\vert \vp_{\ell}: m j  - 2 m + 1 \leq \ell \leq i \Big) \Big\Vert_r^r \nn \\
	& \hspace{0.8cm} =  \e \left\{  \Big| \e \Big[ X_i - G_k \Big( \mathscr{F}_{ m j  - 2m }^* , \vp_{ m j  - 2m + 1 }, \ldots, \vp_{i} \Big)  \Big\vert \mathscr{F}_i  \Big] \Big|^r  \right\} \nn \\
	& \hspace{0.8cm} \leq \e \Big[ \Big| X_i  -  G_k \Big( \mathscr{F}_{ m j  - 2m }^* , \vp_{ m j  - 2m + 1 }, \ldots, \vp_{i}  \Big)  \Big|^r  \Big] \nn \\
	& \hspace{0.8cm} \leq  \Big[ \Delta_r ( i - m j + 2 m   ) \Big] ^{r} ,
\end{align*}
which together with the assumption \eqref{GMC} and the fact that $ m ( j - 1) + 1 \leq i \leq m j  $ for $ i \in  H_j $ yields
\begin{align}
	\parallel Y_j- \tilde{Y}_j \parallel_r & \leq \sum_{ i \in H_j  } \Big\Vert   X_i - \e \Big(  X_i \Big\vert \vp_{\ell} : i - 2m + 1 \leq \ell \leq i  \Big)  \Big\Vert_r \label{difference} \\
	& \leq m a_1 e^{ -  a_2 m^{\tau}} .   \nn
\end{align}
The above bound shows that $ \{Y_j\}_{1 \leq j \leq k } $ can be well approximated by the one-dependent sequence $ \{ \tilde{Y}_j \}_{1 \leq j \leq k} $. We can derive Theorem \ref{thm_GMC} by aggregating the following two propositions. The proofs of Propositions \ref{prop_approx} and \ref{prop_error} will be provided in Sections \ref{pf-prop8} and \ref{pf-prop9} in the Supplementary Material, respectively. 

\begin{proposition}  \label{prop_approx}
	Under conditions of Theorem \ref{thm_GMC}, we have there exists a positive constant $ d_0 $ depending on $\tau, \alpha, w_1, a_1$ and $a_2$ such that 
	\begin{align} \label{main-p} 
		\P ( \tilde{T}_k \geq x ) 
		= [ 1 - \Phi (x) ] \Big( 1 + O \Big( \frac { 1 + x^2  }  { n^{ ( 1 - \alpha )/ 4  } } +  \frac { 1 +  x^2  } { n^{ \alpha } } \Big)  \Big),   
	\end{align}
	Uniformly for $ x \in ( 0, d_0 \min \{ n^{ (1 - \alpha)/8 }, n^{\alpha/ 2}\} )   $.
\end{proposition}

\begin{proposition}  \label{prop_error}
	Under conditions of Theorem \ref{thm_GMC}, we have for $ x > 0 $, 
	\begin{align}
		\P (T_k \geq x ) &  \leq \P ( \tilde{T}_k \geq x  - C_1 n^{-1} ) + C_2 ( e^{ - a_2 n^{\alpha \tau}}  + e^{ - O (n^{1 - \alpha})}) \label{error_1p} \\
		\mbox{and} \quad \P (T_k \geq x ) &  \geq \P ( \tilde{T}_k \geq x + C_1  n^{-1} ) - C_2 ( e^{ - a_2 n^{\alpha \tau}}  + e^{ - O (n^{1 - \alpha})}). \label{error_2p}
	\end{align}
\end{proposition}

Applying Proposition \ref{prop_approx} to $ \P ( \tilde{T}_k \geq x+ O(  n^{-1} ) ) $ yields
\begin{align*}
	\P \Big(  \tilde{T}_k &\geq x + O (n^{-1})  \Big)   \\
	& =  \Big[ 1 - \Phi \Big( x +O(n^{-1}) \Big)  \Big]  \Big( 1 + O \Big( \frac { 1 + x^2  } { n^{ ( 1 - \alpha  ) / 4   } } + \frac { 1 + x^2  } { n^{ \alpha  } } \Big) \Big) \nn \\
	&  =  \big[ 1 - \Phi (x) \big] \Big( 1 + O \Big( \frac { 1 + x^2  } { n^{ ( 1 - \alpha  ) / 4   } } + \frac { 1 + x^2  } { n^{ \alpha  } } \Big) \Big)
\end{align*}
uniformly for $ x \in (0, d_0 \min  \{ n^{ (1 - \alpha)/8 }, n^{\alpha/ 2}\}) $.
Note that the error term $e^{ - a_2 n^{\alpha \tau}}  + e^{ - O (n^{1 - \alpha})}  $ decays at an exponential rate, which is always faster than the polynomial rate. By substituting the above result into \eqref{error_1p} and \eqref{error_2p} leads to the desired result \eqref{GMC_result}. The proof of Theorem \ref{thm_GMC} is completed.

\section*{Acknowledgements}
The authors would like to thank the Co-Editors, Associate Editor, and referees for their constructive comments that have helped improve the paper significantly.

\begin{supplement} 
	%\sname{Supplement A}
	\textbf{Supplement to ``Refined Cram\'er Type Moderate Deviation Theorems for General Self-normalized Sums with Applications to Dependent Random Variables and Winsorized Mean"}.
	%	\slink[doi]{COMPLETED BY THE TYPESETTER}
	%	\sdatatype{.pdf}
	The supplement \cite{GSS2021} contains all the technical details of the proofs.
\end{supplement}
%% if your bibliography is in bibtex format, uncomment commands:
\bibliographystyle{imsart-nameyear} % Style BST file (imsart-number.bst or imsart-nameyear.bst)
\bibliography{general_CTMD}       % Bibliography file (usually '*.bib')

\newpage
\appendix
\setcounter{page}{1}
\setcounter{section}{0}
\renewcommand{\theequation}{A.\arabic{equation}}
\setcounter{equation}{0}

\begin{center}{\bf \large Supplementary Material to ``Asymptotic Distributions of High-Dimensional Distance Correlation Inference"}
	
	\bigskip
	
	Lan Gao, Qi-Man Shao and Jiasheng Shi
\end{center}

\noindent This Supplementary Material contains all technical details of proofs. Section \ref{SecA} provides the proofs of all the propositions, Theorems \ref{thm-winsor-2}, \ref{thm-winsor} and \ref{thm4.3}, and Corollary \ref{cor}. Section \ref{SecB} presents the proofs of lemmas. 
\section{Proof of Propositions} \label{SecA}

\subsection{Proof of Proposition \ref{prop-small-x}} \label{proof-prop1}
The proof borrows some techniques from \cite{wang2011refined}. The main idea is to first truncate the random variables, and then to apply Berry-Esseen bound for U-statistics and Berry-Esseen bound for sum of independent random variables for the upper bound and the lower bound, respectively. Define 
\begin{gather*}
	\hatX_i = X_i \mathbbm{1} ( | (1 + x ) X_i | \le 1 ) , \quad \hatS_n = \sum_{i=1}^n \hatX_i , \\
	\hatY_i = Y_i \mathbbm{1} ( | (1 + x) Y_i | \le 1 )  ,  \quad  {\hatV_n}^2 = \sum_{i=1}^n \hatY_i^2.
\end{gather*}
It is easy to see
\begin{align}
	& \big|\P (S_n > x V_n+c) - [ 1 - \Phi(x + c) ] \big| \label{truncation}  \\
	& \leq \big|\P (\hatS_n > x \hatV_n+c) - [ 1 - \Phi(x + c) ] \big|  \nn \\
	&\qquad + \sum_{i=1}^n \P \( | (1 + x ) X_i| > 1 \) + \sum_{i=1}^n \P \( | ( 1 + x ) Y_i| > 1 \)  \nn  \\
	& \leq \big|\P (\hatS_n > x \hatV_n+c) - [ 1 - \Phi(x + c) ] \big| + (1 + x)^3 L_{3,n}. \nn
\end{align}

\noi Let $ \hatB_n^2 = \sum_{i=1}^n \e \hatY_i^2 $ and $ K_i = \hatY_i^2 - \e \hatY_i^2  $, then it follows from the basic inequality \eqref{ele-ineq} that
\begin{align}
	\P (\hatS_n > x \hatV_n+c) &\ge \P\( 2 x \hatS_n -x^2 \hatV_n^2 > x^2 + 2 x c  \), \label{lb}
\end{align}
and
\begin{align}
	& \P (\hatS_n > x \hatV_n+c) \label{ub} \\
	&\le \P\Bigg( \hatS_n - \f{x}{2 \hatB_n} \( \hatV_n^2 - \hatB_n^2 \) + \f{x}{\hatB_n^3} \( \hatV_n^2 - \hatB_n^2 \)^2 > x \hatB_n + c \Bigg) \nn \\
	&= \P \Bigg( \sum_{i=1}^{n} (U_i - \e U_i) + \f{x}{\hatB_n^3} \sum_{i \ne j} K_i K_j \ge  x \hatB_n + c - \sum_{i=1}^{n} \e U_i \Bigg),    \nn
\end{align}
where $ U_i = \hatX_i - \f{x}{2\hatB_n} K_i +\f{x}{\hatB_n^3} K_i^2 $. Observe that $ \sum_{i=1}^{n} (U_i - \e U_i) + \f{x}{\hatB_n^3} \sum_{i \ne j} K_i K_j $ is a U-statistic. We will apply the Berry-Esseen bound for U-statistics established by \cite{alberink2000berry}.
Obviously, $ | \hatB_n^2  -  1 | \leq    (1 + x ) L_{3, n} \leq c_1 $ by \eqref{G1}, hence $ 3/4 \leq \hatB_n^2 \leq 5/4 $ for $ c_1 \leq 1/2  $. We obtain that for $ 0 < x \leq 3  $,
\begin{gather*}
	\Big| \sum_{i=1}^{n} \e U_i \Big| \leq  (1 + x)^2 \sum_{i = 1}^n \e  |X_i|^3  + \f{x}{\hatB_n^3} \sum_{i = 1}^n \e \hatY_n^4  \le AL_{3,n}, \\
	\Big| \sum_{i=1}^{n}\var( U_i )  -  1 \Big|  \leq A L_{3, n}, \\
	\sum_{i=1}^n \e |U_i|^3 \le \sum_{i=1}^n A\Big( \e |\hat{X}_i|^3 +\frac{x^3}{\hat{B}_n^3}\e |K_i|^3 + \frac{x^3}{\hat{B}_n^9}\e |K_i|^6   \Big) \le AL_{3,n},\\
	\mb{and} \quad  \sum_{ i \neq j } \var( \f{x}{\hatB_n^3} K_iK_j ) = \sum_{ i \ne j } \f{x^2}{\hat{B}_n^6} \e K_i^2 \e K_j^2  \le AL_{3,n}^2.
\end{gather*}

It follows from the Berry-Esseen bound for U-statistics (see Theorem 1 in \cite{alberink2000berry}) that
\begin{align*}
	& \bigg|  \P \Big( \sum_{i=1}^{n} (U_i - \e U_i) + \f{x}{\hatB_n^3} \sum_{i \ne j} K_i K_j \ge  x \hatB_n + c - \sum_{i=1}^{n} \e U_i \Big)   - [ 1 - \Phi (x + c )  ] \bigg| \\
	& \leq A L_{3, n}  + \Big| [ 1 - \Phi ( x \hatB_n + c - \sum_{i = 1}^n \e U_i ) ] - [ 1 - \Phi(x + c) ]  \Big| \\
	& \leq A L_{3, n} + A \Big( ( \hatB_n - 1 ) + \Big| \sum_{i = 1}^n \e U_i \Big| \Big) \leq A L_{3,n },
\end{align*}
which together with \eqref{ub} yields the upper bound
\begin{equation}
	\P (\hatS_n > x \hatV_n+c) - [1 - \Phi(x + c) ] \le A L_{3,n}. \label{ub_xb}
\end{equation}
Regarding the lower bound, denote $ \hat{W_i} = 2x \hatX_i - x^2 \hatY_i^2 $. It is easy to find that for $ 0 <  x \leq 3 $,
\begin{gather*}
	\Big | \sum_{i=1}^n \e \hat{W_i}  +  x^2 \Big|  \leq A x L_{3,n}, \quad \Big| \sum_{i=1}^n \var(\hat{W_i}) -  4x^2 \Big| \leq A  x^2 L_{3,n}, \\
	\mb{and}  \qquad \sum_{i=1}^n \e | \hat{W_i} |^3 \leq A x^3  L_{3,n}.
\end{gather*}
Consequently, the Berry-Esseen bound for sums of independent random variables implies
\begin{align}
	& \big|\P( 2 x \hatS_n -x^2 \hatV_n^2 > x^2 + 2 x c  ) - [1 - \Phi( x + c )] \big|  \\
	&   \leq \Big| \P \Big( \frac { \sum_{i = 1}^n  ( \hat{W}_i -  \e  \hat{W}_i  ) } {  (\sum_{i = 1}^n \var(\hat{W}_i) )^{1/2}  } >    x + c - A L_{3, n}   \Big)  \nn \\
	&  \hspace{0.7cm}   - [ 1 + \Phi ( x + c - A L_{3, n}) ]  \Big|  +  | \Phi ( x + c - A L_{3, n} ) - \Phi ( x + c ) |  \nn \\
	&   \le   \frac { A \sum_{ i = 1 }^n \e | \hat{W}_i |^3 } { \big(\sum_{i = 1}^n \var(\hat{W}_i)  \big)^{3/2}}+ A L_{3,n} \leq A L_{3, n}\nn ,
\end{align}
which together with \eqref{lb} yields
\begin{equation}
	\P (\hatS_n > x \hatV_n+c) - [1 - \Phi(x + c) ] \ge - A L_{3,n}. \label{lb_xb}
\end{equation}
Combining \eqref{truncation}, \eqref{ub_xb} with \eqref{lb_xb}, we obtain that for $ 0 < x \leq 3 $,
\begin{align*}
	\big| \P (S_n > x V_n + c )  -  [1 - \Phi (x+c)] \big| \leq  A  L_{3,n}.
\end{align*}
This completes the proof.

\subsection{Proof of Propositions \ref{prop-1}} \label{pf_prop1}
Before starting to prove Proposition \ref{prop-1}, we first collect some notations related to the conjugated method, which is the main tool to prove Propositions \ref{prop-1}--\ref{prop-2}. For $1 \leq i \leq n$, let $ W_i = 2 x X_i - x^2 Y_i^2 $ and $ (\xi_i, \eta_i) $ be independent random vectors with distribution
\begin{equation}
	V_i \(x,y\) = \f { \e  \{  e^{\la W_i } \mathbbm{1} (X_i\le x, \, Y_i\le y ) \} } { \e e^{ \la W_i }}, \label{change-measure}
\end{equation}
Denote $ \WT{W_i} = 2x \xi_i - x^2 \eta_i^2 $. It holds that
\begin{align*}
	\e \WT{W_i} & = \frac { \e W_i e^{\lambda W_i}  } { \e e^{ \lambda W_i }  }, \\
	\var \WT{W_i} & = \frac {  \e  W_i^2   e^{ \lambda W_i }  } { \e e^{ \lambda W_i }   } - ( \e \WT{W_i} )^2,  \\
	\e | \WT{W_i} |^3 & = \frac { \e | W_i |^3  e^{ \lambda W_i }    } { \e e^{ \lambda W_i }   }.
\end{align*}
We establish the expansion of the above moments in the Lemma below. The proof of Lemma \ref{lemma-1} will be given in Section \ref{pf-lemma-1}.
\begin{lemma} 	\label{lemma-1}	
	Let $ W_i = 2x X_i - x^2 Y_i^2 $. For $ \f14 \leq  \la \leq  \f34 $ and $ x > 0  $ satisfying \eqref{G2}, there exists an absolute constant $A$ such that
	\begin{align}
		\e e^{\la W_i} &= 1 + 2 \la^2 x^2 \e X_i^2 - \la x^2 \e Y_i^2 + \f43 \la^3 x^3 \e X_i^3 - 2 \la^2 x^3 \e X_iY_i^2 + O_1 R_{x,i}, \label{0-mo} \\
		&= \exp \left\{ 2 \la^2 x^2 \e X_i^2 - \la x^2 \e Y_i^2 + \f 4 3 \la^3 x^3 \e X_i^3 - 2 \la^2 x^3 \e X_iY_i^2 + O_2 R_{x,i}        \right\},  \nn
	\end{align}
	\begin{equation}
		\e W_ie^{ \la W_i} = 4 \la x^2 \e X_i^2 - x^2 \e Y_i^2 + 4\la^2 x^3 \e X_i^3 - 4 \la x^3\e X_iY_i^2 + O_3 R_{x,i}, \label{1-mo}
	\end{equation}
	\begin{equation}
		\e W_i^2 e^{\la W_i} = 4x^2 \e X_i^2 + 8 \la x^3 \e X_i^3 - 4 x^3 \e X_iY_i^2 + O_4 R_{x,i},  \label{2-mo}
	\end{equation}
	\begin{equation}
		\e |W_i|^3 e^{\la W_i} = O_5 x^3\( \e |X_i|^3 + \e |Y_i|^3 \) + O_6 R_{x,i}, \label{3-mo}
	\end{equation}
	where $ | O_i  | \leq A $ for $ i = 1, \ldots, 6 $.
\end{lemma}

The next lemma is Lemma 2.1 in \cite{wang2011refined} with some modifications.
\begin{lemma}     \label{lemma-0}
	We have for $x$ satisfying \eqref{G2} that
	\begin{align*}
		(1 + x) ^4 (\e X_i^2)^2 &  \le 2 \de_{x,i}, \\
		(1 + x)^5 \e X_i^2   \e |X_i|^3 & \le 2 \de_{x,i}, \\
		( 1 + x )^6 (\e |X_i|^3)^2 & \le \de_{x,i},
	\end{align*}
	and similar results hold for $Y_i$. In addition, if $x$ also satisfies (\ref{G1}), then
	\begin{equation*}
		(1  + x )^4 L_{3,n}^2  \le 2\de_x.
	\end{equation*}
\end{lemma}

By Lemmas \ref{lemma-1} and \ref{lemma-0}, it is readily seen that under \eqref{G2},
\begin{align*}
	\e \WT{W_i} &= x^2 \( 4\la \e X_i^2- \e Y_i^2 \) + x^3 \( 4\la^2 \e X_i^3 - 4\la \e X_iY_i^2 \) + O_1 R_{x,i},  \\
	\var \WT{W_i} &= 4x^2 \e X_i^2+ x^3 \( 8\la \e X_i^3 - 4\e X_iY_i^2 \)+ O_2 R_{x,i},
\end{align*}
and
\begin{equation*}
	\e \big| \WT{W_i} \big|^3 = O_3 x^3 \( \e \left|X_i\right|^3 + \e \left|Y_i\right|^3 \) + O_4 R_{x,i}.
\end{equation*}

\noi Put $ m_n = \sum_{i=1}^n \e \WT{W_i}$, $ \si_n^2 = \sum_{i=1}^n Var \WT{W_i} $, $ v_n = \sum_{i=1}^n \e \big| \WT{W_i} \big|^3 $. Consequently, we obtain
\begin{align}
	m_n &= \( 4\la - 1 \) x^2 + x^3 \( 4\la^2 \sum_{i=1}^n \e X_i^3 - 4\la \sum_{i=1}^n \e X_iY_i^2 \) + O_1 R_x, \label{expectation} \\
	\si_n^2 &= 4x^2 + x^3 \( 8\la \sum_{i=1}^n \e X_i^3 - 4\sum_{i=1}^n \e X_iY_i^2 \) + O_2 R_x, \label{2nd-moment}
\end{align}
and
\begin{equation}
	v_n = O_3 x^3 L_{3,n} + O_4 R_x. \label{3rd-moment}
\end{equation}

\noi Define $ m \(\la\) = \sum_{i=1}^n \log{ \e e^{\la W_i} } $, therefore, $ m_n = m' (\la) $ and $\si_n^2 = m''( \la )$. Before proving Proposition \ref{prop-1}, let us present the following lemma that will be applied in the proof. 
\begin{lemma} 	\label{lemma-2}
	For $ x $ satisfying \eqref{G1} and \eqref{G2}, if $|\de \(x\)| < \f12 x^2 $ and $|\de_0 \(x\)|<\f12 x^2 $, then the equation
	\begin{equation*}
		m' \(\la\) = x^2 + \de \(x\),
	\end{equation*}
	has a unique solution $\la_\de$. In addition, $ \la_{\de} $ satisfies $ \f14 < \la_\de < \f34 $ and
	\begin{equation}
		\left| \la_\de - \( \f12 + \f {\de(x)}{4x^2} \) + x \( \la_\de^2 \sum_{i=1}^{n} \e X_i^3 - \la_\de \sum_{i=1}^{n} \e X_iY_i^2 \) \right| \le A x^{-2} R_x,
		\label{pro-solu}
	\end{equation}
	and
	\begin{equation}
		\left| \la_\de - \la_{\de_0} - \frac{ \delta\left(x\right) -
			\de_0 (x) }{ 4x^2 } \right| \le A \{ x^{-2} R_x + \left| \de (x) - \de_0 (x) \right| x^{-1} L_{3,n} \}.
		\label{e6}
	\end{equation}
	Moreover, we have
	\begin{equation}
		m\( \la_\de \) =  \( 2\la_\de^2 - \la_\de \) x^2 + x^3\( \f43 \la_\de^3  \sum_{i=1}^{n} \e X_i^3 - 2\la_\de^2 \sum_{i=1}^{n} \e X_iY_i^2 \) + O_1 R_x ,
		\label{e4}
	\end{equation}
	and
	\begin{equation}
		\left| m \( \la_\de \) - m \( \la_{\de_0} \) \right| \le A \( R_x + \left| \de \(x\) - \de_0 \(x\) \right| \).
		\label{e5}
	\end{equation}
\end{lemma}
The proof of Lemma \ref{lemma-2} will be shown in Section \ref{pf-lemma-2}. Now we are ready to prove Proposition \ref{prop-1}. 
The main idea is to use the conjugate method (see (4.9) in \cite{petrov1965probabilities}).
Assume $ \la_1 $ is the solution to $ m'(\la) = x^2 +2xc $.
Write $ \WT {S_n} = \sum_{i=1}^n \WT {W_i} $, $ U_n = \(\WT {S_n} - m_n \)/\si_n$. It is well-known that by the conjugate method
\begin{align}
	\P \( 2x S_n - x^2 V_n^2 \ge x^2 + 2xc \)
	=& \exp \{ m(\la_1) \} \e e^{-\la_1 \WT {S_n}} \mathbbm{1}( \WT {S_n} \ge x^2 + 2xc ) \nn \\
	=& \exp \{ m(\la_1) - \la_1 m_n\} \e e^{-\la_1 \si_n U_n} \mathbbm{1} ( U_n \ge 0 ). \nn
	% &\le \exp \{ m(\la_1) - \la_1 m_n\} \{ I_1 + I_2 \} ,
\end{align}

Let $Z$ be a standard normal random variable. Define
\begin{align*}
	I_1  &=  \sup \limits_{ x \in \mathbbm{R} } | \P(U_n \leq x )  -  \P (Z \leq x)|, \\
	I_2  & = \e e^{-\la_1 \si_n Z} \mathbbm{1} ( Z \ge 0 )  .
\end{align*}
Integration by parts gives
\begin{align*}
	\big| \e e^{-\la_1 \si_n U_n} \mathbbm{1} ( U_n \ge 0 )  - \e e^{-\la_1 \si_n Z} \mathbbm{1} ( Z \ge 0 )  \big|  \leq  2 I_1.
\end{align*}
As a result,
\begin{align}
	& \big|  \P \( 2x S_n - x^2 V_n^2 \ge x^2 + 2xc \)  - \exp \{ m(\la_1) - \la_1 m_n\} I_2 \big|  \label{con-re} \\
	& \hspace{1cm} \leq  2 \exp \{ m(\la_1) - \la_1 m_n\}  I_1. \nn
\end{align}

\noi Applying the Berry-Esseen theorem to $I_1$ and by the fact of \eqref{2nd-moment} and \eqref{3rd-moment}, it is easy to find that for $x > 3 $ satisfying \eqref{G1},
\begin{equation}
	I_1 \le A \f { v_n }{ \si_n^3 } \le A ( L_{3,n} + x^{-3} R_x ).
	\label{I1-bound}
\end{equation}
Regarding $ I_2 $, we have
\begin{equation*}
	I_2 = \frac{e^{\lambda_1^2 \sigma_n^2/2}}{\sqrt{2\pi}} \int_{ \lambda_1 \sigma_n }^{\infty} e^{ - t ^2 / 2 } dt = \f{1}{ \sqrt{2\pi} } { 1-\Phi(\la_1 \si_n) \over \varphi (\la_1 \si_n) },
\end{equation*}
where $ \phi (x) $ is the standard normal density function.
Let $ \psi (x) = \f { 1 - \Phi(x) }{ \varphi(x) } $. For $ x > 3 $,
\begin{equation}
	x^{ -1 }/2 \le \psi ( x )  \le x^{-1} \quad \mb{ and } \quad |\psi ' (x)| = |x \psi (x) - 1| \le x^{-2}. \label{property1}
\end{equation}
By \eqref{2nd-moment} and the fact that $ 1/4 < \lambda_1 < 3/4 $,
$$
| \la_1 \si_n - 2 \la_1 x | = \f { \la_1 |\si_n^2 - 4x^2| }{ \si_n + 2x } \le A \( x^2 L_{3n} + x^{-1} R_x \),
$$
therefore, under \eqref{G1} for small constant $ c_1 $, there exists some $\theta^*$ between $2\lambda_1 x$ and $\lambda_1 \sigma_n$, such that
\begin{align}
	I_2 &= \f {1}{ \sqrt{2\pi} } \{ \psi \( 2 \la_1 x\) + \psi' \( \theta^* \) \( \la_1 \si_n - 2 \la_1 x \) \} \label{I2-bound}  \\
	&= \f {1}{ \sqrt{2\pi} }  \psi \( 2 \la_1 x\) \left\{ 1+O_1 ( x L_{3,n} +  x^{-2} R_x ) \right\} \nn  \\
	&= e^{ 2 \la_1^2 x^2 } \[ 1-\Phi \( 2 \la_1 x \) \] \left\{ 1 + O_1 (  x L_{3,n} +  x^{-2} R_x ) \right\}.
	\nn
\end{align}

Since $ 1 /4 < \lambda_1 < 3/4 $, we have for $x > 3 $
\begin{align*}
	e^{ 2 \la_1^2 x^2 } \[ 1-\Phi \( 2 \la_1 x \) \]  \geq  A x^{  -1 },
\end{align*}
which yields
\begin{align*}
	I_1  \leq  A  e^{ 2 \la_1^2 x^2 } \[ 1-\Phi \( 2 \la_1 x \) \]   ( x L_{3, n} + x^{ -2 } R_x ).
\end{align*}
Consequently, by \eqref{con-re}--\eqref{I2-bound} and $ m_n = m'(\la_1) = x^2 + 2 x c  $, we arrive at
\begin{align}
	\P (2x &S_n - x^2V_n^2 \ge x^2 + 2xc )   \label{mid-re} \\
	&= \exp \big\{ m( \la_1 ) + ( 2\la_1^2 - \la_1 ) x^2 - 2 \la_1 xc  \big\} \nn \\
	& \quad \quad \times \[ 1 - \Phi (2\la_1 x) \]\left\{ 1 + O_1 x L_{3,n} +O_2 x^{-2} R_x \right\} \nn \\
	&= \exp \big\{ m( \la_1 ) + \f12 (x+c)^2  -\la_1 ( x^2 + 2xc ) \big\} \nn \\
	& \quad \quad \times \[ 1 - \Phi (x+c) \]  e^{O_1 R_x} \{ 1 + O_2 x L_{3,n} \}, \nn
\end{align}
and the last equality holds because
\begin{equation*}
	|  \la_1  -  \gamma | \le A_1 x L_{3,n} + A_2 x^{-2} R_x.
\end{equation*}
By applying \eqref{pro-solu} and \eqref{e4} in Lemma \ref{lemma-2}, where $ \gamma = (1 + c /x ) /2  $.  Noticing that $ x^4 L_{3,n}^2 \le 2\de_x < 2R_x $ by Lemma \ref{lemma-0}, we obtain for $x$ satisfying \eqref{G1} and \eqref{G2} that
\begin{equation*}
	(2\lambda_1^2-\lambda_1)x^2+\frac{1}{2}(x+c)^2-\lambda_1(x^2+2xc)= 2x^2 (\lambda_1 -\gamma)^2 \leq O_1 R_x.
\end{equation*}
Thus \eqref{e4} implies
\begin{align}
	& \; m( \la_1 ) + \f12 (x+c)^2  -\la_1 ( x^2 + 2xc ) \label{transform} \\
	=\; &  \f43 \ga^3 x^3 \sum_{i=1}^n \e X_i^3 - 2 \ga^2 x^3 \sum_{i=1}^n \e X_i Y_i^2 + O_1 R_x. \nn
\end{align}
Finally, we arrive at the desired result \eqref{prop1-re} by substituting \eqref{transform} into \eqref{mid-re}. The proof is completed.

\subsection{Proof of Proposition \ref{prop-2}} \label{proof-prop3}
The proof again relies on the conjugate method and a randomized concentration inequality for independent random variables. Recalling the definitions of $\la_1,  \widetilde{S_n}  $ and $ U_n  $ at the very beginning of Section \ref{pf_prop1}, we have by the conjugate method
\begin{align}
	&\P \( 2x S_n - x^2 V_n^2 \ge x^2 + 2xc - x \De_{n} \) \label{P_1+P_2}\\
	&\hspace{1cm}  = \exp \{ m(\la_1) \} \e e^{-\la_1 \WT {S_n}} \mathbbm{1}(  \WT {S_n} \ge x^2 + 2xc - x \WT\De_{n} ) \nn \\
	&\hspace{1cm}  = \exp \{ m(\la_1) - \la_1 m_n\} \e e^{-\la_1 \si_n U_n} \mathbbm{1} ( U_n \ge -\f {x \WT\De_{n}}{\si_n}) \nn \\
	& \hspace{1cm} :=  P_1 + P_2 , \nn
\end{align}
where
\begin{align*}
	P_1 &= \exp \{ m(\la_1) - \la_1 m_n\} \e e^{- \la_1 \si_n U_n } \mathbbm{1} ( U_n > 0 ),\\
	P_2 &= \exp \{ m(\la_1) - \la_1 m_n\} \e e^{- \la_1 \si_n U_n } \mathbbm{1} ( -\f {x \WT\De_{n}}{\si_n} \le U_n \le 0 ),
\end{align*}
and $\WT\De_{n} =  \min \{ 2 x ( \sum_{i = 1}^n \eta_i^2 - 1 )^2 , x^{ -1 }  (2 \lor 72 R_x) \} $. We already established in \eqref{mid-re} that
\begin{align}
	P_1 & = \exp \{ m(\la_1) + \( 2\la_1^2 - \la_1\) x^2 - 2 \la_1 xc \} \label{P1} \\
	& \hspace{2cm} \times [ 1 - \Phi ( 2 \la_1 x  ) ] ( 1 + O_1 x L_{3,n} + O_2 x^{-2} R_x ).  \nn
\end{align}
As for $ P_2 $, because $ x \WT\De_{n} \le 2 + 72 R_x $,
\begin{equation}
	P_2 \le e^{ 2 + 72 R_x  } \exp \{ m\( \la_1 \) - \la_1 m_n \} \P \Big( -\f {x \WT\De_{n}}{\si_n} \le U_n \le 0 \Big). \label{P2}
\end{equation}
As $ U_n  $ is a normalized sum of independent random variables, we apply the randomized concentration inequality by \cite{shao2016cramer} (Theorem 4.1 therein) to estimate the concentration probability involved in $ P_2 $. As $  \Delta_n \leq 2 x ( V_n^2 -1  )^2 $, it is sufficient to bound the concentration probability as if $ \Delta_n = 2 x ( V_n^2 -1  )^2  $ and $ \widetilde{\Delta}_n = 2 x ( \sum_{i = 1}^n \eta_i^2 - 1 )^2 $. Denote $ Z_i = Y_i^2 - \e Y_i^2 $ and define $\De_{n}^{(i)}$ as $ \De_{n}^{(i)}= 2x(\sum_{ j \ne i } Z_j)^2  $. $\WT\De_{n}^{(i)}$ is defined in the same manner as $\De_{n}^{(i)}$ but with $ \{Y_i \}_{ 1 \leq i \leq n  }$ replaced by $\{ \eta_i  \}_{1 \leq i \leq n}$. We obtain
\begin{equation}
	\P \Big( -\f {x \WT\De_{n}}{\si_n} \le U_n \le 0 \Big) \le 17 \f {v_n}{ \si_n^3 } + \f {5x}{ \si_n }  \e \big| \WT\De_{n} \big| + \f {2 x}{ \si_n^2}  \sum_{i=1}^n \e \big| \WT {W_i} \{ \WT\De_{n} - \WT\De_{n}^{(i)} \} \big|. \label{BE}
\end{equation}
Furthermore, it follows from the distribution of $ (\xi_i, \eta_i ) $ defined in \eqref{change-measure} that
\begin{align}
	\e \big| \WT\De_{n} \big| &= \exp \{ - m (\la_1) \} \e \left\{ |\De_{n}| e^{ \la_1 \sum_{i=1}^n W_i }  \right\} \label{esti-1} \\
	&   \le 2 x \exp \{ - m (\la_1) \} \e \left\{  \( V_n^2 -1 \)^2 e^{ \la_1 \sum_{i=1}^n W_i }  \right\} , \nn
\end{align}
and
\begin{align}
	& \e \big| \WT {W_i} \{ \WT\De_{n} - \WT\De_{n}^{(i)} \} \big| \label{esti-2}  	 \\
	& = \exp \{ - m (\la_1) \} \e \left\{ \big| W_i \{ \De_{n} -  \De_{n}^{(i)}  \}  \big| e^{ \la_1 \sum_{i=1}^n W_i } \right\}  \nn \\
	& \leq  2 x \exp \{ - m (\la_1) \} \e \Big \{ |W_i| \Big| ( V_n^2 - 1 )^2 - \Big(\sum_{ j \ne i } Z_j\Big)^2 \Big| e^{ \la_1 \sum_{i=1}^n W_i }  \Big \}   \nn \\
	& =  2 x \exp \{ - m (\la_1) \} \e \Big\{ |W_i| \Big| Z_i^2 + 2 Z_i \sum_{j \ne i} Z_j \Big| e^{ \la_1 \sum_{i=1}^n W_i }  \Big\}. \nn 	
\end{align}
The lemma below presents the bounds of \eqref{esti-1} and \eqref{esti-2} and the proof will be presented in Section \ref{pf-lemma0}.
\begin{lemma} \label{lemma-3}
	For $ x $ satisfying \eqref{G1} and \eqref{G2} and $ 1 /4 < \la_1 < 3 / 4  $, there exists an absolute constant $A$ such that for $ x > 3 $,
	\begin{align}
		\e \left\{ \( V_n^2 -1 \)^2 e^{ \la_1 \sum_{i=1}^n W_i } \right\} &\le A_1 \exp \{ m( \la_1 ) \} x^{-2} R_x e^{ A_2 R_x },  \label{error1} \\
		\sum_{ i =1 }^n \e \Big\{ \big|W_i\big| \Big| Z_i^2 + 2 Z_i \sum_{j \ne i} Z_j \Big| e^{ \la_1 \sum_{i=1}^n W_i }  \Big\} &\le A_1 \exp \{ m( \la_1 ) \} x^{-2} R_x e^{ A_2 R_x }. \label{error2}
	\end{align}
	In general, assume $0 < r < r_0 < 1  $ for a  constant $r_0$. For a number $    \omega >   r_0 $, there exist constants $A_1 $ and $ A_2$ depending on $r_0   $ and $ \omega $ such that
	\begin{align}
		& \e \Big\{ \big( V_n^2 -1 \big)^2  e^{ \sum_{i=1}^n (2rxX_i- \omega r x^2Y_i^2) } \Big\}         \label{error1generalversion}  \\
		&   \le \f{A_1 R_x}{x^2} \exp\Big\{  (2r^2- \omega r )x^2-2 \omega r^2 x^3\sum_{ i =1 }^n  \e X_iY_i^2  +   \frac{4}{3}r^3x^3\sum_{ i =1 }^n \e X_i^3+ A_2 R_{x}   \Big\}  .  \nn
	\end{align}
\end{lemma}

Consequently, by assembling \eqref{P2}--\eqref{error2} and \eqref{expectation}--\eqref{3rd-moment}, it holds that
\begin{align}
	P_2 &\le A \exp \{ A R_x \}  \exp \{ m( \la_1 ) - \la_1 m_n \}  ( L_{3,n} +  x^{-1} R_x )  .  \label{Porp2p2}
\end{align}
Moreover, by observing that $ m_n = m' (\la_1 ) = x^2 + 2 x c  $ and $  1 - \Phi ( 2 \la_1 x ) \geq  C x^{ -1 } e^{ - 2 \la_1^2 x^2   } $ for $ x > 2 $ and $ 1/ 4 <  \la_1 < 3 / 4 $, we obtain
\begin{align}
	P_2 & \leq  A \exp \{ A R_x \} \exp \{ m( \la_1 ) + (2 \la_1^2 - \la_1) x^2 - 2 \la_1 xc \}   \\
	& \quad \quad \times \[ 1 - \Phi( 2 \la_1 x ) \] \{ A_4 x L_{3,n} + A_5 R_x \},  \nn
\end{align}
which combined with \eqref{P1} yields
\begin{align*}
	&\P ( 2x S_n - x^2 V_n^2 \ge x^2 + 2xc - x \De_{n} )  \\
	&\le \exp \{ m( \la_1 ) + (2 \la_1^2 - \la_1) x^2 - 2 \la_1 xc \}   \[ 1 - \Phi( 2 \la_1 x ) \] e^{A R_x} (1 + A  x L_{3,n} ) \\
	&= \exp \big\{ m( \la_1 ) + \f12 (x+c)^2  -\la_1 ( x^2 + 2xc ) \big\} \\
	& \hspace{2cm} \times \[ 1 - \Phi (x+c) \]  e^{A_3 R_x}  ( 1 + A_4  x L_{3,n} ),
\end{align*}
where the last equality is derived by a similar procedure to \eqref{mid-re}. Finally, we arrive at the desired result \eqref{prop2-re} by using \eqref{transform}. The proof is completed.

\subsection{Proof of Proposition \ref{prop-3}} \label{proof-prop4}
We use the truncation technique to estimate the error term. Let $ B = 2500 \lor  200c_0 $. Here the constants $2500$ and $200$ are set large enough for simplicity of proof, because we are not pursuing the best possible constants. The integrating region can be partitioned into three parts as follows,
\begin{equation}
	\P\Big( S_n > x V_n + c, | V_n^2 - 1 | >   x^{-1} ( 1 \lor 6 R_x^{1/2})  \Big) \le \sum_{i=1}^3 \P \( ( S_n, V_n^2 ) \in \Omega_i \) ,
\end{equation}
where $ \{ \Omega_i, i=1,2,3 \} $ are given by
\begin{align*}
	\Omega_1 &= \big\{ (u,v) \in (R,R^+): ~u > x \sqrt{v} + c,~  1 +  x^{-1} ( 1 \lor 6 R_x^{1/2}) < v \le B \big\}, \\
	\Omega_2 &= \big\{ (u,v) \in (R,R^+):~ u > x \sqrt{v} + c,~v < 1 - x^{-1} ( 1 \lor 6 R_x^{1/2}) \big\}, \\
	\Omega_3 &= \big\{ (u,v) \in (R,R^+): ~ u > x \sqrt{v} + c, ~ v > B  \big\}.
\end{align*}

\vskip 0.5em
For the first part, we choose
\begin{equation*}
	r_1= x + c ,\quad   t_1 =  \f{2}{5}x(x + c).
\end{equation*}
By Chebyshev's inequality we obtain
\begin{align}
	&\P \big( \big( S_n , V_n^2 \big) \in \Omega_1 \big)  \label{er-1} \\
	&\le x^2 \exp \Big\{ - \inf_{ (u,v) \in \Omega_1 } ( r_1 u - t_1 v ) \Big\}\e \left\{ \( V_n^2 - 1 \)^2 e^{ r_1 S_n - t_1 V_n^2 } \right\}. \nn
\end{align}
Obviously,
\begin{align}
	& \inf_{ (u,v) \in \Omega_1 } ( r_1 u - t_1 v )  \label{inf} \\
	&= r_1 \Big( c + x \sqrt { 1 +   x^{-1} ( 1 \lor 6 R_x^{1/2}) } \Big) - t_1 \big( 1 + x^{-1} ( 1 \lor 6 R_x^{1/2})  \big) . \nn
\end{align}
It follows from applying (\ref{error1generalversion}) in Lemma \ref{lemma-3} with $ r  = r_1 / (2 x ) $ and $ t =t_1 / x^2$ (It is easy to verify the conditions for \eqref{error1generalversion} are satisfied when $|c| \leq x / 5 $) that
\begin{align}
	&\e \{ \( V_n^2 - 1 \)^2 e^{ r_1 S_n - t_1 V_n^2 } \}  \label{er-11}	 \\
	& \leq  \frac{A_1 R_x } { x^2 } \exp \Big\{  \f { r_1^2 }{2} - t_1 + \f {r_1^3}{6} \sum_{i=1}^n \e X_i^3 - r_1 t_1 \sum_{i=1}^n \e [X_i Y_i^2] + A_2 R_x \Big \}. \nn
\end{align}

\noi By plugging \eqref{inf}, \eqref{er-11} and the value of $ r_1,~t_1 $ into \eqref{er-1}, we arrive at
\begin{align}
	& \P ( ( S_n, V_n^2  ) \in \Omega_1 )  \label{Omega1} \\
	& \leq A_1 R_x \Psi_x^*  \exp \Big\{  - \frac { ( x + c )^2 } { 2 } + \f{2}{5} \gamma^2 x^3 \sum_{i = 1}^n \e [X_i Y_i^2] \nn \\
	& \quad   -  2 x \gamma ( 1 \lor 6 R_x^{1/2} ) \Big[\big( \sqrt{ 1 + x^{-1} ( 1 \lor 6 R_x^{1/2} ) } + 1 \big)^{-1} - \f{2}{5}   \Big] + A_2 R_x \Big\} , \nn
\end{align}
where $ \gamma = ( x + c ) / (2 x ) $. For $ x > 3  $ satisfying \eqref{G1} with a small constant $ c_1 \leq \frac {1} {324}$, it holds that $ x^{-1 } (1 \lor 6 R_x^{1/2} ) \leq 1/3 $. Note also that $ x^3  |\e [ X_i Y_i^2 ] | \leq \frac 23  x^3 L_{3, n} \leq x R_x^{1/2} $ by Lemma \ref{lemma-0}. Hence
\begin{align*}
	& \P ( ( S_n, V_n^2  ) \in \Omega_1  ) \\
	& \leq A_3 x  R_x \Psi_x^* \big[ 1 - \Phi (x + c ) \big] \exp \Big\{ 0.4 \gamma^2 x R_x^{1/2}  -  2 x \gamma  ( 1 \lor 6 R_x^{1/2} ) \cdot 0.06  + A_2 R_x \Big\}.
\end{align*}
Furthermore, observe that $ 1 \lor 6 R_x^{1/2}  \geq 1/2 + 3 R_x^{1/2}$ and $ 2 /5 \leq  \gamma \leq 3/5 $, therefore
\begin{align}
	& \P ( ( S_n , V_n^2 ) \in \Omega_1 ) \label{om1} \\
	& \leq A_3 x R_x \Psi_x^* \big[ 1 - \Phi (x + c) \big] \exp \{  - 0. 024x  -  0.048 x R_x^{1/2} + A_2 R_x \} \nn \\
	& \leq A_4 R_x \Psi_x^*  \big[ 1 - \Phi (x + c) \big] e^{ A_2 R_x }. \nn
\end{align}

\vskip 0.5em 	
As for the second error term, by choosing $ r_2 =  x + c $ and $ t_2 = 2 x (x + c) $, we have under \eqref{G1} that
\begin{equation}
	\inf_{ (u,v) \in \Omega_2 } ( r_2 u - t_2 v ) = r_2 \Big( c + x \sqrt { 1 - x^{-1} (1  \lor   R_x^{1/2} ) } \Big) - t_2 \Big( 1 -  x^{-1} (1  \lor   R_x^{1/2} )  \Big). \label{inf2}
\end{equation}
In the same manner as the proof above, combining (\ref{er-1}), (\ref{om1}) and (\ref{inf2}) with $r_1,t_1$ replaced by $r_2,t_2$, we obtain
\begin{align}
	&\P \( \( S_n , V_n^2 \) \in \Omega_2 \) \label{Omega2} \\
	&\qquad \le A_1 x R_x \[ 1 - \Phi (x + c ) \] \Psi_x^*  \nn \\
	&\hspace{1.5cm} \times \exp \Big\{ - \f x2 \lor 3xR_x^{1/2}   - 6 \ga^2  x^3 \sum_{i=1}^n \e X_i Y_i^2 + A_2 R_x \Big\} \nn \\
	&\qquad \le A_1 R_x [ 1 - \Phi(x + c) ] \Psi_x^* e^{ A_2 R_x }. \nn
\end{align}

\vskip 0.5em
Next we deal with $ \P \( \( S_n , V_n^2 \) \in \Omega_3 \) $. Recall the notations
\begin{gather*}
	\hatX_i = X_i \mathbbm{1} ( | (1 + x ) X_i | \le 1 ) , \quad \hatS_n = \sum_{i=1}^n \hatX_i , \\
	\hatY_i = Y_i \mathbbm{1} ( | (1 + x) Y_i | \le 1 )  ,  \quad  {\hatV_n}^2 = \sum_{i=1}^n \hatY_i^2.
\end{gather*}
We also denote
\begin{gather*}
	\bar{X}_i = X_i \mathbbm{1} ( | (1 + x ) X_i | >  1 ) , \quad \bar{S}_n = \sum_{i=1}^n \bar{X}_i , \\
	\bar{Y}_i = Y_i \mathbbm{1} ( | (1 + x) Y_i | >  1 )  ,  \quad  {\bar{V}_n}^2 = \sum_{i=1}^n \bar{Y}_i^2.
\end{gather*}
It is evident that
\begin{align*}
	& \P \( \( S_n , V_n^2 \) \in \Omega_3 \)   \\
	\le& \P \( \hatS_n > \f{  x V_n + c  }{10}, V_n > \sqrt{ B } \) +
	\P \( \bar{S}_n > \f{ 9\( x V_n + c  \)  }{10}, V_n > \sqrt{ B } \) \\
	\le& K_1 + K_2 + K_3,
\end{align*}
where
\begin{align*}
	K_1 &= \P \( \hatS_n > \f{\sqrt{B}}{10} x + \f {c}{10},\quad \hatV_n^2 > \f B2 \), \\
	K_2 &= \P \( \hatS_n > \f{\sqrt{B}}{10} x + \f {c}{10},\quad \bar{V}_n^2 > \f B2 \), \\
	K_3 &= \P \( \bar{S}_n > \f{ 9\( x V_n + c \) }{10}, V_n > \sqrt{ B } \).
\end{align*}

\noi By Chebyshev's inequality and recalling that  $ B = \max \{ 2500, 200 c_0 \} $ we have
\begin{equation}
	K_1 \le \f{1}{ \( \f B2  - 1 \)^2 } \exp \left\{  - \f 32 x^2 - \f{1}{20} x c \right\} \e \[ \(\hatV_n^2 - 1 \)^2 e^{ \f {x \hatS_n}{2} } \]. \label{K1}
\end{equation}
Let $ \hat{Z}_i = \hatY_i^2 - \e \hatY_i^2 $, then it holds that
\begin{align}
	&  \e \[ \(\hatV_n^2 - 1 \)^2 e^{ \f {x \hatS_n}{2} } \]  \label{err}  \\
	&  \hspace{0.5cm}= \e \[ \(\sum_{ i =1 }^n \hat{Z}_i - \sum_{ i =1 }^n \e \bar{Y}_i^2 \)^2 e^{ \f {x \hatS_n}{2} } \] \nn \\
	&  \hspace{0.5cm}\le  2 \e \[ \(\sum_{ i =1 }^n \hat{Z}_i \)^2 e^{ \f {x \hatS_n}{2} } \] +2 \e \[ \(\sum_{ i =1 }^n \e \bar{Y}_i^2 \)^2 e^{ \f {x \hatS_n}{2} } \]  \nn  \\
	&  \hspace{0.5cm} \le  2 \e \[ \(\sum_{ i =1 }^n \hat{Z}_i \)^2 e^{ \f {x \hatS_n}{2} } \] +2 \e \[ x^{-4}\(\sum_{ i =1 }^n \e x^3 \bar{Y}_i^3 \)^2 e^{ \f {x \hatS_n}{2} } \]  \nn  \\
	&   \hspace{0.5cm}\le 2 \sum_{i=1}^n \f { \e \hat{Z}_i^2 e^{\f{ x \hatX_i }{2}}} { \e e^{\f{ x \hatX_i }{2}} }
	\prod_{ j=1 }^n \e e^{\f{ x \hatX_j }{2}} + 2 x^{-4} R_x^2
	\prod_{i=1}^{n} \e e^{\f{ x \hatX_j }{2}}  \nn   \\
	& \hspace{1cm} + 4 \sum_{ i \ne j }
	\f { \e \hat{Z}_i e^{ \f{ x \hatX_i }{2} } } { \e e^{\f{ x \hatX_i }{2}} }
	\f { \e \hat{Z}_j e^{ \f{ x \hatX_j }{2} } } { \e e^{\f{ x \hatX_j }{2}} } \prod_{ j=1 }^n \e e^{ \f{ x \hatX_j }{2} }. \nn
\end{align}

\noi Furthermore, by Taylor expansion and the same argument used in Lemma \ref{lemma-1}, we can obtain
\begin{align}
	\e e^{\f{ x \hatX_j }{2}} &=
	\exp \Big\{\f18 x^2 \e X_i^2 + \f{1}{48} x^3 \e X_i^3	+ O_1 \de_{x,i} \Big\} ,  \label{0-m}   \\
	\e \hat{Z}_i e^{\f{ x \hatX_j }{2}} &= O \( x \( \e |X_i|^3 + \e |Y_i|^3 \) \)
	+ O \( x^{-2} \de_{x,i} \), \label{1-m}  \\
	\e \hat{Z}_i^2 e^{\f{ x \hatX_j }{2}} &= O \( x^{-4} \de_{x,i} \), \quad {\rm and}\quad \e \bar{Y}_i^2 e^{\f{ x \hatX_j }{2}} = O \( x^{-2} \de_{x,i} \). \label{2-m}
\end{align}
By plugging \eqref{0-m}$-$\eqref{2-m} into \eqref{err}, and the fact that $ x^4 L_{3,n}^2 \le 4\de_x  $, it is readily seen that
\begin{equation}
	\e \[ \(\hatV_n^2 - 1 \)^2 e^{ \f {x \hatS_n}{2} } \] \le \frac{A  R_x}{x^2} \exp \left\{\f18 x^2 + \f{1}{48} x^3 \sum_{i=1}^n \e X_i^3	+ A_1 R_x \right\}. \label{err-bd}
\end{equation}

\noi Substituting \eqref{err-bd} into \eqref{K1} yields for $ |c| \le x/5 $,
\begin{align}
	K_1 &\le C_1 x^{-2} R_x \exp \Big\{ - \f{11}{8} x^2 - \f{xc}{20} + \f{1}{48} x^3 \sum_{i=1}^n \e X_i^3
	+ A_1  R_x  \Big\} \label{bound-K1}  \\
	&\le C_2 x^{-2} R_x \exp \left\{ - \f { (x + c)^2 }{2} - \f 14 x^2 + A_1 R_x \right\} \nn  \\
	&\le C_3 \[ 1 - \Phi( x + c) \] R_x \exp \left\{ - \f 18 x^2 + A_1 R_x  \right\} \nn\\
	&\le C_4 R_x  \[ 1 - \Phi( x + c) \] \Psi_x^* e^{ A_1 R_x} ,
	\nn
\end{align}
and the last inequality holds because $ x $ satisfies \eqref{G1} and
\begin{align*}
	\big| x^3 \sum_{i=1}^n \e X_i^3  \big| \le x^3 L_{3,n}  ,\quad \mb{and} \quad  \big| x^3 \sum_{i=1}^n \e X_i Y_i^2 \big| \le x^3 L_{3,n}.
\end{align*}

\noi Similarly, exploiting the upper bounds \eqref{0-m} and \eqref{2-m}, we have
\begin{equation}
	K_2 < \f 2B \exp \left\{ - \f 32 x^2 - \f{xc}{20} \right\}
	\e \( \bar{V}_n^2 e^{ \f{ x \hat{S}_n }{2}} \),
\end{equation} 	
and
\begin{align*}
	\e \( \bar{V}_n^2 e^{ \f{ x \hat{S}_n }{2}} \) &= \sum_{i=1}^n \f{\e \Big( \bar{Y}_i^2   e^{ \f{ x \hat{X}_i }{2}} \Big) }
	{ \e e^{ \f{ x \hat{X}_i }{2}} }\prod_{ j = 1 }^n \e e^{ \f{ x \hat{X}_j }{2}} \nn  \\
	&\le \f{AR_x}{x^2} \exp \left\{ \f18 x^2 + \f{1}{48} x^3 \sum_{i=1}^n \e X_i^3	+ A_1 R_x \right\}.
\end{align*}
Hence, in the same manner as the proof of \eqref{bound-K1}, it follows that
\begin{equation}
	K_2 \le A R_x \[ 1 - \Phi( x + c) \] \Psi_x^* e^{A_1 R_x}. \label{bound-K2}
\end{equation}

\noi Finally, as for the bound for $ K_3 $, denote
\begin{equation*}
	\bar{X}_i^{(1)}=\bar{X}_i\mathbbm{1}\big( 2xX_i\leq \f{X_i^2}{Y_i^2+c_0\e Y_i^2}   \big):=\bar{X}_i\mathbbm{1}(\mathcal{M}),\quad {\rm and} \quad \bar{X}_i^{(2)}=\bar{X}_i- \bar{X}_i^{(1)}.
\end{equation*}
Cauchy inequality leads to
\begin{align*}
	K_3 &= \P \Big( \sum_{i=1}^n \bar{X}_i  > \f { 9 \( x V_n + c \) }{10}, \;\;  V_n > \sqrt{B} \Big) \nn  \\
	&\leq \P\Big( \sum_{i=1}^n \bar{X}_i^{(1)}  > \f {  \( x V_n + c \) }{100}, \;\;  V_n > \sqrt{B} \Big) \nn  \\
	&\qquad + \P\Big( \sum_{i=1}^n \bar{X}_i^{(2)}  > \f { 89 \( x V_n + c \) }{100}, \;\;  V_n > \sqrt{B} \Big) \nn  \\
	&\le \P\Big( \sum_{i=1}^n 2x\bar{X}_i^{(1)}  > \f {  \( x^2 V_n + cx \) }{50}, \;\;  V_n > \sqrt{B} \Big), \nn \\
	&\qquad + \P \Bigg( \sqrt{ \sum_{i=1}^n \f{ (\bar{X}_i^{(2)})^2 }{ Y_i^2 + c_0 \e Y_i^2 } } > \f{89}{100} \f{x V_n + c }{ \sqrt{ V_n^2 + c_0 } }, \;\;  V_n >  \sqrt{B}   \Bigg) \nn \\
	& := K_4+K_5.
\end{align*}
Recalling that $ B = \max \{ 2500, 200 c_0 \} $. we further have for $ |c| \le x/5 $
\begin{align}
	K_4 &\le \P \(  \sum_{i=1}^n 2x\bar{X}_i^{(1)} >  0.996x^2    \)  \label{bound-K3} \\
	&\le  \f {C_1 x^{-2} }{ \exp \{ 0.986 x^2 \} }  \e \[ \sum_{i=1}^n 2x\bar{X}_i^{(1)}  e^{ 0.99 \sum_{i=1}^n 2x\bar{X}_i^{(1)}  } \]  \nn \\
	&\le  \f { C_2x^{-2} }{ \exp \{ 0.986 x^2 \} }  \sum_{i=1}^n \e \Big[ e^{2x\bar{X}_i^{(1)} }\mathbbm{1}\big( |(1  + x)X_i|>1, \mathcal{M}  \big) \Big]  \prod_{ j \ne i } \e e^{ 2x\bar{X}_j^{(1)} } \nn\\
	&\le A R_x \[ 1 - \Phi (x + c) \] \Psi_x^* \exp{\{ R_x  \} }, \nn
\end{align}
\begin{align}
	K_5 &\le \P \(  \sum_{i=1}^n \f{ (\bar{X}_i^{(2)})^2 }{ Y_i^2 + c_0 \e Y_i^2 }> ( 0.884x )^2   \)  \nn \\
	&\le  \f { C_1x^{-2} }{ \exp \{ 0.773 x^2 \} }  \e \[ \sum_{i=1}^n \f {(\bar{X}_i^{(2)})^2}{Y_i^2 + c_0 \e Y_i^2}  e^{ 0.99 \sum_{i=1}^n \f {(\bar{X}_i^{(2)})^2}{Y_i^2 + c_0 \e Y_i^2}  } \]  \nn \\
	&\le \f { C_2 x^{-2} }{ \exp \{ 0.773 x^2 \} }  \sum_{i=1}^n \e \Big[ e^{\f {(\bar{X}_i^{(2)})^2}{Y_i^2 + c_0 \e Y_i^2}}\mathbbm{1}\big( |x'X_i|>1, \mathcal{M}^c  \big) \Big] \prod_{ j \ne i } \e e^{ \f {(\bar{X}_i^{(2)})^2}{Y_j^2 + c_0 \e Y_j^2}  } \nn\\
	&\le A R_x \[ 1 - \Phi (x + c) \] \Psi_x^* \exp{\{ R_x  \} }, \nn
\end{align}
where we have used the condition $ x L_{3,n} \le c_1 $, for $c_1$ being some small enough constant, and that
\begin{align*}
	\prod_{ j \ne i } \e e^{ 2x \bar{X}_i^{(1)}} &\le  \prod_{ j \ne i } \left\{ 1 + \e e^{ 2x \bar{X}_i^{(1)}} \mathbbm{1} ( |x' X_i | > 1, \mathcal{M} ) \right\}   	\\
	&\le \exp \left\{ \sum_{i=1}^n \e e^{ 2x \bar{X}_i^{(1)} } \mathbbm{1} ( |x' X_i | > 1, \mathcal{M} ) \right\}   \\
	&\le \exp \{ r_x \} \le \exp \{ R_x \},\\
	\prod_{ j \ne i } \e e^{ \f {(\bar{X}_i^{(2)})^2}{Y_j^2 + c_0 \e Y_j^2}} &\le  \prod_{ j \ne i } \left\{ 1 + \e e^{\f {X_i^2}{Y_i^2 + c_0 \e Y_i^2} } \mathbbm{1} ( |x' X_i | > 1, \mathcal{M}^c ) \right\}   	\\
	&\le \exp \left\{ \sum_{i=1}^n \e e^{\f {X_i^2}{Y_i^2 + c_0 \e Y_i^2} } \mathbbm{1} ( |x' X_i | > 1, \mathcal{M}^c ) \right\}   \\
	&\le \exp \{ r_x \} \le \exp \{ R_x \}, \\
	R_x\geq r_x & = \sum_{i=1}^n \e \Big[ e^{2x\bar{X}_i^{(1)} }\mathbbm{1}\big( |x'X_i|>1, \mathcal{M}  \big) \Big]\\
	& \qquad + \sum_{i=1}^n \e \Big[ e^{\f {(\bar{X}_i^{(2)})^2}{Y_i^2 + c_0 \e Y_i^2}}\mathbbm{1}\big( |x'X_i|>1, \mathcal{M}^c  \big) \Big].
\end{align*}

\noi Therefore by \eqref{bound-K1}, \eqref{bound-K2} and \eqref{bound-K3}, we conclude
\begin{equation}
	\P \( \( S_n , V_n^2 \) \in \Omega_3 \) \le A R_x \[ 1 - \Phi( x + c) \] \Psi_x^* e^{ A_1 R_x }.  \label{Omega3}
\end{equation}
Consequently, the desired bound \eqref{prop3-re} follows from \eqref{Omega1}, \eqref{Omega2} and \eqref{Omega3}. The proof is completed.

\subsection{Proof of Proposition \ref{lemma-bigb}} \label{pf_prop6.1}
The main idea is to apply the general Theorem \ref{thm-main}. Recall the notation in \eqref{notationonedependent}, based on the one-dependent structure, $ \{ (X_j, Y_j) \}_{j=1}^k $ is a sequence of independent random vectors. Further, the moment condition in Theorem \ref{thm-2} indicates
\begin{align*}
	\f 12 \Big( \f{\e  S_{n 1} ^2 }{ \e  V_{n 1}^2 }- 1 \Big)   =\rho_n+ O \big(\f{ a^2}{ n^{  \alpha  }} \big).
\end{align*}
Hence
\begin{equation}
	\begin{split}
		& \P \big( S_{n 1 }  \geq  x V_{n 1 } + d_1 n^{ - \alpha / 2  } x B_n  \big) \\
		& =  \P \left (  \frac { S_{n 1 } } { ( \e  S_{ n 1 }^2 )^{ 1/2 } } > \frac {x  \big( 1 + O ( a^2 n^{ - \alpha  } )  \big) } { \sqrt{ 1 + 2 \rho_n} }   \frac {  V_{ n 1 } } {  ( \e    V_{ n 1}^2  )^{ 1/2 }  }  +  O \big( d_1 n^{ -\alpha /2   } x \big) \right) .  
	\end{split}
	\label{normalized} 
\end{equation}

\noi We now focus on bounding the error terms involved in Theorem \ref{thm-1}. Note that  $\{ \xi_i \}_{1\leq i \leq n}$ are one-dependent sequence of random variables, and Rosenthal's inequality yields for $2<p\le 4$,
\begin{align*}
	\e |X_i|^p &= \e \Big[ \Big| \sum_{j\in H_i, odd} \xi_j + \sum_{j\in H_i, even} \xi_j   \Big|^p \Big] \\
	& \leq 2^p \Big\{  \e \Big[ \Big| \sum_{j\in H_i, odd} \xi_j   \Big|^p  \Big]+  \e \Big[ \Big| \sum_{j\in H_i, even} \xi_j   \Big|^p  \Big]  \Big\}  \\
	& \leq A \Bigg[ \Big( \sum_{j\in H_i, odd}\e  \xi_j^2  \Big)^{p/2} +\Big( \sum_{j\in H_i, even}\e  \xi_j^2  \Big)^{p/2}+  \sum_{j\in H_i}\e | \xi_j | ^p     \Bigg] \\
	& \leq A  a_1^p n^{p \alpha /2 }.
\end{align*}
In  addition,
\begin{gather*}
	\e |Y_i|^p \le \big( \e |Y_i|^4 \big)^{\f{p}{4}}    \le \Big[ k\sum_{j\in H_i} \e \xi_j^4  \Big]^{\f{p}{4}} \le  A  a_1^p n^{\f {p \alpha}{2} } , \quad
	B_{n1}^2 := \sum_{j=1}^k \e Y_i^2 \ge  A  a_2^2 n , \\
	\Big| \sum_{j = 1}^k \e  X_j Y_j^2   \Big| + \Big| \sum_{i = 1}^k \e  X_i^3  \Big|\le 7\sum_{i = 1}^k \sum_{j\in H_i} \big| \e \xi_j^3  \big|\leq Aa_1^3 n.
\end{gather*}
Hence
\begin{gather*}
	L_{3, n} \le  A  a^3  n^{- \f {1-\alpha}{2}}, \quad
	\de_x \le  A  a^4  \f {( 1 + x )^4}{ n^{ (1 - \alpha)} } , \quad
	\Psi_{x}^{*} \leq \exp\big(A  a^3  x^3 n^{ - 1/2  } \big).
\end{gather*}
Thefore, applying Theorem \ref{thm-main} to \eqref{normalized} yields
\begin{align}
	& \P \big( S_{n 1 }  \geq  x V_{n 1 } + d_1 n^{ - \alpha / 2  } x B_n  \big)   \label{st1} \\
	&  =  \Big[ 1 - \Phi( \f{x}{ \sqrt{ 1 + 2 \rho_n } } ) \Big]   e^{ O_1  \sum r_{x, j}  } \Big( 1 + O_2 \big(  \frac { a^4 x^4  } { n^{ 1- \alpha } } +   \frac { a^2 x^2  } { n^{ \f{\alpha}{2} } } +\f{a^3x}{n^{\f{1-\alpha}{2}}} \big) \Big)  \nn
\end{align}
for $ x \in ( 2, a^{ - 1 }  \min \{ n^{ \alpha / 4 }, n^{ ( 1 - \alpha  ) / 4  } \} ) $, where
\begin{equation*}
	r_{x,j} \leq \e \bigg[ \exp \bigg\{ \f {  (\sum_{i \in H_j} \xi_i)^2 } { \sum_{i \in H_j} (\xi_i^2 + c_0 \e \xi_i^2) } \bigg\} \IF \Big ( x \Big|\sum_{i \in H_j } \xi_i \Big| >  c(\rho)   a_2 n^{1/2} \Big) \bigg],
\end{equation*}
for some pre-assigned constant $c_0>0$ and constant $c(\rho)$. The derivation of the upper bound of $ r_{x, j} $ is rather complicated.  We exploit a technical lemma \ref{lemma-rxi} whose proof will be provided in Section \ref{sec:lemma-rxi} to bound $ r_{x, j} $. 
\begin{lemma}		\label{lemma-rxi}
	Let $ \{ X_i\}_{1\leq i\leq n}$ be a sequence of independent random variables with $ \e X_i = 0 $ and $ \e |X_i|^p < \infty $ for $ p > 2 $. Denote $ S_n = \sum_{i=1}^{n} X_i $, $ V_n^2 = \sum_{i=1}^{n} X_i^2 $ and $ B_n^2 = \e V_n^2 $. Assume $ b > 0 , A > 0  $. Then there exist positive constants $c_0$ and $ K $ depending on $ p $ and $ A $, such that
	\begin{equation}
		\e  \Big[  e^{ \f {A S_n^2} {V_n^2 + c_0 B_n^2}} \IF ( b |S_n| > 1)  \Big] \le K b^{p} \Big( \sum_{i=1}^{n} \e |X_i|^p + (\e S_n^2)^{\f p2} \Big) .  \label{rxi}
	\end{equation}
\end{lemma}

Since $ \{\xi_i \}_{1 \leq i \leq n }$ are one-dependent and Lemma \ref{lemma-rxi} is built for independent random variables, we separate the odd-indexed terms and the even-indexed terms to apply Lemma \ref{lemma-rxi}. Specifically,
\begin{align}
	r_{x,j} &\le \e \bigg[ \exp \Big\{ \f { 4 |\sum_{i \in H_j, odd} \xi_i |^2} {\sum_{i \in H_j, odd}(\xi_i^2 + c_0\e \xi_i^2) } \Big \} \IF \Big ( \Big|\sum_{i \in H_j,odd } \xi_i \Big| > \f {c(\rho ) a_2 n^{1/2}}{ 2 x } \Big  ) \bigg]  \label{separation} \\
	& \quad + \e \bigg[ \exp \Big\{ \f { 4 |\sum_{i \in H_j, even} \xi_i |^2} {\sum_{i \in H_j, even}(\xi_i^2 + c_0\e \xi_i^2) }  \Big\} \IF \Big ( \Big|\sum_{i \in H_j, even} \xi_i \Big| > \f {c(\rho ) a_2 n^{1/2}}{ 2 x } \Big ) \bigg] \nn\\
	& \le C(\rho) \f { x^4 } { a_2^4  n^2 } \Big[ \sum_{i \in H_j} \e |\xi_i|^4 + ( \sum_{i \in H_j, odd} \e \xi_i^2  )^2 + ( \sum_{i \in H_j, even} \e \xi_i^2  )^2 \Big] \nn\\
	& \le C(\rho)  \f{a^4 x^4 }{ n^{2-2 \alpha} }  , \nn
\end{align}
where the second inequality is derived by applying Lemma \ref{lemma-rxi}. Therefore, $ r_x = \sum_{j=1}^k r_{x,j} \leq  C(\rho) a^4 x^4 / n^{1 - \alpha} $, which together with \eqref{st1} yields the desired result \eqref{big-block}. 

In the same manner as the proof of \eqref{Omega1}, using Chebyshev's inequality and Taylor expansion, it is easy to obtain \eqref{error_low1}. To avoid redundancy, here we omit the proof.

\subsection{Proof of Proposition \ref{lemma-smallb}} \label{pf_prop6.2}
%We apply Bernstein inequality for sum of independent random variables. First for \eqref{error_up},
We start with the error bound (\ref{error_up}). Denote $\omega_j=\hat{\xi}_{j(l+1)}-\e \hat{\xi}_{j(l+1)}$ for $1\leq j\leq k$. Therefore, $\{ \omega_j\}_{1\leq j\leq k}$ are independent random variables with $\omega_j \leq 2\tau$ and $\e\omega_j=0$ for each $j$. Recalling the definitions in \eqref{smallb-def} and $ \tau = B_{n 2 } / x $.
For $ d_1 = \kappa_{\rho} a^2 $ with $\kappa_{\rho} \geq 10$, routine calculation shows
\begin{align*}
	\sum_{j=1}^{k} \e \omega_j^2 &\leq B_{n2}^2 \leq a_1^2 n^{1-\alpha} , \quad d_1 n^{ - \alpha /2  } x B_n -\sum_{j=1}^{k}\e \hat{\xi}_{j(l+1)}  \geq (\kappa_{\rho} - 1) a_1 x n^{\frac{1-\alpha}{2}} .
\end{align*}
%		 \noi Therefore, when $ d_1 >  a^2  +  a  $, we obtain
%		 $$
%		 B_{n2}^2 \leq \tau ( d_1 n^{ - \alpha /2  } x B_n -\sum_{j=1}^{k}\e \hat{\xi}_{j(l+1)}).
%		 $$
%		 Also note that
%		 \begin{align*}
%		 d_1 n^{ - \alpha /2  } x B_n -\sum_{j=1}^{k}\e \hat{\xi}_{j(l+1)}  \leq ( d_1 + 1  ) a_1 x n^{ \frac  { 1 - \alpha  } { 2 }  }.
%		 \end{align*}
Hence, it follows from Bernstein inequality (see for example, Theorem A in \cite{fan2015exponential}) that
\begin{align*}
	\P \big(\hatS_{n2} > d_1 n^{ - \alpha /2  } x B_n \big)
	%		  \\
	%		 =&   \P \Big( \sum_{j=1}^{k}\omega_j \geq d_1 n^{ - \alpha /2  }  x B_n - \sum_{j=1}^{k}\e \hat{\xi}_{j(l+1)} \Big) \\
	&\leq   \P\Big( \sum_{j=1}^{k}\omega_j \geq (\kappa_{\rho}-1)a_1xn^{(1-\alpha)/2} \Big) \\
	\leq \exp\Big( &-\frac{\kappa_{\rho}^2a_1^2 x^2 n^{1-\alpha}} { 3 B_{n2}^2+3\tau\kappa_{\rho} a_1 x n^{\f{(1-\alpha)}{2}} }  \Big)  \leq  \exp\Big( -\frac{\kappa_{\rho} x^2 } { 6  }   \Big).
\end{align*}
%		 Hence for $ d_1 >  8  a^2  $,
%		 \begin{align*}
%		 \P (\hatS_{n2} > d_1 n^{ - \alpha /2  } x B_n + x (\hatV_n - V_{n1}))   \leq \exp \Big\{  -    \frac  { d_1   x^2  }{ 6 a^2 }    \Big \}.
%		 \end{align*}
\noi As for inequality \eqref{error_low2}, similarly,
\begin{align*}
	&\P \Big(\hatS_{n2} < - d_1 n^{ - \alpha /2  } x B_n + x (\hatV_n - V_{n1}), ~V_{n1}^2 > \f14 B_n^2 \Big) \\
	= & \P \Big(\sum_{j=1}^{k} \hat{\xi}_{j(l+1)} < - d_1 n^{ - \alpha /2  } x B_n + x\frac{\sum_{j=1}^{k} \hat{\xi}_{j(l+1)}^2}{\hat{V}_n+V_{n1}}, V_{n1}^2 > \frac{1}{4} B_n^2 \Big) \\
	\leq & \P \Big( \sum_{j=1}^{k} \big( \gamma_j - \e \gamma_j \big) > d_1 n^{ - \alpha /2  }  x B_n- \sum_{j=1}^{k} \e \gamma_j \Big) ,
\end{align*}
where $\gamma_j=\frac{x}{B_n}\hat{\xi}_{j(l+1)}^2 - \hat{\xi}_{j(l+1)}$, $1\leq j\leq k$, are independent random variables with $\gamma_j-\e \gamma_j \leq 4\tau$ for each $j$. Note that
\begin{align*}
	\sum_{j=1}^{k} \e (\gamma_j-\e \gamma_j)^2 &\leq \sum_{j=1}^{k} \e \gamma_j^2 \leq  2\sum_{j=1}^{k} \Big(\frac{x^2\tau^2}{B_n^2} +1\Big)\e \hat{\xi}_{j(l+1)}^2\leq 4a_1^2 n^{1-\alpha} , \\
	\Big |\sum_{j=1}^{k} \e \gamma_j \Big|
	&= \Big| \frac{x}{B_n} \sum_{j=1}^{k} \e \hat{\xi}_{j(l+1)}^2 +\sum_{j=1}^{k} \e \xi_{j(l+1)} \IF (|\xi_{j(l+1)}|>\tau )\Big| \\
	%&\leq \sum_{j=1}^{k} \big|\frac{x}{B_n}\e \hat{\xi}_{j(l+1)}^2 + \frac{1}{\tau} \e \xi_{j(l+1)}^2 I {\{ |\xi_{j(l+1)}|>\tau \}} \big| \\
	&\leq \frac{1}{\tau} \sum_{j=1}^{k} \e \xi_{j(l+1)}^2 = xB_{n2}.
\end{align*}
Again, Bernstein inequality leads to
%		 \begin{equation*}
%		 (d_1 a_2 -a_1)xn^{\frac{1-\alpha}{2}}  \leq d_1 n^{ - \alpha /2  }  x B_n - \sum_{j=1}^{k} \e \gamma_j  \leq  ( d_1 + 1  ) a_1 x n^{ \frac { 1 - \alpha } { 2 } }.
%		 \end{equation*}
%		 Hence by Bernstein inequality and by similar procedure in proof of \eqref{error_up}, we obtain for $ d_1 >     8  a^2   $,
\begin{align*}
	\P \Big(\hatS_{n2} < - d_1 n^{ - \alpha /2  } x B_n + x (\hatV_n - V_{n1}), ~V_{n1}^2 > \f14 B_n^2\Big)   \leq \exp \Big(   -  \frac  { \kappa_{\rho} x^2  }{  14 }     \Big).
\end{align*}
The proof is completed.

\subsection{Proof of Proposition \ref{prop-error}}	\label{pf_prop7}

Observe that for any $  a,b, \delta_2\in\mathbbm{R}$ satisfying $3b\geq a^2$ and $\delta_2\geq 0$, we have for arbitrary $x^2\geq 3$,
\begin{align*}
	x\sqrt{b+\delta_2} &\geq \sqrt{x^2\Big(  \f{a^2}{3}+\delta_2    \Big)}=\sqrt{(x^2-3) \f{a^2}{3}+(x^2-3)\delta_2+ a^2+3\delta_2}\\
	&\geq \sqrt{(x^2-3)\delta_2+ a^2+ 2a\sqrt{\delta_2(x^2-3)}}\geq a+\sqrt{(x^2-3)\delta_2}.
\end{align*}
Denote $H= \{ j ( \ell + 1  ): 1\leq j \leq k  \}$, the index set of small block terms. We set
\begin{align*}
	&a:=\hat{\xi}_{j(\ell+1)}+\xi_{j(\ell+1)-1}+\xi_{j(\ell+1)+1},\quad \hat{S}_n^{(j)}:=\hat{S}_n-a,\\
	&b:=\hat{\xi}_{j(\ell+1)}^2+\xi_{j(\ell+1)-1}^2+\xi_{j(\ell+1)+1}^2,\quad \delta_2:= \big(\hat{V}_n^{(j)}\big)^2 = \hat{V}_n^2-b,
\end{align*}
then $\{ \hat{S}_n>x\hat{V}_n  \}\subset \{  \hat{S}_n^{(j)} > (\sqrt{x^2-3}) \hat{V}_n^{(j)}        \}$ by the above analysis. Moreover, $ \hat{S}_n^{(j)} $ and $ \hat{V}_n^{(j)} $ are independent of $\xi_{j (l +1 )}$ since $ \{\xi_i, 1 \leq i \leq n\}  $ are one-dependent random variables. Thus it follows from \eqref{trun-re} that for $  x \in (c_{\rho} \sqrt{\log n}, d_0 a^{-1} n^{1/8} )$,
\begin{align}
	&\P \Big( \hatS_n \geq x \hatV_n, \max_{1 \le j \le k} |\xi_{j(l+1)}| > \tau \Big)  \label{trun-error1} \\
	& \le \sum_{i \in H} \P \Big( \hatS_n^{ (i ) } \geq \sqrt{x^2-3} \hatV_n^{( i )}  \Big) \P \big(|\xi_{ i }| > \tau\big) \nn \\
	& \le A   \f{\sum_{ i \in H }  \e|\xi_{  i  } |^4}{\tau^4}  \Bigg[ 1 - \Phi\Big(\f{\sqrt{x^2-3}}{ \sqrt{ 1 + 2 \rho_n^{ (i) } } } \Big) \Bigg]   \Big( 1 + O ( \frac {a ^4  x^2  } { n^{ 1/4 } }  ) \Big) \nn \\
	& \le  A  a^4 \f{(1+x)^4} {n^{1/2}}\Bigg[ 1 - \Phi\Big( \f{\sqrt{x^2-3}}{ \sqrt{ 1 + 2 \rho_n } } \Big) \Bigg]  \Big( 1 + O ( \frac {a ^4  x^2  } { n^{ 1/4 } }  ) \Big) \nn \\
	& \le A e^{3/2(1-2\rho)} a^4 \f{(1+x)^4} {n^{1/2}}\Bigg[ 1 - \Phi\Big( \f{x}{ \sqrt{ 1 + 2 \rho_n } } \Big) \Bigg]  \Big( 1 + O ( \frac {a ^4  x^2  } { n^{ 1/4 } }  ) \Big). \nn
\end{align}
where $\rho_{ n }^{(i)}$ is defined parallel to $\rho_n$ and
\begin{align*}
	\rho_{ n }^{(i)}&:=  \frac 1 2 \Big(  \frac { \e [ ( \hat{S}_n^{(i)})^2 ]   } { \e [ ( \hat{V}_n^{(i)} )^2 ] }  - 1 \Big) = \rho_n \Big(1+O\big( \f{ a^2}{n} \big)  \Big)
\end{align*}
because $ \e \xi_i^4 \leq a_1^4  $ and $ \e \xi_i^2 \geq a_2^2$ for $ i \geq 1 $. Then \eqref{error-bound1} is derived.

Regarding the proof of \eqref{error-bound2}, a more technical iterative argument will be exploited. We have
\begin{align}
	&\P \Big( S_n \geq x V_n, \max_{1 \le j \le k} |\xi_{j(l+1)}| > \tau \Big)  \nn\\
	& \le \sum_{i_1 \in H} \P \Big(S_n^{(i_1)} \geq \sqrt{x^2-3} V_n^{(i_1)}\Big) \P \big(|\xi_{i_1}| > \tau\big) \nn\\
	& \le \sum_{i_1 \in H} A \f{a^4 (1+x)^4} {n} \P \Big(S_n^{(i_1)} \geq \sqrt{x^2-3} V_n^{(i_1)}\Big) \nn \\
	& \leq \sum_{i_1\in H} A \f{a^4 (1+x)^4} {n} \Big[ \P \Big(\hat{S}_n^{(i_1)} \geq \sqrt{x^2-3} \hat{V}_n^{(i_1)}\Big) \nn\\
	&\qquad \qquad \qquad +\P \Big(S_n^{(i_1)} \geq \sqrt{x^2-3} V_n^{(i_1)} , \max_{j\in H\setminus \{ i_1\}} |\xi_{j(l+1)}| > \tau \Big) \Big] \nn \\
	& \leq \sum_{i_1\in H} A \f{a^4 (1+x)^4} {n} \Big[ \P \Big(\hat{S}_n^{(i_1)} \geq \sqrt{x^2-3} \hat{V}_n^{(i_1)}\Big) \nn\\
	&\qquad \qquad \qquad +\sum_{ i_2\in H \setminus \{i_1 \} }\P \Big(S_n^{(i_1,i_2)} \geq \sqrt{x^2-6} V_n^{(i_1,i_2)} \Big) \P \big(|\xi_{ i_2} | > \tau \Big) \Big] \nn \\
	& \leq  \sum_{i_1\in H} A \f{a^4 (1+x)^4} {n} \P \Big(\hat{S}_n^{(i_1)} \geq \sqrt{x^2-3} \hat{V}_n^{(i_1)}\Big) \nn\\
	& \qquad +\sum_{i_1\in H} \sum_{ i_2\in H \setminus \{i_1 \} } \Big( A \f{a^4 (1+x)^4} {n} \Big)^2  \P \Big(S_n^{(i_1,i_2)} \geq \sqrt{x^2-6} V_n^{(i_1,i_2)} \Big). \label{trun-error2}
\end{align}
Put $ u = [\f{x^2}{6}] $ so that $\sqrt{x^2-3u}  \approx  x / \sqrt 2$ and we repeat the above procedure iteratively up to $ u $ times,
\begin{align*}
	&\P \Big( S_n \geq x V_n, \max_{1 \le j \le k} |\xi_{j(l+1)}| > \tau \Big)  \nn\\
	& \leq  \sum_{j=1}^{u-1} \sum_{i_1 \in H  } \cdots \sum_{i_j \in H  }  \Big( A \f{a^4 (1+x)^4} {n} \Big)^j \P \Big(\hatS_n^{(i_1,\ldots,i_j)} \geq \sqrt{x^2-3j} \hatV_n^{(i_1,\ldots,i_j)} \Big)  \nn\\
	&~ \quad +\sum_{i_1 \in H  } \cdots \sum_{i_u  \in H  } \Big( A \f{a^4 (1+x)^4} {n}\Big)^u \P \Big(S_n^{ (  i_1,\ldots,i_u ) } \geq \sqrt{x^2-3u} V_n^{ ( i_1,\ldots,i_u ) }\Big).
\end{align*}
As we have chosen $ \alpha = 1/2 $ and $|H|=k\asymp n^{1-\alpha}=\sqrt{n}$, it follows by \eqref{trun-re} that for arbitrary $i_1,\cdots,i_j\in H$, and $1\leq j\leq x^2/6$, we have for $ x \in (c_{\rho } \sqrt{\log n } , d_0 a^{- 2} n^{1/8}) $,
\begin{align*}
	& \sum_{i_1 \in H  } \cdots \sum_{i_j \in H  }   \Big( A \f{a^4 (1+x)^4} {n} \Big)^j \P \Big(\hatS_n^{(i_1,\ldots,i_j)} \geq \sqrt{x^2-3j} \hatV_n^{(i_1,\ldots,i_j)} \Big) \\
	& \leq  \Big( A \f{a^4 (1+x)^4} {n^{1/2}} \Big)^j  \bigg[ 1-\Phi\bigg(  \f{\sqrt{x^2-3j}}{\sqrt{1+2\rho_{ n }^{(i_1,\cdots,i_j)}}}  \bigg) \bigg]     \Big( 1+O\Big(  \f{a^4 x^2}{n^{1/4}}   \Big)  \Big) \\
	& \leq   A \Big( A \f{a^4 (1+x)^4} {n^{1/2}} \Big)^{j} e^{3j/(2+4\rho_{ n })} \bigg[ 1-\Phi\Big(  \f{x}{\sqrt{1+2\rho_{ n }}}  \Big) \bigg]     \Big( 1+O\Big(  \f{a^4 x^2}{n^{1/4}}   \Big)  \Big) \\
	& \leq   A \f{a^4(1+x)^4} {n^{1/2}} \Big( A \f{a^4 (1+x)^4} {n^{1/2}}e^{3/(2+4\rho_{ n })} \Big)^{j-1}  \bigg[ 1-\Phi\Big(  \f{x}{\sqrt{1+2\rho_{ n }}}  \Big) \bigg]    \\
	& \leq   A \f{a^4(1+x)^4} {n^{1/2}}\cdot \big( \f{1}{2} \big)^{j-1} \bigg[ 1-\Phi\Big(  \f{x}{\sqrt{1+2\rho_{ n }}}  \Big) \bigg]
\end{align*}
for sufficiently large $n$ and some constant $A''$ depending on $a,\rho$, where $\rho_{ n }^{(i_1,\cdots,i_j)}$ is defined as
\begin{align*}
	\rho_{ n }^{(i_1,\cdots,i_j)} : =  \frac 1 2 \Big(  \frac { \e [ (\hat{S}_n^{( i_1, \ldots, i_j )})^2  ]  } { \e [ ( \hat{V}_n^{(i_1, \ldots, i_j )} )^2 ] }  - 1 \Big) = \rho_n \Big(1+O\big( \f{ a^2x^2}{n} \big)  \Big)
\end{align*}
Therefore, by combining the upper bounds above we conclude	
\begin{align*}
	&\P \Big( S_n \geq x V_n, \max_{1 \le j \le k} |\xi_{j(l+1)}| > \tau \Big)  \\
	& \le  \sum_{j=1}^{u-1} A\f{a^4(1+x)^4} {n^{1/2}}\cdot \big( \f{1}{2} \big)^{j-1} \bigg[ 1-\Phi\Big(  \f{x}{\sqrt{1+2\rho_{ n }}}  \Big) \bigg] + \Big(A \f{ a^ 4 (1+x)^4}{n^{1/2}}\Big)^{u}  \nn \\
	& \le  A_1 \f{a^4 x^4}{n^{1/2}}\bigg[ 1-\Phi\Big(  \f{x}{\sqrt{1+2\rho_{ n }}}  \Big) \bigg]\Big( 1 + \Big(A \f{ a^ 4 (1+x)^4}{n^{1/2}} e^{3/(1+2\rho_n)} \Big)^{x^2/3} \Big)\nn \\
	& \le A_2  \f{a^4 x^4}{n^{1/2}}\bigg[ 1-\Phi\Big(  \f{x}{\sqrt{1+2\rho_{ n }}}  \Big) \bigg]  \nn
\end{align*}
for $ x \in (c_{\rho } \sqrt{\log n } , d_0 a^{- 2} n^{1/8}) $. The proof is completed.

\subsection{Proof of Proposition \ref{prop-thm3.2}} \label{pf-thm3.2-prop}
We apply Theorem \ref{thm-1} to estimate $ \P ( S_{k,1} \ge x V_{k,1} - \vp x B_n ) $. As a preparation, we first calculate the relevant moments. By Lemma \ref{lemma-mo1}, there exist some positive constants $ A_1 $,$ A_2 $ depending on $ a_1, a_2 $ and $ \tau $ such that for $ 1 \leq u \leq k_1  $,
\begin{equation}
	\e \xi_u^4 +\e \eta_u^4 \le A_1 n^{2 (\alpha +\alpha_1) } \mu_1^4,\quad  \e |\xi_u |^3 +\e |\eta_u |^3 \le A_2 n^{\f32 (\alpha +\alpha_1 ) } \mu_1^3. \label{4-m}
\end{equation}
In addition, by Lemma \ref{lemma-mo2} and some routine calculations we can obtain
\begin{align}
	|\e \xi_u^3|  \le  A_3 n^{(\alpha + \alpha_1)} c_1^3 , \quad  |\e \xi_u \eta_u^2|  \le  A_3 n^{(\alpha + \alpha_1)} c_1^3,  \label{3-mm}
\end{align}
for a positive constant $ A_3 $ depending on $ a_1, a_2 $ and $ \tau $.
To apply the moderate deviation theorem for general self-normalized sums, it remains to estimate the error term
\begin{equation*}
	r_{x, u } :=  \e \bigg[ \exp \Big( \f{  \xi_u ^2 } { \eta_u ^2 + c_0 \e \eta_u ^2 } \Big) \IF \Big( \f{1 + x} {B_{n,1}} |\xi_u| > 1 \Big) \bigg] .
\end{equation*}
We split the block $ I_u  $ into odd-indexed terms and even-indexed terms to construct a weakly dependent structure so that Lemma \ref{lemma-indep} can be applied. To be specific, let  $ \mathcal{J}_1 = \{  j  : j~ \mbox{ is odd} \} $ and $ \mathcal{J}_2 = \{  j  : j ~ \mbox{is even} \} $.
%	  Let $ v = [ \frac { 1 - \alpha } { 2 \alpha \tau  } ] + 3 $, and for $ s = 0, 1, \ldots, v - 1 $, define
%	 \begin{equation*}
%	     \mathcal{J}_s:= \{  j: j \equiv s \;(mod\; v)        \},
%	 \end{equation*}
Correspondingly, we define for $ t = 1, 2 $
\begin{equation*}
	\mathcal{G}_t :=\sum_{j\in I_u \cap \mathcal{J}_t} Y_j, \quad \mathcal{K}_t^2 := \sum_{j\in I_u \cap \mathcal{J}_t} \big( Y_j^2+c_0 \e Y_j^2\big).
\end{equation*}
In the same manner as \eqref{separation}, we have
%and (\ref{e-1}), we divide $r_{x, u }$ into $v$ parts,
\begin{align}
	r_{x, u }
	%	 &=   \int_{0}^{\infty} e^t \; \P \Big( \f { |\xi_u| } { \sqrt{ \eta_u^2  + c_0 \e \eta_u^2 } } > \sqrt{t}, x'|\xi_u| > B_{n,1} \Big) \; d t            \label{rxl}  \\
	%	 &=  \sum_{s=0}^v \int_{0}^{\infty} e^t \; \P \Bigg( \f { |\sum_{s=0}^v p_s| } { \sqrt{ \sum_{s=0}^v q_s^2 } } > \sqrt{t}, \, \f{x'|\xi_u|}{B_{n,1}} > 1,\, | p_s |\geq \f{1}{v}|\xi_u| \Bigg) \; d t            \nn  \\
	%	 &\leq  \sum_{s=0}^v \int_{0}^{\infty} e^t \; \P \Bigg( \f { |p_s| } { \sqrt{ \sum_{s=0}^v q_s^2 } } > \f{\sqrt{t}}{v}, \, \f{x'|p_s|}{B_{n,1}} > \f{1}{v} \Bigg) \; d t            \nn  \\
	%	 &\leq  \sum_{s=0}^v \int_{0}^{\infty} e^t \; \P \Bigg( \f { |p_s| } { q_s  } > \f{\sqrt{t}}{v}, \, \f{x'|p_s|}{B_{n,1}} > \f{1}{v} \Bigg) \; d t            \nn  \\
	&\le   \sum_{ t = 1}^{ 2 } \e \Big[ \exp \Big( \f{  4  \mathcal{G}_t^2 } { \mathcal{K}_t^2 } \Big) \IF \Big( \f{ (1 + x)| \mathcal{G}_t |} {B_{n,1}}  >  \f{1}{v} \Big) \Big] :  = \sum_{ t = 1}^{ 2 } r_{ x, u, t }. \label{rxl}
\end{align}
As $\{ r_{x, u, t}, t = 1, 2 \} $ share the same bound, we only present the analysis of $ r_{x, u ,1} $. Suppose the cardinality of $ I_{u}\cap \mathcal{J}_1 $ is $  [n^{\alpha_1}  / 2] $ for simplicity. By Lemma \ref{lemma-indep}, $ \{Y_j, j \in I_u \cap \mathcal{J}_1 \} $ can be approximated by a sequence of independent random variables $ \{Y_j^*, 1 \le j \le  [n^{\alpha_1}  /2] \} $. Denote	
\begin{equation*}
	\mathcal{G}_1^*=\sum_{j=1}^{[n^{\alpha_1} / 2 ]} Y_j^*, \quad  (\mathcal{K}_1^{*} )^2 =  \sum_{j=1}^{[n^{\alpha_1} / 2]} \big( (Y_j^*)^2+c_0 \e (Y_j^*)^2\big).
\end{equation*}
It follows from Lemmas \ref{lemma-rxi} and \ref{lemma-indep} that
%	 Then it follows from Lemma \ref{lemma-rxi} and the choice of $v  =  [ \frac { 1 - \alpha } { 2 \alpha \tau  } ] + 3  $ that
%	 Therefore, by applying Lemma \ref{lemma-indep} and Theorem \ref{thm-1}, and noticing the fact of \eqref{mo-ratio}, \eqref{4-m}--\eqref{3-mm} and \eqref{error_ind}, we obtain there exist constants $ d_1 $ and $ c_0 $  depending on  $ c_1/c_2, a_1, a_2, \alpha $ and $\tau $ such that when $ \vp = d_1  n^{ - \alpha_1 / 2  } \log n   $,
%	 uniformly for $ x \in (0, c_0 \min \{ n^{( 1 - \alpha - \alpha_1 )/4 }, n^{\alpha \tau/2}, n^{ \alpha / 2  }, ( \log n  )^{ - 1/2  } n^{ \alpha_1 / 4  }  \}  )$.
%	 Thus \eqref{main1}--\eqref{main2} combined with \eqref{er1}-\eqref{er2} yield for $ \alpha_1 = ( 1 - \alpha  ) /2 $,
%	 uniformly for $ x \in (0, c_0 \min \{   n^{\alpha \tau/2}, n^{ \alpha / 2  }, ( \log n  )^{ - 1/2  } n^{ ( 1 - \alpha ) / 8  }  \}  )$.
%	
%	 Finally, the procedure of handling the truncated errors in \eqref{trun1} and \eqref{trun2} is the same as that of \eqref{trun-error1} and \eqref{trun-error2}, so we omit the proof. Consequently, the same bound in \eqref{trunc-re} holds for $ \P(S_k \ge V_k) $. The proof is completed. \qed
%	
\begin{align}
	r_{x,u,1} &\le \e \Big[ \exp  \Big( \f{ 4 ( \mathcal{G}_1^* )^2 } {(\mathcal{K}_1^{*})^2 } \Big) \IF \Big(\f{ (1 + x )| \mathcal{G}_1^*  | } {B_{n,1}}  > \f {1} {2}  \Big) \Big] \label{rxl_1} \\
	& \quad + \e \Big[ \exp \Big( \f{  4 \mathcal{G}_1^2 } {\mathcal{K}_1^2 } \Big) \IF \Big( \f{ (1 + x )| \mathcal{G}_1|} {B_{n,1}}  >  \f{1}{2} \Big) \IF \big( \, \exists j: Y_j \neq Y_j^*  \big) \Big] \nn \\
	&\le \e \Big[ \exp  \Big( \f{ 4 ( \mathcal{G}_1^* )^2 } {(\mathcal{K}_1^{*})^2 } \Big) \IF \Big(\f{ (1 + x) | \mathcal{G}_1^*  | } {B_{n,1}}  > \f12  \Big) \Big] +  \f{a_1 n^{ \alpha_1 } }{2} e^{ 4 n^{ \alpha_1 }- a_2  n^{\alpha  \tau}  } \nn \\
	& \le \f{ C(1 + x)^4 }{ n^{2(1-\alpha-\alpha_1)} } +  \f{a_1 n^{ \alpha_1 } }{ 2 } e^{ 4 n^{ \alpha_1 }- a_2 n^{ \alpha \tau} } , \nn
\end{align}
%\begin{equation}
%\e \Big[ \exp  \Big\{ \f{4 ( \sum_{j=1}^p Y_j^* )^2 } {\sum_{j=1}^p ((Y_j^*)^2 + c_0 \e (Y_j^*)^2  )  }  \Big\} I \{ \f{1 + x} {B_{n,1}} | \sum_{j=1}^p Y_j^*  | > \f12  \} \Big] \le A_4  \f{ (1 + x)^4 }{ n^{2(1-\alpha-\alpha_1)} } \;, \label{error_ind}
%\end{equation}
where  $ \alpha_1 $ should be chosen such that $  \alpha_1 \leq \alpha \tau $. Therefore, by Theorem \ref{thm-1} and plugging in the bounds given in (\ref{mo-ratio}), (\ref{4-m}), (\ref{3-mm}) and \eqref{rxl_1}, we can obtain there exist a constant  $ c_0 $ depending on  $d_1,  \mu_1/\mu_2, a_1, a_2, \alpha $ and $\tau $ such that when $\epsilon=d_1 n^{-\alpha_1/2}\log n$,
\begin{align}
	&\P \big( S_{k,1} \ge x V_{k,1} \pm \vp x B_n\big) \big/ \big[1- \Phi(x) \big]    \\
	=& 1 + O\bigg(\f{(1 + x)^4}{ n^{1- \alpha- \alpha_1} }  + \f{1 + x }{n^{(1 - \alpha - \alpha_1 )/2 }} + \f{( 1 + x )^2 }{ n^{\alpha }  } + \f{ ( 1 + x )^2 \log n } {n^{ \alpha_1/2 }} \bigg)  \nn
	%& \P( S_{k,1} \ge x V_{k,1} + \vp x B_n) \nn \\
	%&=  [1- \Phi(x)] \big(1 + O(1) (\f{(1 + x)^4}{ n^{1- \alpha- \alpha_1} } + \f{(1 + x )^3}{\sqrt{n}} + \f{1 + x }{n^{(1 - \alpha - \alpha_1 )/2 }} + \f{( 1 + x )^2 }{ n^{\alpha }  } + \f{ ( 1 + x )^2 } {n^{ \alpha_1/2 }} ) \big)   \label{main2}
\end{align}
uniformly for $ x \in (0, c_0 \min \{ n^{( 1 - \alpha - \alpha_1 )/4 }, n^{\alpha \tau/2}, n^{\alpha/2}, (\log n)^{-1/2}n^{\alpha_1/4}\} )$. The proof of Proposition \ref{prop-thm3.2} is completed.

\subsection{Proof of Proposition \ref{prop_approx}}    \label{pf-prop8}

\noi Note that $ \lim_{\ell \rightarrow - \infty} \e ( X_{i} \vert \mathscr{F}_{\ell} ) $ belongs to the tail $ \sigma $-field, hence $  \lim_{\ell \rightarrow - \infty} \e ( X_{i} \vert \mathscr{F}_{\ell} ) = \e (  X_{i} ) = 0 $ by the zero-one law. Therefore, we may write $X_i$ as
\begin{equation*}
	X_i = \sum_{u=0}^{\infty} \mathcal{P}_{i-u} X_i , \qquad  {\rm where}~\ \mathcal{P}_{ i - u } ( \cdot ) =\e(\cdot|\mathscr{F}_{ i - u })-\e(\cdot|\mathscr{F}_{ i - u -1}).
\end{equation*}
Observe that,
\begin{equation*}
	\parallel Y_j \parallel_4 = \Big \Vert  \sum_{u =0}^{\infty}\sum_{ i \in H_j}\mathcal{P}_{i -u }X_i  \Big \Vert_4 \leq \sum_{ u =0}^{\infty}  \Big\Vert \sum_{ i \in H_j}\mathcal{P}_{ i -u } X_i  \Big \Vert_4.
\end{equation*}
For fixed $ u \geq 0 $, denote $Z_i=\mathcal{P}_{ i- u }X_i $ and $ \mathscr{A}_i   =\mathscr{F}_{ i-u }$. It is obvious that $ Z_i \in \mathscr{A}_i $ and $ \mathscr{A}_{i - 1 } \subset  \mathscr{A}_i  $ for any $ i \ge 1   $.
Moreover,
\begin{align*}
	\e ( Z_i |\mathscr{A}_{ i - 1 } ) &  =   \e ( X_i | \mathscr{A}_{i - 1 } )  -   \e ( X_i | \mathscr{A}_{i - 1 } )    = 0,
\end{align*}
thus for any fixed $ u \geq 0  $, $\left\{ ( Z_i ,\mathscr{A}_i ), i \geq 1  \right\} $ is a martingale difference sequence. It follows from Burkholder's martingale inequality that
\begin{align*}
	\e \Big[  \Big|  \sum_{ i \in H_j  }  Z_i  \Big|^4  \Big]  \leq C  \e \Big[ \Big|  \sum_{ i \in H_j  } Z_i^2  \Big|^{  2  }  \Big]  \leq C m     \sum_{ i \in H_j } \e ( Z_i^4 ).
\end{align*}
%	\begin{equation*}
%	\parallel \sum_{l\in H_j}\mathcal{P}_{l-i}X_l \parallel_r=\parallel \sum_{k=1}^{m}Z_k \parallel_r \leq C_r \cdot \left( \sum_{k=1}^{m} \parallel Z_k \parallel_r^2 \right)^{1/2} \quad {\rm with}\; C_r=\f{18r}{\sqrt{1-\frac{1}{r}}}<\infty .
%	\end{equation*}
\noi In addition, by conditional Jensen's inequality and recalling the definition in \eqref{theta_r}, we obtain
\begin{align*}
	\e  ( Z_i^4  ) &  =  \e \Big[    \Big| \e\big( X_i  \big|\mathscr{F}_{ i - u } \big)-\e \big( X_ i \big|\mathscr{F}_{ i - u -1} \big)  \Big| ^4 \Big] \\
	& =  \e\Big\{ \Big|  \e\big[ X_i - G_i \big( \ldots,\vp_{ i - u - 1 },\epsilon_{  i  - u }^{\ast},\vp_{ i - u +1},\ldots,\vp_{i }\big) \big| \mathscr{F}_{ i - u } \big] \Big|^4  \Big\} \\
	& \leq    \e \big[  \big( X_i - G_i \big( \ldots,\vp_{ i - u - 1 },\epsilon_{  i  - u }^{\ast},\vp_{ i - u +1},\ldots,\vp_{i }\big)  \big)^4    \big ]   = \big [ \theta_4 (u) \big]^4 .
\end{align*}
Thus, we obtain by condition \eqref{GMC} and Remark \ref{remark-theta} that
\begin{align*}
	\parallel Y_j \parallel_4  &  \leq \sum_{ u =0}^{\infty}  \Big\Vert \sum_{ i \in H_j} Z_i \Big \Vert_4 \leq C  \sum_{ u = 0 }^{\infty} m^{ 1/ 2 }  a_1' e^{ - a_2' u^{ \tau  } }  \leq C_1 m^{ 1/2 }
\end{align*}
for $ C_1  $ being a constant depending on $ a_1, a_2 $ and $ \tau $. Further, it follows from Minkowski's inequality and \eqref{difference} that for sufficiently large $n$,
\begin{gather}
	\parallel \tilde{Y}_j \parallel_4  \leq  \parallel Y_j \parallel_4 + \parallel \tilde{Y}_j-Y_j \parallel_4 \leq C_1 m^{1/2}+ m a_1 e^{-a_2m^{\tau}} \leq \tilde{C}_1 m^{1/2} ,  \label{dd1}\\
	\parallel \tilde{Y}_j \parallel_2  \geq \parallel Y_j \parallel_2 - \parallel \tilde{Y}_j-Y_j \parallel_2 \geq c_1 m^{1/2} - ma_1e^{-a_2m^{\tau}} \geq \tilde{C}_2 m^{1/2} , \label{dd2}
\end{gather}
where $\tilde{C}_1$ and $\tilde{C}_2=\omega_1-n^{\f{\alpha}{2}}a_1e^{-a_2n^{\alpha\tau}}> \omega_1/2 $ are constants depending on $\tau, \alpha, \omega_1,  a_1$ and $a_2$.
\vskip 1em

Denote $\rho_k =  \big( \sum_{j =1}^{k-1}\e( \tilde{Y}_j \tilde{Y}_{j +1} ) \big) \big/ \big( \sum_{j =1}^{k}\e (  \tilde{Y}_j ^2 )   \big) $. Note that
\begin{equation*}
	\tilde{C}_2^2\cdot n =\tilde{C}_2^2\cdot  km \leq \sum_{i=1}^{k}\e(\tilde{Y}_i^2) \leq \tilde{C}_1^2 \cdot km = \tilde{C}_1^2 \cdot n,
\end{equation*}
and
\begin{align*}
	& \bigg|\sum_{ j =1}^{k-1}\e \Big(\tilde{Y}_j \tilde{Y}_{ j +1}\Big) -\sum_{ j =1}^{k-1}\e\Big(Y_j Y_{j +1}\Big)  \bigg|\\
	& \hspace{0.8cm} \leq  \sum_{ j =1}^{k-1}\bigg|\e\Big[Y_ j \Big(\tilde{Y}_{ j +1}-Y_{ j +1}\Big)\Big] +\e\Big[Y_{ j +1}\Big(\tilde{Y}_{ j }-Y_{ j }\Big)\Big]  \\
	&\hspace{3cm}  +\e\Big[\Big(\tilde{Y}_j -Y_j \Big)\Big(\tilde{Y}_{j +1}-Y_{j +1}\Big) \Big] \bigg| \\
	& \hspace{0.8cm} \leq \sum_{j =1}^{k-1}   \Big( \big\Vert Y_j \big\Vert_2 \cdot  \big\Vert \tilde{Y}_{j +1}-Y_{ j +1}\big\Vert_2 + \big\Vert Y_{ j +1}\big\Vert_2 \cdot  \big\Vert \tilde{Y}_{ j }-Y_{j }\big\Vert_2  \\
	& \hspace{3cm} +  \big\Vert \tilde{Y}_j -Y_j \big\Vert_2 \cdot  \big\Vert \tilde{Y}_{j +1}-Y_{j +1}\big\Vert_2   \Big)  \\
	& \hspace{0.8cm} \leq  k \left( 2 \tilde{C}_1 m^{1/2} \cdot ma_1e^{-a_2m^{\tau}} + m^2a_1^2e^{-2a_2m^{\tau}}   \right).
\end{align*}
Immediately we have
\begin{equation}
	\bigg|\sum_{i=1}^{k-1}\e (\tilde{Y}_i\tilde{Y}_{i+1}  ) -\sum_{i=1}^{k-1}\e ( Y_iY_{i+1} )  \bigg| \bigg/ \bigg( \sum_{i=1}^{k}\e ( \tilde{Y}_i^2 )  \bigg)\leq C_3  e^{-C_4  a_2m^{\tau}} , \label{dif-1}
\end{equation}
with constants $C_3$ and $C_4$ depending on $\tau, \alpha, \omega_1,  a_1$ and $a_2$. Now, for further investigation into the ratio of $ \sum_{ j = 1 }^{ k - 1 } \e ( Y_j Y_{j + 1} ) $ against $ \sum_{j =1 }^{ k - 1 } \e   (\tilde{Y}_i^2 ) $, we make use of the martingale structure again by representing $Y_{j}$ and $Y_{j+1}$ as
\begin{equation*}
	Y_j = \sum_{ i \in H_j }\sum_{u =0}^{\infty}\mathcal{P}_{i - u } X_i , \qquad  Y_{j+1}=\sum_{ i \in H_{j+1}}\sum_{ u =0}^{\infty}\mathcal{P}_{ i - u } X_i .
\end{equation*}
Thus
\begin{equation*}
	\e(Y_jY_{j+1})= \sum_{ i_1 \in H_j } \sum_{ i_2 \in H_{j+1} } \sum_{u_1=0}^{ \infty}\sum_{ u_2=0}^{ \infty} \e \Big[ \Big( \mathcal{P}_{ i_1 - u_1 } X_{i_1} \Big) \cdot \Big( \mathcal{P}_{i_2- u_2 } X_{i_2} \Big)  \Big].
\end{equation*}
Note that $\mathcal{P}_{i_2-u_2}X_{i_2} \in \mathscr{F}_{i_2-u_2}$ and for any $\ell_1 < i_1 - u_1$, $\ell_2 < i_2-u_2$,
\begin{equation*}
	\e\Big[ \mathcal{P}_{i_1-u_1}X_{i_1}\big|\mathscr{F}_{\ell_1}  \Big] = 0 , \quad \e\Big[ \mathcal{P}_{i_2-u_2}X_{i_2}\big|\mathscr{F}_{\ell_2}  \Big] = 0.
\end{equation*}
Therefore, for $ i_1- u_1\neq i_2- u_2$, we have $ \e \left[ \left( \mathcal{P}_{ i_1 - u_1 } X_{i_1} \right) \cdot \left( \mathcal{P}_{i_2- u_2 } X_{i_2} \right)  \right] = 0 $. Consequently,
\begin{align*}
	\big|\e(Y_jY_{j+1}) \big| & = \bigg| \sum_{ i_1\in H_j } \sum_{ i_2\in H_{j+1} } \sum_{u_1=0}^{ \infty}  \e \Big[ \Big( \mathcal{P}_{i_1-u_1}X_{i_1} \Big)\cdot \Big( \mathcal{P}_{i_1-u_1} X_{i_2} \Big)  \Big] \bigg| \\
	&\leq \sum_{i_1\in H_j}\sum_{i_2\in H_{j+1}} \sum_{u_1=0}^{+\infty} 4a_1'^2\cdot \exp\Big[ -a_2' \Big(u_1^{\tau}+(i_2-i_1+u_1)^{\tau} \Big)  \Big],
\end{align*}
where we have used the fact that $ \Vert   \mathcal{P}_{ i_1 - u_1 } X_{i_1} \Vert_2  \leq \theta_4 ( u_1 ) $ and $ \Vert   \mathcal{P}_{ i_1 - u_1 } X_{i_2} \Vert_2  \leq \theta_4 ( i_2 - i_1 + u_1 ) $.
%	As a matter of fact, we would have
%	\begin{align*}
%	\big|\e(Y_jY_{j+1})  \big|\geq \sum_{l_1\in H_j}\sum_{l_2\in H_{j+1}} \sum_{i_1=0}^{+\infty} a_1^2\cdot \exp\left\{ -a_2[i_1^{\tau}+(l_2-l_1+i_1)^{\tau} ]  \right\}
%	\end{align*}
%	if $X_l=\sum_{i=-\infty}^{l}a_1e^{-a_2(l-i)^{\tau}}\epsilon_{l-i}$.

It is trivial to see from the integrability of the function $ e^{ - a x^{\tau} } $ that
\begin{align}
	&\sum_{ i_1\in H_j}\sum_{ i_2\in H_{j+1}} \sum_{u_1=0}^{ \infty} 4 a_1'^2\cdot \exp\Big[ -a_2' \Big(u_1^{\tau}+(i_2-i_1+u_1)^{\tau} \Big)  \Big] \label{dif-2} \\
	\leq& \sum_{ u =0}^{ \infty}  \sum_{ \ell =1 }^{2 m}  \ell \cdot  \exp\Big[ - a_2' \Big( u ^{\tau}+( \ell + u )^{\tau} \Big) \Big]  \nn   \\
	\leq&  \left( \sum_{ \ell =1}^{ 2 m }  \ell  e^{-a_2' \ell^{\tau}} \right)  \left(\sum_{ u =0}^{\infty} e^{-  a_2'  u ^{\tau}} \right) \leq C_5  ,    \nn
\end{align}
for some constant $ C_5 $ depending on $  \alpha, a_1, a_2 $ and $ \tau $. Finally, as a result of \eqref{dd2}, \eqref{dif-1} and \eqref{dif-2}, we obtain
\begin{equation}
	| \rho_k |:= \bigg|\sum_{i=1}^{k-1}\e (\tilde{Y}_i\tilde{Y}_{i+1}  )  \bigg| \bigg/ \bigg( \sum_{i=1}^{k}\e ( \tilde{Y}_i^2 )  \bigg) \leq  C_6 m^{ -1  }  \label{rho-k}
\end{equation}
for large enough $n$, where $C_6$ is a constant depending on $\tau, \alpha, \omega_1,   a_1 $ and $a_2$.

Now we are ready to apply Theorem \ref{thm-2} to the self-normalized sum
\begin{equation*}
	\tilde{T}_k=\frac{\sum_{j=1}^{k}\tilde{Y}_j}{ ( \sum_{ j = 1 }^k \tilde{Y}_j^2 )^{ 1/2 } }.
\end{equation*}
Notice that $ \{ \tilde{Y}_j\}_{ 1 \leq j \leq k } $ are one-dependent random variables and the bounds in \eqref{dd1}, \eqref{dd2} and \eqref{rho-k}, i.e., the conditions in Theorem \ref{thm-2}, are satisfied. Recall that $ m = [ n^{ \alpha  } ] $ and $ k = [ n^{ 1 - \alpha }  ] $. We obtain by \eqref{re-dep} that there is a constant $ d_0 $ depending on $ \tau, \alpha, \omega_1, a_1 $ and $ a_2  $, such that
\begin{align*}
	\P ( \tilde{T}_k \geq x ) & =  \left[  1 - \Phi \left( \frac { x  } { \sqrt{ 1 + O ( m^{ -1  } ) } } \right) \right] \Big( 1 + O \Big( \frac { 1 + x^2  } { n^{ ( 1 - \alpha  ) / 4 } }   \Big)  \Big)   \\
	& = \Big[ 1 - \Phi (x) \Big] \Big( 1 + O \Big( \frac { 1 + x^2  }  { n^{ ( 1 - \alpha )/ 4  } } +  \frac { 1 +  x^2  } { n^{ \alpha } } \Big)  \Big),  \nn
\end{align*}
uniformly for $ x \in (0, d_0 \min \{ n^{ ( 1 - \alpha  )/ 8 }, n^{ \alpha / 2 }  \} ) $. This completes the proof of Proposition \ref{prop_approx}. 

\subsection{Proof of Proposition \ref{prop_error}}  \label{pf-prop9}
For simplicity, we denote $ V_k^2 = \sum_{j = 1}^k Y_j^2  $ and $ \tilde{V}_k^2 = \sum_{j  = 1}^k \tilde{Y}_j^2  $. Recall the constant $\tilde{C}_2 $ in the lower bound shown in \eqref{dd2}. Set
\begin{align*}
	\mathcal{B} = \left\{  \max_{ 1 \leq j  \leq k } \big| Y_i-\tilde{Y}_i  \big|\leq \f{1}{n^2}, ~\tilde{V}_k^2\geq \f{\tilde{C}_2^2}{2}\cdot n  \right\}.
\end{align*}
Within the set $ \mathcal{B} $, it holds that 
\begin{align*}
	| V_k  - \tilde{V}_k |  & \leq  \tilde{V}_k^{ -1 } | V_k^2 - \tilde{V}_k^2  | \leq  n^{ -2  } \tilde{V}_k^{ -1 }   \sum_{ j  = 1 }^k ( 2 | \tilde{Y}_j  | + n^{ -2  } ) \nn \\
	& \leq    2 n^{ - 2  } k^{ 1/2 }  + \sqrt 2  \tilde{C}_2^{ -1 }   n^{ - 7 / 2 - \alpha  } \leq  ( 2 \tilde{C}_2^{ -1  } + 2  ) n^{ - \f 32 - \f {\alpha} { 2 } },
\end{align*}
hence when $ n  $ is sufficiently large,
\begin{align*}
	\big| T_k - \tilde{T}_k   \big| & = \frac { \Big| \Big( \sum_{ j = 1 }^k Y_j  \Big) \tilde{V_k} -  \Big( \sum_{ j = 1 }^k \tilde{Y}_j \Big) V_k \Big|  } { V_k \tilde{V}_k } \nn \\
	& \leq  \frac {   \sum_{ j = 1 }^k \big| Y_j  - \tilde{Y}_j  \big|  }  { V_k  } + \frac { \Big( \sum_{ j = 1 }^k \big| \tilde{Y}_j \big| \Big)   } {  \tilde{V}_k }\cdot \frac{\big| V_k - \tilde{V}_k \big| }{V_k} \nn \\
	& \leq C \Big(  n^{ - 1/2 } k n^{  -2 } + k^{ 1/2 } n^{ - 1/2 } n^{ - \f 32 - \f { \alpha } { 2 } } \Big) \nn \\
	& \leq  2 C n^{ - \f 32 - \alpha } \leq  2 C n^ { - 1 }.
\end{align*}

\noi Therefore,
\begin{align}
	\P \big( T_k \geq x  \big)  &  \leq \P \Big(  T_k \geq x , ~ \big| T_k - \tilde{T}_k  \big| \leq 2 C  n^{ - 1  }  \Big) + \P ( \mathcal{B}^c ) \label{upper}  \\
	&  \leq \P \Big(  \tilde{T}_k \geq x - 2 C  n^{ -1  }  \Big)  + \P ( \mathcal{B}^c  ),\nn
\end{align}
and similarly,
\begin{align}
	\P \big( T_k \geq x  \big) & \geq  \P \big( T_k \geq x,  ~ \mathcal{B} \big)\geq \P \Big( \tilde{T}_k \geq x + 2 C  n^{ -1 } , ~ \mathcal{B}  \Big) \label{lower} \\
	& \geq \P \Big( \tilde{T}_k \geq x + 2 C n^{ -1 }   \Big) - \P ( \mathcal{B}^c  ).    \nn
\end{align}
As for the bound of  $ \P ( \mathcal{B}^c  ) $, notice that
\begin{align*}
	\P ( \mathcal{B}^c  )  &\leq  \sum_{ j = 1 }^k  \P \big(  \big| Y_j - \tilde{Y}_j  \big| >   n^{  - 2   }  \big)  +  \P \big( \tilde{V }_k^2 \leq \tilde{C}_2^2 n / 2   \big).
	%	&\leq  \sum_{ j = 1 }^k  \P \Big(  \big| Y_j - \tilde{Y}_j  \big| >   n^{  - 2   }  \Big)  \\
	%	&\qquad \;\;+  \P\bigg( \sum_{ \ell =1}^{k/2} Y_{2 \ell-1}^2 \leq \f{\tilde{C}_2^2}{4} n \bigg)+\P\bigg( \sum_{ \ell =1}^{k/2} Y_{2 \ell}^2 \leq \f{\tilde{C}_2^2}{4} n \bigg) .
\end{align*}
By Chebyshev's inequality and \eqref{difference} with $r=2$, we obtain
\begin{align*}
	\sum_{ j = 1 }^k  \P \big(  \big| Y_j - \tilde{Y}_j  \big| >   n^{  - 2   }  \big)  &  \leq  n^{ 1 - \alpha } n^{ 4 + 2 \alpha  }  a_1^2 e^{ - 2a_2 n^{ \alpha \tau  } } \leq C e^{ - a_2 n^{ \alpha \tau  }   }.
\end{align*}

\noi In addition,
\begin{align*}
	\P ( \tilde{ V }_k^2 \leq \tilde{C}_2^2 n / 2   ) & \leq  \P \Big( \sum_{ \ell = 1  }^{ k / 2  } \tilde{Y}_{ 2 \ell - 1  }^2  <  \tilde{C}_2^2 n / 4  \Big)  +   \P \Big( \sum_{ \ell = 1  }^{ k / 2  } \tilde{Y}_{ 2 \ell   }^2  <  \tilde{C}_2 ^2 n / 4  \Big) ,
\end{align*}
where, without loss of generality, we assumed that $ k /2  $ is an integer.
Denoting $ W_{\ell} = - \tilde{Y}_{2 \ell -1}^2+\e(\tilde{Y}_{2 \ell -1}^2)$, $1\leq \ell \leq k/2$,  \eqref{dd1} suggests that
\begin{gather*}
	W_{\ell} \leq  \e ( \tilde{Y}_{ 2 \ell - 1 }^2 ) \leq \tilde{C}_1^2 m ,\quad \quad \e ( W_{ \ell } )=0 ,     \\
	{\rm and}\quad \sum_{ \ell = 1 }^{k/2} \e ( W_{\ell}^2 )=\sum_{i=1}^{k/2}\var(\tilde{Y}_{2i-1}^2)  \leq \sum_{\ell=1}^{k/2}\e(\tilde{Y}_{ 2 \ell - 1 }^4) \leq  \tilde{C}_1^4  k m^2 / 2  .
\end{gather*}

\noi Observe that $ \{  Y_{ 2 \ell - 1  }, 1 \leq \ell \leq k / 2  \} $ is a sequence of independent random variables. It follows from Bernstein's inequality and \eqref{dd2} that
\begin{align*}
	\P\Big( \sum_{ \ell =1}^{k/2} &Y_{2 \ell-1}^2 \leq \f{\tilde{C}_2^2}{4} n \Big)\\
	&= \P\Big( \sum_{\ell=1}^{k/2} W_{\ell} \geq -\f{\tilde{C}_2^2}{4}  n + \sum_{i=1}^{k/2} \e \big(\tilde{Y}_{2i-1}^2 \big) \Big) \leq \P\Big( \sum_{\ell=1}^{k/2}W_{\ell} \geq \f{\tilde{C}_2^2}{4} n \Big) \\
	&\leq \exp \Big\{  -\f{ \big((\tilde{C}_2^2 n)/4\big)^2 } { 2 \big(  \tilde{C}_1^2 m \cdot \tilde{C}_2^2 n / 4  +  \tilde{C}_1^4 k m ^2 / 2   \big) }  \Big\}  \leq \exp \{ - C n^{1-\alpha} \},
\end{align*}
for some constant $C $ depends on $\tau, \alpha, \omega_1,   a_1$ and $a_2$. By the same token,
\begin{equation*}
	\P\Big( \sum_{ \ell =1}^{k/2}Y_{2 \ell}^2\leq \f{\tilde{C}_2^2}{4} n \Big) \leq  \exp \{ - C n^{1-\alpha} \}.
\end{equation*}
\vskip 0.5em
\noi Overall, we obtain
\begin{align}
	\P ( \mathcal{B}^c ) \leq C e^{  - a_2 n^{ \alpha \tau  }   } +   2e^{ - C  n^{1-\alpha}}.
\end{align}
The proof of Proposition \ref{prop_error} is completed.

\subsection{Proof of Theorem \ref{thm-winsor-2}} \label{pf-thm-4.2}

Define $ f (x) =  x \mathbbm{1} \{ | x | \leq \tau  \}  + \tau \mathbbm{1} \{ x > \tau \} - \tau \mathbbm{1} \{ x < - \tau \}$.
In view of the expression shown in \eqref{general-self} and the notations presented above Theorem \ref{thm-winsor-2}, the main idea is to apply Theorem \ref{thm-1} and  estimate the corresponding bias-corrected term $ \Psi_x^* $ and error terms $ L_{3, n} $ and $ R_x $ under fourth moment, which are defined by 
\begin{align*}
	L_{3, n} & =    \frac { \e | f (Y) - \tilde{\mu} |^3 } {  \sqrt n  \sigma_1^3 }  + \frac { \e | f (Y) -  {\mu} |^3 } {  \sqrt n  \sigma_2^3 }    \\
	\Psi_{x}^* & = \exp \Big\{  ( \frac {\sigma_2} {\sigma_1} )^3  x^3 n^{ - 1/2}  \Big ( \frac  43 \gamma^3     \frac { \e ( f (Y) - \tilde{\mu} )^3 } {     \sigma_1^3  } \\
	&  \hspace{4cm} -  2  \gamma^2    \frac { \e [ (  f (Y) - \tilde{\mu} ) ( f (Y) -  {\mu} )^2 ] } {  \sigma_1   \sigma_2^2  } \Big )\Big \} \\
	R_x & =  \frac { (1 + x)^4  } {  n} \Big(   \frac { \e | f (Y) - \tilde{\mu} |^4 } {    \sigma_1^4 }  + \frac { \e | f (Y) -  {\mu} |^4 } {   \sigma_2^4 }   \Big)  \\ 
	& \quad + n  \e \bigg[ \exp \Big\{ \frac {  ( f ( Y ) - \tilde{\mu} )^2  } { ( f ( Y ) - \mu )^2   +  \sigma_2^2  }  \cdot \frac {\sigma_2^2} {\sigma_1^2 } \Big \}   \mathbbm{1} \big( (1 + x)  \frac { | f (Y) - \tilde{\mu} | } {\sqrt{ n }  \sigma_1 }   > 1  \big) \bigg],
\end{align*} 
where $ \gamma = \f 12 (1 + \f c x )  $.

Throughout the rest of the proof, $ A  $ is an absolute positive constant which may vary at each appearance. When the higher moment $ \e [Y^4 ]  < \infty$ is assumed, we can figure out more accurate estimate for the bias-corrected term $ \Psi_x^* $. Let us collect the bounds for some crucial quantities in the proposition below. The proof of Proposition \ref{le-bounds-2} will be given in Section \ref{pf-le-bounds-2}. 
\begin{proposition} \label{le-bounds-2}
	If $ \e |Y|^4  < \infty $, then we have there exists an absolute positive constant $A$ such that  
	\begin{gather}
		|  \sigma_1^2 - \sigma^2 | \leq \frac { A \e | Y |^4   } { \tau^2 }, \quad  | \sigma_2^2 - \sigma^2 | \leq \frac { A \e | Y |^4   } { \tau ^2 } ,  \quad	| \mu - \tilde{\mu} | \leq \frac {  \e | Y |^4  } { \tau^3  }.  \label{eq_22}
	\end{gather}
	Moreover, when $\tau$ satisfies \eqref{cond-2}, we have 
	\begin{equation}
		\Big| \frac { \sigma_2^2 } {\sigma_1^2}  - 1 \Big|  \leq \frac { A (\e | Y |^4 )^2 } { \sigma^2 \tau^6  }, \label{eq_32}
	\end{equation}
	and there exist absolute positive constants $A$ and $ a_0 $ such that 
	\begin{equation}
		\Big|  \frac { 1 - \Phi \big( \frac { \sigma_2 } { \sigma_1 }  x + c \big) }  { 1 - \Phi (x) }  - 1 \Big| \leq  \frac { A ( 1 + x   ) \sqrt n \e | Y |^4  } { \sigma  \tau^3 }  \label{eq_42}
	\end{equation}
	for $ 0 <  x \leq a_0 \frac {\tau^3 \sigma} { \sqrt n \e | Y |^4 } $, and
	\begin{equation}
		R_x \leq \frac { A  (1 + x)^4 \e | Y |^4  } { n  \sigma^4  }.  \label{eq_52}
	\end{equation}
	In addition, there exist positive constants $a_1 $ and $ C $ depending on $ \sigma $, $ \e  |Y|^3  $ and $ \e |Y|^4 $ such that 
	\begin{equation}
		\Psi_x^*  = \exp \Big\{  -  \frac { x^3  \e ( Y - \mu )^3  } { 3 \sqrt n \sigma^3  }  \Big \} \Big[ 1 + O_1 \Big( \frac {  x^3  } { \sqrt n \tau  }  + \frac { x^2   } {\tau^3   } \Big) \Big] \label{eq_62}
	\end{equation}
	for $ 2 < x \leq  a_1 \min \{ \frac { \tau^3  } { \sqrt n },  (\sqrt n \tau )^{1/3}, \tau^{3/2} \}  $, where $ O_1 $ is a bounded quantity such that $ | O_1 | \leq C $.
\end{proposition} 

We next consider two cases of $ 0 < x \leq 2  $ and $ x > 2  $, separately. 
First for $ 0 < x \leq 2  $, it is immediate that $ 1 - \Phi (x)  \geq 1 - \Phi(2) $. Moreover, it follows from Proposition \ref{prop-small-x}  and \eqref{eq_42} that 
\begin{align*}
	\P ( S_{\tau, n}^* > x  )  = [   1 - \Phi (x) ]   \Big[ 1 + O \Big( \frac { 1 } { \sqrt n } +  \frac {\sqrt n } {\tau^3 } \Big) \Big]. 
\end{align*}
Since the bias term $ \exp \{ - \frac { x^3 \e (Y - \mu)^3 } { 3 \sqrt n \sigma^3 } \} = 1 + O ( \frac {1 } {\sqrt n} ) $ for $ 0 < x \leq  2 $, we obtain 
\begin{align}
	\P ( S_{\tau, n}^* > x  )  = [   1 - \Phi (x) ]  \exp \Big\{ - \frac { x^3 \e (Y - \mu)^3 } { 3 \sqrt n \sigma^3 } \Big\}   \Big[ 1 + O \Big( \frac { 1 } { \sqrt n } +  \frac {\sqrt n } {\tau^3 } \Big) \Big]  \label{eq-14}
\end{align}
for $ 0 < x \leq 2  $, which completes the proof of \eqref{re-2} for $ 0 < x \leq 2 $. Then we deal with the case of $ x > 2  $. Applying Theorem \ref{thm-1} to \eqref{general-self} yields
\begin{align}
	\P ( S_{\tau, n}^* > x  )  = \Big[ 1 - \Phi ( \frac {\sigma_2} {\sigma_1} x + c   ) \Big] \Psi_x^* e^{ O_1 R_x } \Big( 1 + O_2 (1 + \frac {\sigma_2} {\sigma_1} x ) L_{3, n}  \Big) \label{xlarge-2}
\end{align}
uniformly for $ | c |  \leq x / 5$ and $ x  > 2 $ satisfying \eqref{G1} and \eqref{G2}. Note that \eqref{G1} and \eqref{G2} are satisfied when $ 2 < x \leq  a_2 \frac {\sqrt n \sigma^4} {\e | Y |^4}  $ for some constant $a_2 > 0$. By plugging in the results in Proposition \ref{le-bounds-2}, we can obtain 
\begin{align*}
	\P ( S_{\tau, n}^* > x   ) & =  [   1 - \Phi (x) ]  \exp \Big\{ - \frac { x^3 \e (Y - \mu)^3 } { 3 \sqrt n \sigma^3 } \Big\} \\
	& \quad \times \Big[  1 + O_1 \Big(  \frac { x \sqrt n } { \tau^3  } + \frac {  x ^3  } { \sqrt n \tau  } + \frac { x ^4 } { n } +  \frac { x  } { \sqrt n } + \frac {  x ^2  } {\tau^3 } \Big)    \Big]
\end{align*}
for $ 2 < x \leq a_3  \min \{ \tau^3 n^{-1/2} , ( \sqrt n \tau )^{1/3} , n^{1/4}, \tau^{3/2}\} $. Observe that when $ \tau  $ satisfies \eqref{cond-2}, $ \tau^{3/2} \geq O (n^{1/4}) $ and $ \frac { 1 + x^2 } { \tau^3 }  \leq A \frac { (1 + x ) \sqrt n } { \tau^3  }$ for $ x = O (\sqrt n) $. Moreover, by the basic inequality that $   a + b >  (a^2 b)^{1/3} $ for $ a > 0$ and $ b> 0 $, we have $ \frac { x^3 } { \sqrt n \tau  } \leq \frac { x^4 } { n } + \frac { x \sqrt n  } { \tau^3 } $ and $ ( \sqrt n  \tau )^{1/3} \geq \min \{ n^{1/4},  \tau^{3} n^{ - 1/2} \} $. 
Hence the desired result \eqref{re-2} holds for $ x \in (2, c_2  \min \{ n^{1/4},   \tau^3 n^{- 1/2 }   \}  ) $, which together with the result \eqref{eq-14} for $  0 < x \leq 2 $ completes the proof of Theorem \ref{thm-winsor-2}.

\subsection{Proof of Theorem \ref{thm-trimmed-2}} \label{pf-thm-trimmed-1} 
Define $ g(Y) = (Y - \mu) \IF\{ | Y |  \leq \tau \}  $. Recall that
\begin{equation}
	\P (U_{\tau, n}^* > x ) = \P \Big( \frac { S_n^{\circ}  - \delta } { V_n^{\circ} } > \frac { \sigma_4 } { \sigma_3 } x \Big), \label{self-trim}
\end{equation}
where
\begin{align*}
	S_n^{\circ} &= \sum_{i = 1}^n \frac {  g(Y_i) - \mu_0 } { \sqrt n \sigma_3 }, \quad  (V_n^{\circ} )^2 = \sum_{i = 1}^n \frac { [ g(Y_i)]^2 } { n \sigma_4^2 } \quad
	\mbox{and} \quad \delta = \frac { \sqrt n \mu_0 } { \sigma_3 }.
\end{align*}
The main idea is to apply Theorem \ref{thm-1} and  estimate the corresponding bias-corrected term $ \Psi_x^* $ and error terms $ L_{3, n} $ and $ R_x $ under fourth moment, which are defined by 
\begin{align*}
	L_{3, n} & =    \frac { \e | g (Y) -  {\mu}_0 |^3 } {  \sqrt n  \sigma_3^3 }  + \frac { \e | g (Y)   |^3 } {  \sqrt n  \sigma_4^3 }    \\
	\Psi_{x}^* & = \exp \Big\{  ( \frac {\sigma_4} {\sigma_3} )^3  x^3 n^{ - 1/2}  \Big ( \frac  43 \gamma^3     \frac { \e ( g(Y) - {\mu}_0 )^3 } {     \sigma_3^3  } \\
	&  \hspace{4cm} -  2  \gamma^2    \frac { \e [ (  g (Y) - {\mu}_0 ) ( g(Y) )^2 ] } {  \sigma_3   \sigma_4^2  } \Big )\Big \} \\
	R_x & =  \frac { (1 + x)^4  } {  n} \Big(   \frac { \e | g (Y) - {\mu}_0 |^4 } {    \sigma_3^4 }  + \frac { \e | g (Y)   |^4 } {   \sigma_4^4 }   \Big)  \\ 
	& \quad + n  \e \bigg[ \exp \Big\{ \frac {  ( g ( Y ) - {\mu}_0 )^2  } { ( g(Y) )^2   +  \sigma_4^2  }  \cdot \frac {\sigma_4^2} {\sigma_3^2 } \Big \}   \mathbbm{1} \big( (1 + x)  \frac { | g (Y) - {\mu}_0 | } {\sqrt{ n }  \sigma_3 }   > 1  \big) \bigg],
\end{align*} 
where $ \gamma = \f 12 (1 + \f \delta x )  $.

When the higher moment $ \e [Y^4 ]  < \infty$ is assumed, we can figure out more accurate estimate for the bias-corrected term $ \Psi_x^* $. Let us collect the bounds for some crucial quantities in the proposition below. The proof of Proposition \ref{le-bounds-t1} will be given in Section \ref{pf-le-bounds-3}. 
\begin{proposition} \label{le-bounds-t1}
	If $ \e |Y|^4  < \infty $, then we have there exists an absolute positive constant $A$ such that  
	\begin{gather}
		|  \sigma_3^2 - \sigma^2 | \leq \frac { A \e | Y |^4   } { \tau^2 }, \quad  | \sigma_4^2 - \sigma^2 | \leq \frac { A \e | Y |^4   } { \tau ^2 } ,  \quad	|  \mu_0 | \leq \frac {  \e | Y |^4  } { \tau^3  }.  \label{eq_22_t}
	\end{gather}
	Moreover, when $\tau$ satisfies \eqref{cond-2}, we have 
	\begin{equation}
		\Big| \frac { \sigma_4^2 } {\sigma_3^2}  - 1 \Big|  \leq \frac { A (\e | Y |^4 )^2 } { \sigma^2 \tau^6  }, \label{eq_32_t}
	\end{equation}
	and there exist absolute positive constants $A$ and $ a_0 $ such that 
	\begin{equation}
		\Big|  \frac { 1 - \Phi \big( \frac { \sigma_4 } { \sigma_3 }  x + \delta \big) }  { 1 - \Phi (x) }  - 1 \Big| \leq  \frac { A ( 1 + x   ) \sqrt n \e | Y |^4  } { \sigma  \tau^3 }  \label{eq_42_t}
	\end{equation}
	for $ 0 <  x \leq a_0 \frac {\tau^3 \sigma} { \sqrt n \e | Y |^4 } $, and
	\begin{equation}
		R_x \leq \frac { A  (1 + x)^4 \e | Y |^4  } { n  \sigma^4  }.  \label{eq_52_t}
	\end{equation}
	In addition, there exist positive constants $a_1 $ and $ C $ depending on $ \sigma $, $ \e  |Y|^3  $ and $ \e |Y|^4 $ such that 
	\begin{equation}
		\Psi_x^*  = \exp \Big\{  -  \frac { x^3  \e ( Y - \mu )^3  } { 3 \sqrt n \sigma^3  }  \Big \} \Big[ 1 + O_1 \Big( \frac {  x^3  } { \sqrt n \tau  }  + \frac { x^2   } {\tau^3   } \Big) \Big] \label{eq_62_t}
	\end{equation}
	for $ 2 < x \leq  a_1 \min \{ \frac { \tau^3  } { \sqrt n },  (\sqrt n \tau )^{1/3}, \tau^{3/2} \}  $, where $ O_1 $ is a bounded quantity such that $ | O_1 | \leq C $.
\end{proposition} 

We next consider two cases of $ 0 < x \leq 2  $ and $ x > 2  $, separately. 
First for $ 0 < x \leq 2  $, it is immediate that $ 1 - \Phi (x)  \geq 1 - \Phi(2) $. Moreover, it follows from Proposition \ref{prop-small-x}  and \eqref{eq_42_t} that 
\begin{align*}
	\P ( U_{\tau, n}^* > x  )  = [   1 - \Phi (x) ]   \Big[ 1 + O \Big( \frac { 1 } { \sqrt n } +  \frac {\sqrt n } {\tau^3 } \Big) \Big]. 
\end{align*}
Since the bias term $ \exp \{ - \frac { x^3 \e (Y - \mu)^3 } { 3 \sqrt n \sigma^3 } \} = 1 + O ( \frac {1 } {\sqrt n} ) $ for $ 0 < x \leq  2 $, we obtain 
\begin{align}
	\P ( U_{\tau, n}^* > x  )  = [   1 - \Phi (x) ]  \exp \Big\{ - \frac { x^3 \e (Y - \mu)^3 } { 3 \sqrt n \sigma^3 } \Big\}   \Big[ 1 + O \Big( \frac { 1 } { \sqrt n } +  \frac {\sqrt n } {\tau^3 } \Big) \Big]   \label{eq-14-t}
\end{align}
for $ 0 < x \leq 2  $. Next we deal with the case of $ x > 2  $. Applying Theorem \ref{thm-1} to \eqref{self-trim} yields
\begin{align}
	\P ( U_{\tau, n}^* > x  )  = \Big[ 1 - \Phi ( \frac {\sigma_4} {\sigma_3} x + \delta   ) \Big] \Psi_x^* e^{ O_1 R_x } \Big( 1 + O_2 (1 + \frac {\sigma_4} {\sigma_3} x ) L_{3, n}  \Big) \label{xlarge-2-t}
\end{align}
uniformly for $ | \delta |  \leq x / 5$ and $ x  > 2 $ satisfying \eqref{G1} and \eqref{G2}. Note that \eqref{G1} and \eqref{G2} are satisfied when $ 2 < x \leq  a_2 \frac {\sqrt n \sigma^4} {\e | Y |^4}  $ for some constant $a_2 > 0$. By plugging in the results in Proposition \ref{le-bounds-t1}, we can obtain 
\begin{align*}
	\P ( U_{\tau, n}^* > x   ) & =  [   1 - \Phi (x) ]  \exp \Big\{ - \frac { x^3 \e (Y - \mu)^3 } { 3 \sqrt n \sigma^3 } \Big\} \\
	& \quad \times \Big[  1 + O_1 \Big(  \frac { x \sqrt n } { \tau^3  } + \frac {  x ^3  } { \sqrt n \tau  } + \frac { x ^4 } { n } +  \frac { x  } { \sqrt n } + \frac {  x ^2  } {\tau^3 } \Big)    \Big]
\end{align*}
for $ 2 < x \leq a_3  \min \{ \tau^3 n^{-1/2} , ( \sqrt n \tau )^{1/3} , n^{1/4}, \tau^{3/2}\} $. Observe that when $ \tau  $ satisfies \eqref{cond-2}, $ \tau^{3/2} \geq O (n^{1/4}) $ and $ \frac { 1 + x^2 } { \tau^3 }  \leq A \frac { (1 + x ) \sqrt n } { \tau^3  }$ for $ x = O (\sqrt n) $. Moreover, by the basic inequality that $   a + b >  (a^2 b)^{1/3} $ for $ a > 0$ and $ b> 0 $, we have $ \frac { x^3 } { \sqrt n \tau  } \leq \frac { x^4 } { n } + \frac { x \sqrt n  } { \tau^3 } $ and $ ( \sqrt n  \tau )^{1/3} \geq \min \{ n^{1/4},  \tau^{3} n^{ - 1/2} \} $. 
Hence the desired result holds for $ x \in (2, c_2  \min \{ n^{1/4},   \tau^3 n^{- 1/2 }   \}  ) $, which together with the result \eqref{eq-14-t} for $  0 < x \leq 2 $ completes the proof of Theorem \ref{thm-trimmed-2}.

\subsection{Proof of Theorem \ref{thm-winsor}} \label{pf-thm-4.1}	 
Recall that $ f (x) =  x \mathbbm{1} \{ | x | \leq \tau  \}  + \tau \mathbbm{1} \{ x > \tau \} - \tau \mathbbm{1} \{ x < - \tau \}$. In view of the expression shown in \eqref{general-self} and the notations presented above Theorem \ref{thm-winsor-2}, the main idea is to apply Theorem \ref{thm-1} and  estimate the corresponding error terms $ L_{3, n},  \Psi_x^* $ and $ R_x $, which are defined by 
\begin{align*}
	L_{3, n} & =    \frac { \e | f (Y) - \tilde{\mu} |^3 } {  \sqrt n  \sigma_1^3 }  + \frac { \e | f (Y) -  {\mu} |^3 } {  \sqrt n  \sigma_2^3 }    \\
	\Psi_{x}^* & = \exp \Big\{  ( \frac {\sigma_2} {\sigma_1} )^3  x^3 n^{ - 1/2}  \Big ( \frac  43 \gamma^3     \frac { \e ( f (Y) - \tilde{\mu} )^3 } {     \sigma_1^3  } \\
	&  \hspace{4cm} -  2  \gamma^2    \frac { \e [ (  f (Y) - \tilde{\mu} ) ( f (Y) -  {\mu} )^2 ] } {  \sigma_1   \sigma_2^2  } \Big )\Big \} \\
	R_x & =  \frac { (1 + x)^3  } { \sqrt n} \Big(   \frac { \e | f (Y) - \tilde{\mu} |^3 } {    \sigma_1^3 }  + \frac { \e | f (Y) -  {\mu} |^3 } {   \sigma_2^3 }   \Big)  \\ 
	& \quad + n  \e \bigg[ \exp \Big\{ \frac {  ( f ( Y ) - \tilde{\mu} )^2  } { ( f ( Y ) - \mu )^2   +  \sigma_2^2  }  \cdot \frac {\sigma_2^2} {\sigma_1^2 } \Big \}   \mathbbm{1} \big( (1 + x)  \frac { | f (Y) - \tilde{\mu} | } {\sqrt{ n }  \sigma_1 }   > 1  \big) \bigg],
\end{align*} 
where $ \gamma = \f 12 (1 + \f c x )  $.

We collect the bounds for some crucial quantities in the proposition below. The proof of Proposition \ref{le-bounds} will be omited since it is similar to the proof of Proposition \ref{le-bounds-2}. 
\begin{proposition} \label{le-bounds}
	If $ \e [ | Y |^3 ] < \infty $, then we have there exists an absolute positive constant $A$ such that  
	\begin{gather}
		\e | f (Y) - \mu |^3 \leq A \e  | Y |^3, \quad  \e | f (Y) - \tilde{\mu} |^3 \leq A \e  | Y |^3 \label{eq_1} \\
		|  \sigma_1^2 - \sigma^2 | \leq \frac { A \e | Y |^3   } { \tau }, \quad  | \sigma_2^2 - \sigma^2 | \leq \frac { A \e | Y |^3   } { \tau } ,  \quad	| \mu - \tilde{\mu} | \leq \frac {  \e | Y |^3 } { \tau^2  } .  \label{eq_2}
	\end{gather}
	Moreover, when $\tau$ satisfies \eqref{cond}, we have 
	\begin{equation}
		\Big| \frac { \sigma_2^2 } {\sigma_1^2}  - 1 \Big|  \leq \frac { A (\e | Y |^3 )^2 } { \sigma^2 \tau^4  }, \label{eq_3}
	\end{equation}
	and there exist absolute constants $A$ and $ a_0 $ such that 
	\begin{equation}
		\Big|  \frac { 1 - \Phi \big( \frac { \sigma_2 } { \sigma_1 }  x + c \big) }  { 1 - \Phi (x) }  - 1 \Big| \leq  \frac { A ( 1 + x   ) \sqrt n \e | Y |^3  } { \sigma  \tau^2 }  \label{eq_4}
	\end{equation}
	for $ 0 <  x \leq a_0 \frac {\tau^2 \sigma} { \sqrt n \e | Y |^3 } $,  and 
	\begin{equation}
		R_x \leq \frac { A  (1 + x)^3 \e | Y |^3  } { \sigma^3 \sqrt n }. \label{eq_5}
	\end{equation}
\end{proposition}

We next consider two cases of $ 0 < x \leq 2 $ and $ x > 2  $, separately.  
First for $ 0 < x \leq 2 $, it follows from Proposition \ref{prop-small-x} that 
\begin{align}
	\Big|  \P ( S_{\tau, n}^* > x  )  - \Big[ 1 - \Phi \big( \frac { \sigma_2 } {\sigma_1} x + c \big) \Big] \Big| \leq A L_{3, n}.
\end{align}
From Proposition \ref{le-bounds} that when $ \tau  $ satisfies \eqref{cond} with some large $a_1$, we can see that 
\begin{align}
	& L_{3, n} \leq \frac { A \e  |Y|^3 } { \sigma^3 \sqrt n }, \label{L_3n} \\
	\mbox{and} \quad & \Big|  \frac {  1 - \Phi \big( \frac { \sigma_2 } {\sigma_1} x + c \big)  } { 1 - \Phi (x) } - 1 \Big| \leq \frac { A \sqrt n \e |Y|^3 } { \sigma \tau^2 }.
\end{align}
Therefore, it holds  for $ 0 < x \leq 2 $ that 
\begin{equation}
	\frac {\P ( S_{\tau, n} ^*  > x  ) } { 1 - \Phi(x) } = 1 + O_1  \Big( \f {  \e | Y|^3  } { \sigma^3 \sqrt{n} } + \frac { \sqrt{n} \e | Y|^3  } { \sigma \tau^2 }   \Big),
\end{equation}
with $ | O_1 | \leq A $. The desired result \eqref{re1} has been proved for $ 0 < x \leq 2 $.

Now we proceed to prove for $ x > 2 $. Applying Theorem \ref{thm-1} to \eqref{general-self} yields
\begin{align}
	\P ( S_{\tau, n}^* > x  )  = \Big[ 1 - \Phi ( \frac {\sigma_2} {\sigma_1} x + c   ) \Big] \Psi_x^* e^{ O_1 R_x } \Big( 1 + O_2 (1 + \frac {\sigma_2} {\sigma_1} x ) L_{3, n}  \Big) \label{xlarge}
\end{align}
uniformly for $ | c |  \leq x / 5$ and $ x  > 2 $ satisfying \eqref{G1} and \eqref{G2}. When $ \tau $ satisfies \eqref{cond}, we can obtain from Proposition \ref{le-bounds} that $ | \sigma_1 / \sigma_2  - 1 | \leq A $ and 
\begin{equation}
	\Big|  \frac {  1 - \Phi \big( \frac { \sigma_2 } {\sigma_1} x + c \big)  } { 1 - \Phi (x) } - 1 \Big| \leq \frac { A ( 1 + x   ) \sqrt n \e | Y |^3  } { \sigma  \tau^2 }  \label{phi}
\end{equation} 
for $  x \leq a_2 \frac {\tau^2 \sigma} { \sqrt n \e | Y |^3 } $, and 
\begin{equation}
	R_x \leq  \frac { A  (1 + x)^3 \e | Y |^3  } { \sigma^3 \sqrt n }.  \label{Rx}
\end{equation}  
Moreover, since $ | \frac c x | \leq  1/ 5  $ as $ \tau  $ satisfies \eqref{cond} and $ x > 2 $, we have
\begin{equation}
	\Psi_x^* \leq \exp \Big\{   \frac { A x^3 \e | Y |^3 } {\sigma^3 \sqrt n } \Big\}.  \label{Psix}
\end{equation}
In addition, note that the conditions of $x$ shown in \eqref{G1} and $ \eqref{G2}$ are satisfied when $ 0 <  x \leq a_2 \min \{  \frac { \sqrt n \sigma^3 } {\e |Y|^3 } , \frac {\tau^2 \sigma} { \sqrt n \e | Y |^3 } \} $. Consequently, the desired result \eqref{re1} for $ x > 2 $ are derived by substituting \eqref{L_3n} and \eqref{phi}--\eqref{Psix} into \eqref{xlarge}. This completes the proof of Theorem \ref{thm-winsor}.

\subsection{Proof of Theorem \ref{thm-trimmed}} \label{pf-thm-trimmed-2}
Define $ g(Y) = (Y - \mu) \IF\{ | Y |  \leq \tau \}  $. Recall that
\begin{equation}
	\P (U_{\tau, n}^* > x ) = \P \Big( \frac { S_n^{\circ}  - \delta } { V_n^{\circ} } > \frac { \sigma_4 } { \sigma_3 } x \Big), \label{self-trim-1}
\end{equation}
where
\begin{align*}
	S_n^{\circ} &= \sum_{i = 1}^n \frac {  g(Y_i) - \mu_0 } { \sqrt n \sigma_3 }, \quad  (V_n^{\circ} )^2 = \sum_{i = 1}^n \frac { [ g(Y_i)]^2 } { n \sigma_4^2 } \quad
	\mbox{and} \quad \delta = \frac { \sqrt n \mu_0 } { \sigma_3 }.
\end{align*}
The main idea is to apply Theorem \ref{thm-1} and  estimate the corresponding bias-corrected term $ \Psi_x^* $ and error terms $ L_{3, n} $ and $ R_x $ under fourth moment, which are defined by 
\begin{align*}
	L_{3, n} & =    \frac { \e | g (Y) -  {\mu}_0 |^3 } {  \sqrt n  \sigma_3^3 }  + \frac { \e | g (Y)   |^3 } {  \sqrt n  \sigma_4^3 }    \\
	\Psi_{x}^* & = \exp \Big\{  ( \frac {\sigma_4} {\sigma_3} )^3  x^3 n^{ - 1/2}  \Big ( \frac  43 \gamma^3     \frac { \e ( g(Y) - {\mu}_0 )^3 } {     \sigma_3^3  } \\
	&  \hspace{4cm} -  2  \gamma^2    \frac { \e [ (  g (Y) - {\mu}_0 ) ( g(Y) )^2 ] } {  \sigma_3   \sigma_4^2  } \Big )\Big \} \\
	R_x & =  \frac { (1 + x)^4  } {  n} \Big(   \frac { \e | g (Y) - {\mu}_0 |^4 } {    \sigma_3^4 }  + \frac { \e | g (Y)   |^4 } {   \sigma_4^4 }   \Big)  \\ 
	& \quad + n  \e \bigg[ \exp \Big\{ \frac {  ( g ( Y ) - {\mu}_0 )^2  } { ( g(Y) )^2   +  \sigma_4^2  }  \cdot \frac {\sigma_4^2} {\sigma_3^2 } \Big \}   \mathbbm{1} \big( (1 + x)  \frac { | g (Y) - {\mu}_0 | } {\sqrt{ n }  \sigma_3 }   > 1  \big) \bigg],
\end{align*} 
where $ \gamma = \f 12 (1 + \f \delta x )  $.

We collect the bounds for some crucial quantities in the proposition below. The proof of Proposition \ref{le-bounds-t2} will be omited since it is similar to the proof of Proposition \ref{le-bounds-t1}. 
\begin{proposition} \label{le-bounds-t2}
	If $ \e [ | Y |^3 ] < \infty $, then we have there exists an absolute positive constant $A$ such that  
	\begin{gather}
		\e | g (Y) - \mu_0 |^3 \leq A \e  | Y |^3, \quad  \e | g (Y)  |^3 \leq A \e  | Y |^3, \label{eq_1-t} \\
		|  \sigma_3^2 - \sigma^2 | \leq \frac { A \e | Y |^3   } { \tau }, \quad  | \sigma_4^2 - \sigma^2 | \leq \frac { A \e | Y |^3   } { \tau } ,  \quad	| {\mu}_0  | \leq \frac {  \e | Y |^3 } { \tau^2  }.  \label{eq_2-t}
	\end{gather}
	Moreover, when $\tau$ satisfies \eqref{cond}, we have 
	\begin{equation}
		\Big| \frac { \sigma_4^2 } {\sigma_3^2}  - 1 \Big|  \leq \frac { A (\e | Y |^3 )^2 } { \sigma^2 \tau^4  }, \label{eq_3-t}
	\end{equation}
	and there exist absolute constants $A$ and $ a_0 $ such that 
	\begin{equation}
		\Big|  \frac { 1 - \Phi \big( \frac { \sigma_4 } { \sigma_3 }  x + \delta \big) }  { 1 - \Phi (x) }  - 1 \Big| \leq  \frac { A ( 1 + x   ) \sqrt n \e | Y |^3  } { \sigma  \tau^2 }  \label{eq_4-t}
	\end{equation}
	for $ 0 <  x \leq a_0 \frac {\tau^2 \sigma} { \sqrt n \e | Y |^3 } $,  and 
	\begin{equation}
		R_x \leq \frac { A  (1 + x)^3 \e | Y |^3  } { \sigma^3 \sqrt n }. \label{eq_5-t}
	\end{equation}
\end{proposition}

We next consider two cases of $ 0 < x \leq 2 $ and $ x > 2  $, separately.  
First for $ 0 < x \leq 2 $, it follows from Proposition \ref{prop-small-x} that 
\begin{align}
	\Big|  \P ( U_{\tau, n}^* > x  )  - \Big[ 1 - \Phi \big( \frac { \sigma_4 } {\sigma_3} x + \delta \big) \Big] \Big| \leq A L_{3, n}
\end{align}
From Proposition \ref{le-bounds-t2} that when $ \tau  $ satisfies \eqref{cond} with some large $a_1$, we can see that 
\begin{align}
	& L_{3, n} \leq \frac { A \e  |Y|^3 } { \sigma^3 \sqrt n }, \label{L_3n_t} \\
	\mbox{and} \quad & \Big|  \frac {  1 - \Phi \big( \frac { \sigma_4 } {\sigma_3} x + \delta \big)  } { 1 - \Phi (x) } - 1 \Big| \leq \frac { A \sqrt n \e |Y|^3 } { \sigma \tau^2 }.
\end{align}
Therefore, it holds  for $ 0 < x \leq 2 $ that 
\begin{equation}
	\frac {\P ( U_{\tau, n} ^*  > x  ) } { 1 - \Phi(x) } = 1 + O_1  \Big( \f {  \e | Y|^3  } { \sigma^3 \sqrt{n} } + \frac { \sqrt{n} \e | Y|^3  } { \sigma \tau^2 }   \Big),
\end{equation}
with $ | O_1 | \leq A $. The desired result has been proved for $ 0 < x \leq 2 $.

Now we proceed to prove for $ x > 2 $. Applying Theorem \ref{thm-1} to \eqref{self-trim-1} yields 
\begin{align}
	\P ( U_{\tau, n}^* > x  )  = \Big[ 1 - \Phi ( \frac {\sigma_4} {\sigma_3} x + \delta   ) \Big] \Psi_x^* e^{ O_1 R_x } \Big( 1 + O_2 (1 + \frac {\sigma_4} {\sigma_3} x ) L_{3, n}  \Big) \label{xlarge-t}
\end{align}
uniformly for $ | c |  \leq x / 5$ and $ x  > 2 $ satisfying \eqref{G1} and \eqref{G2}. When $ \tau $ satisfies \eqref{cond}, we can obtain from Proposition \ref{le-bounds-t2} that $ | \sigma_3 / \sigma_4  - 1 | \leq A $ and 
\begin{equation}
	\Big|  \frac {  1 - \Phi \big( \frac { \sigma_4 } {\sigma_3} x + \delta \big)  } { 1 - \Phi (x) } - 1 \Big| \leq \frac { A ( 1 + x   ) \sqrt n \e | Y |^3  } { \sigma  \tau^2 }  \label{phi-t}
\end{equation} 
for $  x \leq a_2 \frac {\tau^2 \sigma} { \sqrt n \e | Y |^3 } $, and 
\begin{equation}
	R_x \leq  \frac { A  (1 + x)^3 \e | Y |^3  } { \sigma^3 \sqrt n }.  \label{Rx-t}
\end{equation}  
Moreover, since $ | \frac \delta x | \leq  1/ 5  $ as $ \tau  $ satisfies \eqref{cond} and $ x > 2 $, we have
\begin{equation}
	\Psi_x^* \leq \exp \Big\{   \frac { A x^3 \e | Y |^3 } {\sigma^3 \sqrt n } \Big\}.  \label{Psix-t}
\end{equation}
In addition, note that the conditions of $x$ shown in \eqref{G1} and $ \eqref{G2}$ are satisfied when $ 0 <  x \leq a_2 \min \{  \frac { \sqrt n \sigma^3 } {\e |Y|^3 } , \frac {\tau^2 \sigma} { \sqrt n \e | Y |^3 } \} $. Consequently, the desired result \eqref{re1} for $ x > 2 $ are derived by substituting \eqref{L_3n_t} and \eqref{phi-t}--\eqref{Psix-t} into \eqref{xlarge-t}. This completes the proof of Theorem \ref{thm-trimmed}.

\subsection{Proof of Theorem \ref{thm4.3}} \label{pf-thm4.3}

Notice that $\Phi^{-1}(1-\alpha/2p)= 2 (1+o(1)) \sqrt{\log 2 p /\alpha } = o(n^{1/6} ) $ lays in the range of Theorem \ref{thm-winsor}, hence by Theorem \ref{thm-winsor}, 
\begin{equation*}
	\sum_{ j = 1 }^p \Big[ \P (\mu_j<L_j)+\P(\mu_j >U_j) \Big] \leq 2\sum_{ j = 1 }^p \P \big(S_{\tau, n} ^*  > \Phi^{-1} (1-\alpha/2p) \big)=\alpha+o(1),
\end{equation*}
which completes the proof of Theorem \ref{thm4.3}.

\subsection{Proof of Proposition \ref{le-bounds-2}} \label{pf-le-bounds-2}
From the definition of $ f (\cdot) $, it is obvious that $ | f (Y ) | \leq | Y |$ and $ | f (Y)  - Y| \leq |Y|  \mathbbm{1} ( | Y | >  \tau  ) $. First for the bound for $ | \sigma_1^2 - \sigma^2 | $, it holds that 
\begin{align*}
	| \sigma_1^2 - \sigma^2  | & = | \e ( f (Y) - \tilde{\mu} )^2 - \e ( Y - \mu )^2 | \\
	& \leq \e \big\{ ( 2 | Y | + 2 \e | Y | ) ( | Y | \mathbbm{1} ( | Y | > \tau )  +  \e [ | Y | \mathbbm{1} ( | Y | > \tau ) ] ) \big\} \\
	& \leq  \frac { A \e | Y |^4  } { \tau^2 }.
\end{align*}
Similarly, we can obtain the same bound for $ | \sigma_2^2 - \sigma^2 | $. Since $ | f (Y)  - Y| \leq |Y|  \mathbbm{1} ( | Y | >  \tau  ) $, we have
\begin{align*}
	| \mu - \tilde{\mu} |  & = | \e f (Y) - \e Y | \\
	& \leq \e | f (Y) - Y | \leq \e [  |Y|  \mathbbm{1} ( | Y | >  \tau  )  ] \leq \frac { \e | Y |^4 } { \tau^3 }.
\end{align*}
As for \eqref{eq_32}, it follows by \eqref{eq_22} that when $ \tau  $ satisfies \eqref{cond-2},
\begin{align}
	\Big| \frac {\sigma_2^2} {\sigma_1^2 } - 1  \Big| & = \frac { | \e[ (f (Y) - \mu)^2  ] - \e [ (f (Y) - \tilde{\mu})^2  ] | } {  \sigma_1^2  }  \nn \\
	& =\frac { ( \mu - \tilde{\mu} )^2 } { \sigma_1^2  }  \leq \frac { ( \mu - \tilde{\mu} )^2   } { \sigma^2 - \frac {A \e | Y |^4  } {\tau^2 } } \nn  \\
	& \leq \frac { A  ( \mu - \tilde{\mu} )^2  } { \sigma^2 } \leq  \frac { A ( \e | Y |^4 )^2 } { \sigma^2 \tau^6  }.  \label{eq_6}
\end{align}
Now we proceed to show \eqref{eq_42}. First for $ 0 < x \leq 2 $, it holds that $ 1 - \Phi (2) \leq 1 - \Phi (x)  \leq 1 $. Then we can deduce that 
\begin{align}
	\Big|  \frac { 1 - \Phi \big( \frac { \sigma_2 } { \sigma_1 }  x + c \big) }  { 1 - \Phi (x) }  - 1 \Big|   
	&  \leq A \Big| \Phi \big( \frac { \sigma_2 } { \sigma_1 }  x + c \big) - \Phi (x)  \Big| \nn \\
	& \leq A \Big[  \big( \frac {\sigma_2} {\sigma_1} - 1 \big) x + c  \Big]. \label{eq_7}
\end{align}
Recalling that $ | c  | = \frac { \sqrt n | \mu - \tilde{\mu} | } {\sigma_1} \leq \frac {A \sqrt n \e | Y |^4 } { \sigma \tau^3 } $, and $ | \frac {\sigma_1} {\sigma_2} - 1 | \leq |  \frac {\sigma_1^2} {\sigma_2^2} - 1  | \leq   \frac { A ( \e | Y |^4 )^2 } { \sigma^2 \tau^6  }$. Therefore, it follows that for $ 0 < x \leq 2 $, 
\begin{align}
	\Big|  \frac { 1 - \Phi \big( \frac { \sigma_2 } { \sigma_1 }  x + c \big) }  { 1 - \Phi (x) }  - 1 \Big|   & \leq A  \Big(  \frac {   (\e | Y |^4 )^2  } { \sigma^2  \tau^6 }    + \frac { \sqrt n \e | Y |^4   } { \sigma \tau^3  }  \Big) \nn \\
	&  \leq \frac {  A \sqrt n \e | Y |^4   } { \sigma \tau^3  } , 
\end{align}
where the last inequality is derived by using the fact that $  \frac { \sqrt n \e | Y |^4   } { \sigma \tau^3  }  \geq  \frac {  A (\e | Y |^4 )^2  } { \sigma^2  \tau^6 }   $ for $\tau $ satisfying \eqref{cond-2}. For $ x > 2  $, it holds that when $ x \leq a_0  \frac { \tau^3 \sigma } { \e | Y |^4  } $, 
\begin{align*}
	\Big|  \frac { 1 - \Phi \big( \frac { \sigma_2 } { \sigma_1 }  x + c \big) }  { 1 - \Phi (x) }  - 1 \Big|   & \leq A  \Big[  \big( \frac { \sigma_2 } {\sigma_1} - 1 \big) x^2 + c x  \Big] \exp\Big\{ A  \Big[  \big( \frac { \sigma_2^2 } {\sigma_1^2} - 1 \big)  x^2 + c x    \Big] \Big\} \nn  \\
	& \leq A  \Big[  \big( \frac { \sigma_2 } {\sigma_1} - 1 \big) x^2 + c x  \Big]   \nn \\
	& \leq \frac { A x \sqrt n \e | Y |^4 } {   \tau^3  \sigma } .
\end{align*}
The above results for the two cases of $ 0 < x \leq 2 $ and $ 2 < x \leq a_0   \frac { \tau^3 \sigma } { \e | Y |^4  }$ yield the desired result \eqref{eq_42}.

Next we proceed to prove \eqref{eq_52}.  Recall that 
\begin{align}
	R_x \leq A \frac { ( 1 + x  )^4 \e | Y |^4 } {  n \sigma^4 } + r_{x},  \label{Rx_bound}
\end{align}
where  
\begin{align*}
	r_x = n   \e \bigg[ \exp \Big\{ \frac {  ( f ( Y ) - \tilde{\mu} )^2  } { ( f ( Y ) - \mu )^2   +  \sigma_2^2  }  \cdot \frac {\sigma_2^2} {\sigma_1^2 } \Big \}   \mathbbm{1} \big( (1 + x)  \frac { | f (Y) - \tilde{\mu} | } {\sqrt{ n }  \sigma_1 }   > 1  \big) \bigg].
\end{align*}
The goal is to upper  bound the exponential part. By \eqref{eq_32} and the range \eqref{cond-2} for $\tau$, we can obtain
\begin{align}
	\frac {  ( f ( Y ) - \tilde{\mu} )^2  } { ( f ( Y ) - \mu )^2   +  \sigma_2^2  }  \cdot \frac {\sigma_2^2} {\sigma_1^2 }  & \leq A \frac { (  f ( Y ) - \tilde{\mu} )^2 } { ( f ( Y ) - \mu )^2   +  \sigma_2^2  } \nn  \\
	& = A \frac { (  f ( Y ) - \tilde{\mu} )^2 } { (  f ( Y ) - \tilde{\mu} )^2 + \kappa_1 ( f (Y) - \tilde{\mu} ) + \kappa_2  } \nn \\
	& = \frac { A } { 1 + \kappa_1 ( f (Y) - \tilde{\mu} )^{-1} + \kappa_2 (  f (Y) - \tilde{\mu}  )^{ - 2} }, \label{eq-13}
\end{align}
where 
\begin{align*}
	\kappa_1 & =  2 ( \tilde{\mu}  - \mu  ) \leq \frac { A \e | Y |^4 } { \tau^3  }, \\
	\kappa_2 & = \sigma_2^2 + ( \tilde{\mu} - \mu )^2 \geq A \sigma^2  - \frac { A ( \e | Y |^4 )^2  } { \tau^6 } \geq A \sigma^2 
\end{align*}
for $ \tau $ satisfying \eqref{cond-2}. Observing that for any $ x \in \mathbbm{R} $, 
\begin{equation*}
	1 + \kappa_1 x + \kappa_2 x^2 \geq 1 - \frac { \kappa_1^2 } { 4 \kappa_2 }, 
\end{equation*}
hence when \eqref{cond-2} hold with large number $ c_1 $ such that $ \kappa_1^2 \leq \kappa_2 $, we obtain
\begin{align}
	\frac {  ( f ( Y ) - \tilde{\mu} )^2  } { ( f ( Y ) - \mu )^2   +  \sigma_2^2  }  \cdot \frac {\sigma_2^2} {\sigma_1^2 } \leq  A. 
\end{align}
Consequently, it follows that 
\begin{align*}
	r_x \leq A \frac { (1 + x)^4 \e ( | f (Y) - \tilde{\mu} |^4 ) } {  n \sigma_1^4   } \leq \frac { A ( 1 + x )^4 \e | Y |^4  } { n \sigma^4    },
\end{align*}  
which combining with \eqref{Rx_bound} yields the desired result \eqref{eq_52}.

Finally, we deal with the proof of \eqref{eq_62} when $ \tau  $ satisfies \eqref{cond-2}. It follows from \eqref{eq_32} that 
\begin{align}
	\Big| \Big( \frac { \sigma_2 } {\sigma_1} \Big)^3 - 1  \Big| \leq \frac { A ( \e | Y |^4  )^2  } { \sigma^2 \tau^6 }. \label{eq-9}
\end{align}
Since $ \gamma = \frac 12 ( 1 + \frac c x ) $ and $ | c | \leq  \frac {A \sqrt n | \mu - \tilde{\mu} | } { \sigma } \leq \frac {  A \sqrt n \e | Y |^4  } { \sigma \tau^3   }  \leq A_0 $ for some absolute constant $A_0$ when $ \tau $ satisfies \eqref{cond-2}, we have $ | c |  / x \leq A_0 / 2 $  for $ x > 2 $ and $\tau $ satisfying \eqref{cond-2}. In addition, 
\begin{align}
	\Big|   \gamma  - \frac 1 2 \Big|   \leq \frac {  A \sqrt n \e | Y |^4 } { x \sigma \tau^3  }, \label{eq-10}
\end{align}
and we can obtain from \eqref{eq_32} that 
\begin{align}
	&  \Big| \frac { \e ( f(Y)  - \tilde{\mu})^3 } { \sigma_1^3 }  - \frac { \e ( Y - \mu )^3 } { \sigma^3 } \Big| \nn \\
	&  \leq   \frac { | \e (  f(Y) - \tilde{\mu} )^3  - \e (Y- \mu)^3   | } { \sigma_1^3  }  + | \e ( Y - \mu )^3 | \Big| \frac { 1 } { \sigma_1^3} - \frac {1} {\sigma^3} \Big| \nn \\
	& \leq A \Big( \frac { \e | Y |^4  } {  \sigma^3 \tau  }  + \frac { | \e ( Y - \mu )^3 | ( \e | Y |^4 )^2 } { \sigma^5 \tau^6 }  \Big) = O \Big( \frac { 1 } { \tau  } \Big). \label{eq-11}
\end{align}
Similarly, 
\begin{align}
	\Big| \frac { \e [ ( f(Y)  - \tilde{\mu}) (f (Y) - \mu)^2 ] } { \sigma_1 \sigma_2^2 }  - \frac { \e ( Y - \mu )^3 } { \sigma^3 } \Big|  
	& = O \Big( \frac { 1 } { \tau  } \Big).  \label{eq-12}
\end{align}
Then the desired result \eqref{eq_62} is derived by plugging in the above bounds in \eqref{eq-9}--\eqref{eq-12}. This completes the proof of Proposition \ref{le-bounds-2}.

\subsection{Proof of Proposition \ref{le-bounds-t1}} \label{pf-le-bounds-3}

Recall the definition that $ g (Y) = (Y -\mu) \IF\{ | Y | \leq \tau \} $.  First for the bound for $ | \sigma_1^2 - \sigma^2 | $, it holds that 
\begin{align*}
	| \sigma_1^2 - \sigma^2  | & = | \e \{ [ g (Y) - {\mu}_0 - ( Y - \mu ) ] [ g (Y) - {\mu}_0 + ( Y - \mu )]  \} | \\
	& \leq A  \tau^{-2}  \e \big\{  \big( | Y |^3 +   | \mu | Y^2   +  \e | Y |^3 + |\mu| \e Y^2  \big) ( |Y| + |\mu| + \e|Y| \big) \big\} \\
	& \leq  \frac { A \e | Y |^4  } { \tau^2 }. 
\end{align*}
Similarly, we can obtain the same bound for $ | \sigma_2^2 - \sigma^2 | $. Moreover, we have
\begin{equation*}
	| \mu_0 |   = | \e g (Y)  |  \leq  \e [ | Y - \mu | \IF \{ |Y|  > \tau \} ]   \leq \frac { \e | Y |^4 } { \tau^3 }.
\end{equation*}
As for \eqref{eq_32_t}, it follows by \eqref{eq_22_t} that when $ \tau  $ satisfies \eqref{cond-2},
\begin{align}
	\Big| \frac {\sigma_4^2} {\sigma_3^2 } - 1  \Big| & = \frac { | \e[ (g (Y) - \mu_0)^2  ] - \e [ (g (Y) )^2  ] | } {  \sigma_3^2  }  \nn \\
	& =\frac {  \mu_0^2 } { \sigma_3^2  }  \leq \frac {  \mu_0^2   } { \sigma^2 - \frac {A \e | Y |^4  } {\tau^2 } } \nn  \\
	& \leq \frac { A   \mu_0^2 } { \sigma^2 } \leq  \frac { A ( \e | Y |^4 )^2 } { \sigma^2 \tau^6  }.  \label{eq_6_t}
\end{align}
Therefore, \eqref{eq_42_t} can be proved by the same procedure of proving \eqref{eq_42}.  

Next we proceed to prove \eqref{eq_52_t}.  Recall that 
\begin{align}
	R_x \leq A \frac { ( 1 + x  )^4 \e | Y |^4 } {  n \sigma^4 } + r_{x},  \label{Rx_bound-t}
\end{align}
where  
\begin{align*}
	r_x = n   \e \bigg[ \exp \Big\{ \frac {  ( g ( Y ) - {\mu}_0 )^2  } { ( g ( Y )  )^2   +  \sigma_4^2  }  \cdot \frac {\sigma_4^2} {\sigma_3^2 } \Big \}   \mathbbm{1} \big( (1 + x)  \frac { | g (Y) - {\mu}_0 | } {\sqrt{ n }  \sigma_3 }   > 1  \big) \bigg].
\end{align*}
The key point is to upper  bound the exponential part. By \eqref{eq_32_t} and the range \eqref{cond-2} for $\tau$, we can obtain
\begin{align}
	\frac {  ( g ( Y ) - {\mu}_0  )^2  } { ( g ( Y )   )^2   +  \sigma_4^2  }  \cdot \frac {\sigma_4^2} {\sigma_3^2 }  & \leq A \frac { (  g ( Y ) - {\mu}_0 )^2 } { ( g ( Y )  )^2   +  \sigma_4^2  } \nn  \\
	& = A \frac { (  g ( Y ) - {\mu}_0 )^2 } { (  g ( Y ) -\mu_0  )^2 + \kappa_1 ( g (Y) - {\mu}_0 ) + \kappa_2  } \nn \\
	& = \frac { A } { 1 + \kappa_1 ( g (Y) - {\mu}_0 )^{-1} + \kappa_2 (  g (Y) - {\mu}_0 )^{ - 2} }, \label{eq-13-t}
\end{align}
where 
\begin{align*}
	\kappa_1 & =  2  {\mu}_0 \quad \mbox{satisfying} ~| \kappa_1 | \leq \frac { A \e | Y |^4 } { \tau^3  }, \\
	\kappa_2 & = \sigma_4^2 + ( {\mu}_0 )^2 \geq A \sigma^2  - \frac { A ( \e | Y |^4 )^2  } { \tau^6 } \geq A \sigma^2 
\end{align*}
for $ \tau $ satisfying \eqref{cond-2}. Observing that for any $ x \in \mathbbm{R} $, 
\begin{equation*}
	1 + \kappa_1 x + \kappa_2 x^2 \geq 1 - \frac { \kappa_1^2 } { 4 \kappa_2 }, 
\end{equation*}
hence when \eqref{cond-2} hold with large number $ c_1 $ such that $ \kappa_1^2 \leq \kappa_2 $, we obtain
\begin{align}
	\frac {  ( g ( Y ) - {\mu}_0  )^2  } { ( g ( Y )   )^2   +  \sigma_4^2  }  \cdot \frac {\sigma_4^2} {\sigma_3^2 }  \leq  A. 
\end{align}
Consequently, it follows that 
\begin{align*}
	r_x \leq A \frac { (1 + x)^4 \e ( | g (Y) - {\mu}_0 |^4 ) } {  n \sigma_3^4   } \leq \frac { A ( 1 + x )^4 \e | Y |^4  } { n \sigma^4    },
\end{align*}  
which combining with \eqref{Rx_bound} yields the desired result \eqref{eq_52_t}.

Finally, we deal with the proof of \eqref{eq_62_t} when $ \tau  $ satisfies \eqref{cond-2}. It follows from \eqref{eq_32_t} that 
\begin{align}
	\Big| \Big( \frac { \sigma_4 } {\sigma_3} \Big)^3 - 1  \Big| \leq \frac { A ( \e | Y |^4  )^2  } { \sigma^2 \tau^6 }. \label{eq-9-t}
\end{align}
Since $ \gamma = \frac 12 ( 1 + \frac \delta x ) $ and $ | \delta | \leq  \frac {A \sqrt n | \mu_0 | } { \sigma } \leq \frac {  A \sqrt n \e | Y |^4  } { \sigma \tau^3   }  \leq A_0 $ for some absolute constant $A_0$ when $ \tau $ satisfies \eqref{cond-2}, we have $ | \delta |  / x \leq A_0 / 2 $  for $ x > 2 $ and $\tau $ satisfying \eqref{cond-2}. In addition, 
\begin{align}
	\Big|   \gamma  - \frac 1 2 \Big|   \leq \frac {  A \sqrt n \e | Y |^4 } { x \sigma \tau^3  }, \label{eq-10-t}
\end{align}
and we can obtain from \eqref{eq_32_t} that 
\begin{align}
	&  \Big| \frac { \e ( g(Y)  - {\mu}_0 )^3 } { \sigma_3^3 }  - \frac { \e ( Y - \mu )^3 } { \sigma^3 } \Big| \nn \\
	&  \leq  A  \frac { | \e (  g(Y) - {\mu}_0 )^3  - \e (Y- \mu)^3   | } { \sigma_3^3  }  + A [  \e | Y - \mu |^3 ] \Big| \frac { 1 } { \sigma_3^3} - \frac {1} {\sigma^3} \Big| \nn \\
	& \leq A \Big( \frac { \e | Y |^4  } {  \sigma^3 \tau  }  + \frac { | \e ( Y - \mu )^3 | ( \e | Y |^4 )^2 } { \sigma^5 \tau^6 }  \Big) = O \Big( \frac { 1 } { \tau  } \Big). \label{eq-11-t}
\end{align}
Similarly, 
\begin{align}
	\Big| \frac { \e [ ( g(Y)  - {\mu}_0 ) (g (Y))^2 ] } { \sigma_3 \sigma_4^2 }  - \frac { \e ( Y - \mu )^3 } { \sigma^3 } \Big|  
	& = O \Big( \frac { 1 } { \tau  } \Big).  \label{eq-12-t}
\end{align}
Then the desired result \eqref{eq_62_t} is derived by plugging in the above bounds in \eqref{eq-9-t}--\eqref{eq-12-t}. This completes the proof of Proposition \ref{le-bounds-t1}.

\subsection{Proof of Corollary \ref{cor}} \label{pf-cor}
We apply Theorem \ref{thm-main} with $ Y_i^2 = \e X_i^2  $, $ c_0 = 0 $ and $ c = 0 $. In this case,
\begin{align*}
	\Psi_x^* = \exp \Big \{   \frac { x^3   \e X_1^3 } { 6 \sigma^3 \sqrt{n}  }   \Big\} \quad , L_{3, n} \leq  \frac { 2   \e | X_1 |^3 } { \sigma^3 \sqrt n    }  , \quad  \delta_{x, 1} \leq \frac { 4 (1 + x)^4 \e X_1^4  }  { n^2 \sigma^4  }.
\end{align*}
For $ x \leq t_0 \sigma \sqrt n / 4 $, we have $ 2 x / (\sigma \sqrt n )  \leq t_0 / 2$ and thus
\begin{align*}
	r_{x, 1}  & = \e \bigg[ \exp \Big\{ \min \Big( \frac { X_1^2  } { \sigma^2   } , \frac { 2 x X_1 } { \sigma \sqrt n  } \Big) \Big\}  \IF ( | (1 + x) X_1 | > \sigma \sqrt n  ) \bigg]  \\
	& \leq  \e  \bigg[ \exp \Big\{  \frac { 2 x X_1 } { \sigma \sqrt n  }   \Big\}  \IF ( | (1 + x) X_1 | > \sigma \sqrt n  ) \bigg] \\
	& \leq  \e  \big[  e^{ t_0 X_1 / 2  } \IF \{  (1 + x ) | X_1| > \sigma \sqrt n \} \big]  \\
	& \leq (1 + x)^4 n^{-2 } \e \big[ | X_1 |^4 e^{ t_0 X_1 / 2  }   \big] \\
	& \leq C t_0^{-4}  \e \big[ e^{ 3 t_0 X_1 / 4  }  \big]  ( 1 + x  )^4 n^{ - 2  } .
\end{align*}
Therefore it follows from Theorem \ref{thm-main} that
\begin{align*}
	\frac { \P ( S_n > x  \sigma \sqrt n  ) } { 1 - \Phi (x) }  = \exp \Big\{ \frac { x^3   \e X_1^3 } { 6 \sigma^3 \sqrt{n}  }  +  O \Big( \frac { ( 1 + x )^4  } { n  } \Big) \Big \} \Big[ 1 + O \Big( \frac { 1 + x } { \sqrt n } \Big) \Big]
\end{align*}
uniformly for $ 0 < x \leq O (n^{1/2}) $. Hence \eqref{classical-sum} holds for $ 0 < x \leq O (n^{1/4}) $. This completes the proof of Corollary \ref{cor}.

\renewcommand{\theequation}{B.\arabic{equation}}
\setcounter{equation}{0}

\section{Additional lemmas and proofs}  \label{SecB}

\subsection{Proof of Lemma \ref{lemma-1}} \label{pf-lemma-1}
Recall the notation $W_i=2x{X_i} - x^{2}{Y_i}^2$. Define
$$
\nu_i = 	\min\left\{  \frac{{X_i}^2}{{Y_i}^2+c_0\e {Y_i}^2} +c_0x^2\e {Y_i}^2 , 2xX_i \right\}  .
$$
Observe that $ W_i \leq \nu_i $.
For $ 1/4 \leq  \lambda \leq 3 /4 $ and $ k \in \{ 0, 1, 2, 3 \}$, by the elementary inequalities $ |s|^{k} e^{- s} \leq  c_1(k) $  and $ s^{k} e^{\lambda s }  \leq c_2(k) e^{s}$ for $ s \geq 0 $, we obtain
\begin{align}
	&\e \{ |W_i|^k   e^{\la W_i} \mathbbm{1}(|(1+x)X_i| \geq 1 ) \} \label{Wtail}\\
	& \hspace{1cm}   \leq   A \e \{e^{(0\lor W_i)} \mathbbm{1}(|(1+x)X_i| \geq 1 ) \} \nn \\
	& \hspace{1cm}  \leq   A \big( \e \{  e^{\nu_i} \mathbbm{1}(|(1+x)X_i| \geq 1 ) \}  + \delta_{x,i} \big) \nn \\
	& \hspace{1cm} \leq A  R_{x,i} ,  \nn
\end{align}
where the last inequality comes from the fact that $ c_0 x^2 \e {Y_i}^2 \le  1/4  $ for $x$ satisfying \eqref{G2}. Additionally, there is a constant $A_1(k, s_0) $ such that $ |s|^k e^{s} \leq A_1 (k, s_0) $ if $ s \leq s_0 $. Hence,
\begin{align*}
	&\e \{ |W_i|^k   e^{\la W_i} \mathbbm{1}(|(1+x)X_i| \leq 1, |(1+x)Y_i| \geq 1 ) \}   \\
	& \leq   A \P(|(1+x)X_i| \leq 1, |(1+x)Y_i| \geq 1 ) \\
	& \leq  A \de_{x,i}.
\end{align*}
For simplicity, denote $\mathcal{U}_i = \{ |(1+x) X_i| \leq 1, |(1+x) Y_i| \leq 1   \}$. We have
\begin{align}
	\e \{ {W_i}^k e^{\la W_i}  \} =  \e \{ {W_i}^k e^{\la W_i} \mathbbm{1}(\mathcal{U}_i ) \} + O(1) R_{x,i}, \label{e1}
\end{align}
where and hereafter $ O(1) $ is a bounded quantity.

By Taylor expansion,
\begin{align}
	\e \{ e^{\la W_i} \mathbbm{1}(\mathcal{U}_i ) \} &=  1+\la \e \{ W_i \mathbbm{1}(\mathcal{U}_i  )\}+\frac{{\la}^2}{2} \e \{  {W_i}^2 \mathbbm{1}(\mathcal{U}_i )\}  \label{0th} \\
	& \quad   +\frac{{\la}^3}{6} \e \{  {W_i}^3 \mathbbm{1}(\mathcal{U}_i  )\}  + O(1) {\la}^4 \e \{  {W_i}^4   \mathbbm{1}(\mathcal{U}_i  )\}. \nn
\end{align}
The analysis for the terms on the right hand side of the above equality is similar. We calculate the third term for example. Notice that
%	\begin{align*}
%	(I). \quad \e \{W_i \mathbbm{1}(|xX_i| \leq 1,|xY_i| \leq 1)\} & = - x^2\e {Y_i}^2 -2x\e \{ X_i \mathbbm{1}(|xX_i| \geq 1) \} \\
%	& \quad \; -2x\e \{ X_i \mathbbm{1}(|xX_i| \leq 1,|xY_i| \geq 1) \} \\
%	& \quad \; +x^2\e \{ {Y_i}^2 \mathbbm{1}(|xY_i| \geq 1) \} \\
%	& \quad \; +x^2\e \{ {Y_i}^2 \mathbbm{1}(|xX_i| \geq 1,|xY_i| \leq 1) ,
%	\end{align*}
%	And since
%	\begin{align*}
%	|x^2\e \{ {Y_i}^2 \mathbbm{1}(|xX_i| \geq 1,|xY_i| \leq 1)|+|2x\e \{ X_i \mathbbm{1}(|xX_i| \geq 1) \}| \leq 3\de_{x,i},\\
%	|2x\e \{ X_i \mathbbm{1}(|xX_i| \leq 1,|xY_i| \geq 1) \}|+|x^2\e \{ {Y_i}^2 \mathbbm{1}(|xY_i| \geq 1) \}| \leq 3\de_{x,i},
%	\end{align*}
%	hence
%	\begin{align}
%	\e \{W_i \mathbbm{1}(|xX_i| \leq 1,|xY_i| \leq 1)\} = -x^2\e {Y_i}^2 +O(1)\de_{x,i}.
%	\end{align}
\begin{align*}
	\e \{{W_i}^2 \mathbbm{1}(\mathcal{U}_i )\} &= 4x^2\e {X_i}^2 - 4x^3\e X_i{Y_i}^2 +x^4\e \{ {Y_i}^4 \mathbbm{1}(\mathcal{U}_i ) \}\\
	& \quad  -4x^2\e \{ {X_i}^2 \mathbbm{1}(|(1+x)X_i| \geq 1) \} \\
	& \quad -4x^2\e \{{X_i}^2 \mathbbm{1}(|(1+x)X_i| \leq 1,|(1+x)Y_i| \geq 1) \} \\
	& \quad +4x^3\e \{ X_i{Y_i}^2 \mathbbm{1}(|(1+x)X_i| \leq 1,|(1+x)Y_i| \geq 1) \} \\
	& \quad +4x^3\e \{ X_i{Y_i}^2 \mathbbm{1}(|(1+x)X_i| \geq 1,|(1+x)Y_i| \leq 1) \} \\
	& \quad +4x^3\e \{ X_i{Y_i}^2 \mathbbm{1}(|(1+x)X_i| \geq 1,|(1+x)Y_i| \geq 1) \}.
\end{align*}
By the basic inequality $ ab^2 \leq a^3 + b^3  $ for $ a > 0 , b> 0 $ and Chebyshev's inequality, we obtain
\begin{align*}
	\e \{{W_i}^2 \mathbbm{1}(\mathcal{U}_i )\} = 4x^2\e {X_i}^2 - 4x^3\e X_i{Y_i}^2 +O(1)\de_{x,i}.
\end{align*}
In the same manner,
\begin{gather*}
	\e \{W_i \mathbbm{1}(\mathcal{U}_i)\} = -x^2\e {Y_i}^2 +O(1)\de_{x,i}, \quad \e \{{W_i}^3 \mathbbm{1}(\mathcal{U}_i )\} = 8x^3\e {X_i}^3 +O(1)\de_{x,i}, \\
	\quad {\rm and} \hspace{1cm} |\e \{{W_i}^4   \mathbbm{1}(\mathcal{U}_i )\}| \leq O(1)\de_{x,i}.
\end{gather*}	
%	\begin{align*}
%	\quad \e \{{W_i}^3 I(|xX_i| \leq 1,|xY_i| \leq 1)\} & = 8x^3\e {X_i}^3 - 8x^3\e \{{X_i}^3 I(|xX_i| \geq 1) \} \\
%	& \quad \; -8x^3\e \{{X_i}^3 I(|xX_i| \leq 1), |xY_i| \geq 1) \} \\
%	& \quad \; -12x^4\e \{ {X_i}^2{Y_i}^2 I(|xX_i| \leq 1,|xY_i| \leq 1) \} \\
%	& \quad \; +6x^5\e \{ X_i{Y_i}^4 I(|xX_i| \leq 1,|xY_i| \leq 1) \} \\
%	& \quad \; -x^6\e \{ {Y_i}^6 I(|xX_i| \leq 1,|xY_i| \leq 1) \}
%	\end{align*}
%As for
%\begin{align*}
%& \quad \; |8x^3\e \{{X_i}^3 I(|xX_i| \geq 1) \}| \leq 8\de_{1,i} \\
%& \quad \; |8x^3\e \{{X_i}^3 I(|xX_i| \leq 1), |xY_i| \geq 1) \}| + |6x^5\e \{ X_i{Y_i}^4 I(|xX_i| \leq 1,|xY_i| \leq 1) \} | \leq 8\de_{2,i} \\
%& \quad \; |x^6\e \{ {Y_i}^6 I(|xX_i| \leq 1,|xY_i| \leq 1) \}| \leq \de_{2,i} \\
%& \quad \; |12x^4\e \{ {X_i}^2{Y_i}^2 I(|xX_i| \leq 1,|xY_i| \leq 1) \} | \\
%& \leq 12[x^4\e \{ {|X_i|}^4I(|xX_i| \leq 1) \}]^{1/2} [x^4\e \{ {|Y_i|}^4I(|xY_i| \leq 1) \}]^{1/2} \\
%& \leq 12(\de_{1,i}+\de_{2,i})
%\end{align*}
%Hence
%\begin{align}
%\e \{{W_i}^3 I(|xX_i| \leq 1,|xY_i| \leq 1)\} = 8x^3\e {X_i}^3 +O(1)\de_i
%\end{align}
%where $|O(1)| \leq 21$.

%	\begin{align*}
%	(IV). \quad |\e \{{W_i}^4 e^{\la \eta_i} I(|xX_i| \leq 1,|xY_i| \leq 1)\}| & \leq \e \{{W_i}^4 e^{\la |W_i|}I(|xX_i| \leq 1,|xY_i| \leq 1)\}\\
%	& \leq 3\e \{16x^4{X_i}^4 I(|xX_i| \leq 1), |xY_i| \leq 1) \} \\
%	& \quad \; + 3\e \{ {x^8{Y_i}^8}  I(|xX_i| \leq 1,|xY_i| \leq 1) \} \\
%	& \leq 48\de_i
%	\end{align*}
Therefore \eqref{Wtail}--\eqref{0th} yield
\begin{align}
	\e e^{\la W_i} &= 1 - \la x^2\e {Y_i}^2 +2 {\la}^2x^2 \e {X_i}^2  \label{e2} \\
	&\qquad -2{\la}^2x^3 \e {X_i}{Y_i}^2 +\frac{4}{3}{\la}^3x^3 \e {X_i}^3 +O(1)R_{x,i}. \nn
\end{align}
Because $x$ satisfies \eqref{G2},
\begin{align}
	|\e e^{\la W_i} - 1| & \leq 1/2 .  \nn
\end{align}
Furthermore, by Lemma \ref{lemma-0} we have
\begin{align*}
	|\e e^{\la W_i} - 1|^2 & \leq A \big[ x^4\e \{{Y_i}^4 \mathbbm{1}(|(1+x)Y_i|\leq 1)\} + x^4 \e \{{X_i}^4 \mathbbm{1}(|(1+x)X_i|\leq 1)\} \nn\\
	& \qquad + x^6 \e \{({X_i}^6+{Y_i}^6)\mathbbm{1}(\mathcal{U}_i)\} + R_{x,i} \big] \nn\\
	& \leq A R_{x, i},
\end{align*}
Because $ | \log (1+a) - a | \leq  a^2$ whenever $|a|\leq 1/2$, it follows from \eqref{e2} that
\begin{align*}
	& \log \e e^{\la W_i} = - \la x^2\e {Y_i}^2 +2 {\la}^2x^2 \e {X_i}^2 -2{\la}^2x^3 \e {X_i}{Y_i}^2 +\frac{4}{3}{\la}^3x^3 \e {X_i}^3 +O(1)R_{x,i},
\end{align*}
which completes the proof of \eqref{0-mo}.
\ignore{And by \eqref{e1}, for $k=1,2$ we have
	\begin{align}
		\e W_i^k e^{\la W_i} & = \e \{W_i^k e^{\la W_i}\mathbbm{1}(\mathcal{U})\} +O(1)R_{x,i} \nn \\
		& = \e \Big\{ W_i^k \Big(\sum_{j=0}^{3-k} \frac{\la^j}{j!} W_i^j \Big) \mathbbm{1}(\mathcal{U}) \Big\} \nn \\
		&\quad \; +\e \Big\{W_i^k \Big(e^{\la W_i}-\sum_{j=0}^{3-k} \frac{\la^j}{j!} W_i^j \Big) \mathbbm{1}(\mathcal{U})  \Big\} +O(1)R_{x,i} \nn \\
		& := P_{k1} +P_{k2}+O(1)R_{x,i}, \nn
		%	\e {W_i}^2 e^{\la W_i} & = \e \{{W_i}^2 e^{\la W_i}I(|xX_i|\leq 1,|xY_i|\leq 1)\} +O(1)R_{x,i} \nn \\
		%	& = \e \{{W_i}^2 (1+\la W_i )I(|xX_i|\leq 1,|xY_i|\leq 1)  \} \nn \\
		%	&\quad \; +\e \{{W_i}^2 (e^{\la W_i}-1-\la W_i )I(|xX_i|\leq 1,|xY_i|\leq 1)  \} +O(1)R_{x,i} \nn \\
		%	& = P_{21} +P_{22}, \nn \\
	\end{align}
	and
	\begin{align}
		\e |W_i|^3 e^{\la W_i} & = \e \{|W_i|^3 \mathbbm{1}(\mathcal{U})\}  +\e \{|W_i|^3 (e^{\la W_i}-1) \mathbbm{1}(\mathcal{U})\}+O(1)R_{x,i} \nn \\
		& := P_{31} +P_{32}+O(1)R_{x,i}, \nn
	\end{align}
	where we have
	\begin{align*}
		|P_{12}|+|P_{22}|+|P_{32}| & \leq O(1)\e \{|W_i|^4 e^{\la |W_i|}\mathbbm{1}(\mathcal{U})\}  \leq O(1) \e \{ |W_i|^4 \mathbbm{1}(\mathcal{U})\} \} \\
		& \leq O(1) [\e \{ x^4 {X_i}^4 \mathbbm{1}(\mathcal{U}) \} + \e \{ x^8 {Y_i}^8 \mathbbm{1}(\mathcal{U}) \}] \leq O(1)\de_{x,i}.
	\end{align*}
}
By the same token,
\begin{align*}
	\e W_ie^{\la W_i} & = \e \{W_i(1+\la W_i +\frac{{\la}^2}{2}{W_i}^2) \mathbbm{1}(\mathcal{U}_i)  \} +O(1)R_{x,i} \\
	& = -x^2 \e {Y_i}^2 +4\la x^2 \e {X_i}^2 -4\la x^3 \e X_i{Y_i}^2 +4{\la}^2 x^3 \e {X_i}^3 +O(1)R_{x,i}, \\
	\e {W_i}^2 e^{\la W_i} & = \e \{{W_i}^2 (1+\la W_i ) \mathbbm{1}(\mathcal{U}_i)  \}+O(1)R_{x,i} \\
	& = 4x^2\e {X_i}^2 -4x^3\e X_i{Y_i}^2 +8\la x^3 \e {X_i}^3 +O(1)R_{x,i}, \\
	\e |W_i|^3 e^{\la W_i} & = \e \{|W_i|^3  \mathbbm{1}(\mathcal{U}_i )\}+O(1)R_{x,i} \\
	& = O(1) x^3 (\e |X_i|^3 +\e |Y_i|^3)+ O(1)R_{x,i}.
\end{align*}
By a similar procedure to the proof of \eqref{0-mo}, we arrive at \eqref{1-mo}--\eqref{3-mo}. The proof of Lemma \ref{lemma-1} is completed.

\subsection{Proof of Lemma \ref{lemma-2}} \label{pf-lemma-2}
Recall the definition $ m \(\la\) = \sum_{i=1}^n \log{ \e e^{\la W_i} } $. Note that $m(\la_{\de})$ is well-defined under condition \eqref{maincon2} and $m''(\la)=\sum_{i=1}^{n} Var \wt W_i	>0$ for $\la>0$, $x \neq 0$ and nondegenerate $X_i$ and $Y_i$. It follows from Lemmas \ref{lemma-1} and \ref{lemma-0} that for $1/4 \leq \la \leq 3/4$ and $x$ satisfying \eqref{G1}--\eqref{G2},
\begin{align}
	m'(\la) & = \sum_{i=1}^{n} \e W_i e^{\la W_i} /\e e^{\la W_i} \label{term3} \\
	& = (4\la -1)x^2 +4{\la}^2 x^3 \sum_{i=1}^{n} \e {X_i}^3 -4\la x^3 \sum_{i=1}^{n} \e X_i{Y_i}^2 +O(1)R_{x} ,
	\nn
\end{align}
where $ | O(1) | \leq A $ for some absolute constant $A$.
Therefore, under \eqref{G1} with sufficiently small constant $ c_1 $, we have for $ | \de (x) | \leq x^2 /2  $ that
\begin{align*}
	m'(1/4) <  x^2+\de (x) <  m'(3/4),
\end{align*}
which combined with the fact $ m'' (\la) > 0  $ implies that the equation $m'(\la)=x^2+\de (x)$ has a unique solution $\la_{\de}$ such that $1/4 < \la_{\de} < 3/4$. Furthermore, by virtue of \eqref{term3}, it holds that
\begin{align*}
	\la_{\de} = \frac{1}{2} + \frac{\de (x)}{4x^2} - {\la_{\de}}^2 x \sum_{i=1}^{n} \e {X_i}^3 + \la_{\de} x \sum_{ i =1 }^n \e X_i{Y_i}^2 + O (1) x^{ - 2 } R_{x},
\end{align*}
and hence \eqref{pro-solu} follows. Again by \eqref{term3} and $1/4  < \la_{\de} <  3/4  $, $ 1/4  < \la_{\de_0} <  3/4   $, we obtain
\begin{align*}
	| \la_{\de} - \la_{\de_0} | \leq A \big(  x^{-2 } |  \de (x ) - \de_0 (x)  |   +  x L_{3, n}  +  x^{- 2} R_x      \big) .
\end{align*}
Therefore
\begin{align*}
	\Big|  \la_{\de} - \la_{\de_0} - \frac { \de (x)  - \de_0 (x)  } { 4 x ^2  } \Big| &  \leq  A \big(  |\la_{\de} - \la_{\de_0} | x L_{3, n} +  x^{ - 2  } R_x \big) \\
	&   \leq A \big(  x^{ -2 } R_x +  | \de (x)  - \de_{0} (x) | x^{ -1 } L_{2, n} \big).
\end{align*}
Thus we complete the proof of \eqref{e6}.
The result \eqref{e4} directly follows from (\ref{0-mo}), and \eqref{e5} follows from \eqref{e6} and \eqref{e4}. This completes the proof of Lemma \ref{lemma-2}.

\subsection{Proof of Lemma \ref{lemma-3}} \label{pf-lemma0}
In the sequel,  $ C, C_1, C_2, \ldots $ are positive constants that may depend on $ \omega $ and $ r_0 $ and may take different values at each appearance. Note that \eqref{error1} is a special case of \eqref{error1generalversion}.   Denote $ Z_i = Y_i^2 - \e Y_i^2 $ and $ G_i = 2 r x X_i - w r  x^2 Y_i^2  $. Regarding \eqref{error1generalversion}, we have
\begin{align}
	&\e \left\{\left(V_n^2-1\right)^2 e^{\sum_{i=1}^{n} G_i  }\right\} \label{in-1} \\
	& \hspace{0.8cm} =\e \bigg\{\Big(\sum_{i=1}^{n}Z_i\Big)^2 e^{\sum_{i=1}^{n} G_i } \bigg\} \nonumber \\
	& \hspace{0.8cm} = \e e^{\sum_{i=1}^{n} G_i } \bigg[\sum_{i=1}^n\frac{ \e Z_i^2e^{ G_i }}{\e e^{G_i }}
	+\sum_{i\ne j}\frac{\e Z_ie^{ G_i}}{\e e^{G_i}}\frac{\e Z_je^{G_j}}{\e e^{G_j}}\bigg] \nn \\
	&\hspace{0.8cm}  \leq  \e e^{\sum_{i=1}^{n} G_i} \bigg[\sum_{i=1}^n\frac{ \e Z_i^2e^{G_i}}{\e e^{G_i}} + \bigg( \sum_{i=1}^n \frac{\e Z_ie^{G_i}}{\e e^{G_i}} \bigg)^2  \bigg] .\nn
\end{align}
We first treat $ \e [ Z_i e^{G_i} ] $. Recall that $ 0 < r <   r_0 $ for some number $r_0$. For $x$ satisfying \eqref{G2}, $ x^2 \e Y_i^2  \leq 1/16 $ and $ c_0 x^2 \e Y_i^2 \leq 1/4 $. It follows from the basic inequality $ | e^{s} - 1 |  \leq |s| e^{s \lor 0} $ that
\begin{align}
	& x^2 \big| \e [ Z_i e^{G_i} \IF ( x | X_i | \leq 1 , x | Y_i  | \leq 1  ) ] \big|  \label{bound-mo1} \\
	& \hspace{0.5cm} = x^2 \big|  \e    \big[ Z_i \big(1  + O ( x | X_i |  + x^2 Y_i^2  )  \big) \IF  ( x | X_i | \leq 1 , x | Y_i  | \leq 1  ) \big]    \big| \nn \\
	&  \hspace{0.5cm}  \leq  C \big\{ x^2 \e [ (Y_i^2 + \e Y_i^2 ) \IF \{ x | Y_i  | > 1  \} ]  \nn \\
	& \hspace{1.5cm}  +   x^2  \e [ (Y_i^2 + \e Y_i^2 ) \IF \{ x | X_i  | > 1 , x | Y_i  | \leq  1  \} ] \nn \\
	& \hspace{1.5cm}  + x^2 \e [ (Y_i^2 + \e Y_i^2 )  ( x| X_i | + x^2 Y_i^2  ) \IF ( x| X_i | \leq 1, x | Y_i | \leq 1  ) ] \big\}  \nn \\
	&  \leq C x^3 ( \e | X_i|^3 + \e | Y_i |^3  + \e [ | X_i  Y_i^2 |  ]  ) + C x^4  (\e Y_i^2)^2     \nn  \\
	& \leq C x^3 ( \e |X_i|^3 + \e |Y_i |^3  ) ,\nn
\end{align}	
where the last inequality results from Lemma \ref{lemma-0} and the basic inequality $ a b^2 \leq a^3 + b^3  $ for $ a, b > 0 $. In addition, as $ G_i \leq 2 r < 2 $ when $ x | X_i  | \leq 1  $, we obtain for $x > 3  $
\begin{align}
	& x^2 | \e [ Z_i e^{G_i} \IF (x | X_i  | \leq 1 ,   x | Y_i | > 1 ) ] | \label{bound-mo2} \\
	&  \leq C x^2   \e [ Y_i^2 \IF ( x | Y_i  |  > 1  )  ] +  C x^2  \e [ Y_i^2  ] \P ( x | Y_i  | > 1 )  \nn \\
	& \leq C x^3 \e [ |  Y_i |^3 \IF ( x | Y_i | > 1  ) ] \leq C \delta_{x, i} \nn.
\end{align}	
Moreover, as $   x^2 \e Y_i^2 \leq  1/16 $,
\begin{align}
	& x^2 | \e [ Z_i e^{G_i} \IF  ( x | X_i  | > 1 ) |   \label{bound-mo3} \\
	& \leq x^2 \e [ Y_i^2 e^{ G_i  } \IF  ( x | X_i  | > 1  ) ] +  \e [ e^{ G_i  } \IF  ( x | X_i  | > 1  )  ] \nn.
\end{align}	
As for the second error term, we have
\begin{align*}
	& \e [e^{ G_i } \IF ( x | X_i  | > 1  ) ]  \\
	&   = \e [ e^{G_i} \IF (  x | X_i  | > 1, X_i \leq  0 )  ] + \e [e^{ G_i} \IF (  x | X_i  | > 1, X_i > 0 ) ] \nn\\
	&  \leq x^3 \e [ | X_i^3  | \IF( x | X_i | > 1 ) ] + \e [ e^{ G_i } \IF ( x  X_i   > 1 ) ] \nn.
\end{align*}	
Observe that when $ x X_i > 1 $,
\begin{align*}
	e^{ G_i  } = e^{ 2 r x  X_i - \omega r  x^2 Y_i^2  } \leq e^{ 2 r x X_i } \leq e^{ 2 x X_i }
\end{align*}	
and by Cauchy inequality for $0  < r < r_0 < \omega $,
\begin{align}
	2 r  x X_i  - \omega r x^2 Y_i^2           & \leq 2 (\omega r)^{1/2} x X_i - \omega r x^2 Y_i^2  \label{cauchy}  \\
	& \leq \frac {  x^2 X_i^2 } { x^2 Y_i^2 + c_0  x^2 \e Y_i^2 } +  \omega r  c_0 x^2 \e  Y_i^2 . \nn
\end{align}	
Therefore
\begin{align}
	&  \e [ e^{G_i}  \IF ( x |X_i| > 1 ) ] \label{mo-31} \\
	& \leq x^3 \e [ | X_i |^3 \IF ( x | X_i | > 1  ) ] +  C_1  \e \Big[ e^{ \min \big\{ \frac { X_i^2  } { Y_i^2 + c_0 \e Y_i^2  } , 2 x X_i\big \} } \IF ( x | X_i | > 1  ) \Big] \nn \\
	& \leq C_2 R_{x, i}. \nn
\end{align}
As for the first error term in \eqref{bound-mo3},  it follows from the basic inequality $ |x|^k e^{  - \delta  x }  \leq C(k, \delta)$ for some constant $ C(k, \delta) $ depending on $k$ and $ \delta $ that
\begin{align*}
	& x^2 \e [ Y_i^2  e^{ G_i } \IF (x | X_i  |  > 1 ) ] \\
	\leq&  C_1    \e [ e^{ 2 r x X_i  -  ( \omega/2  + r_0 /2 ) r x^2 Y_i^2 } \IF ( x | X_i  | > 1   ) ] \nn \\
	\leq&  C_1    x^3 \e [ x | X_i |^3 \IF ( x | X_i | > 1  ) ] \\
	& \hspace{0.5cm}  +   C_1     \e [ e^{ 2 r x X_i  -  ( \omega/2  + r_0 /2 ) r x^2 Y_i^2 } \IF ( x  X_i   > 1   ) ] \nn.
\end{align*}
As $ 0 < r < r_0 < \omega/2 + r_0 /2 $, by a similar procedure to \eqref{cauchy} we obtain
\begin{align}
	x^2 \e [ Y_i^2  e^{ G_i } \IF (x | X_i  |  > 1 )  ] \leq  C_1   R_{x, i  }. \label{mo-32}
\end{align}
Consequently, it follows from the bounds \eqref{bound-mo1}--\eqref{bound-mo3} and \eqref{mo-31}--\eqref{mo-32} that
\begin{align}
	| \e [ Z_i e^{ G_i } ] | \leq  C \big[ x ( \e | X_i |^3 + | Y_i |^3    ) + x^{ - 2} R_{x, i} \big].  \label{mo-1th}
\end{align}

By similar arguments, we obtain
\begin{align}
	\e [ Z_i^2 e^{ G_i  } ] \leq  C x^{ -4 } R_{x, i}.  \label{mo-2th}
\end{align}

In addition,
\begin{align}
	\e [ e^{ G_i   } ] & = \prod_{i = 1}^n \Big(  \e [ e^{ G_i  } \IF ( x | X_i  |  > 1 ) ] + \e [ e^{ G_i  } \IF ( x | X_i  | \leq 1 , x | Y_i  | \leq 1  ) ] \label{mo-0th} \\
	& \hspace{2cm}  + \e [ e^{ G_i  }  \IF ( x | X_i  | \leq 1 , x  | Y_i |  > 1  ) ] \Big) \nn .
\end{align}
Note that
\begin{align*}
	\e [ e^{ G_i  }  \IF ( x | X_i  | \leq 1 , x  | Y_i |  > 1  ) ] \leq e^{ 2 r_0 }  x^3 \e [ | Y_i |^3 \IF ( x | Y_i  | \geq 1  ) ]  \leq  C \delta_{x, i},
\end{align*}	
and \eqref{mo-31} has shown
\begin{align*}
	\e [ e^{ G_i  } \IF ( x | X_i  |  > 1 ) ] \leq C R_{x, i},
\end{align*}	
moreover, by inequality $ | e^{s} - ( 1 + s + s^2 /2 + s^3 / 6  ) | \leq s^4 e^{s \lor 0}$ and $ G_i \leq 2 r_0  $ when $  x | X_i  | \leq 1 $,
\begin{align}
	& \e [ e^{ G_i  } \IF ( x |X_i | \leq 1 , x | Y_i  | \leq 1 ) ]  \label{0th-expan} \\
	&   =   \e \big[  \big( 1 + G_i + G_i ^2 /2 + G_i ^3 / 6 + O( G_i^4 ) \big)  \IF ( x | X_i  | \leq 1 , x | Y_i   | \leq 1  ) \big]\nn  \\
	&  = 1 + ( 2 r^2 - \omega r^2  ) x^2 \e X_i^2 - 2 \omega r^2 x^3 \e [ X_i Y_i^2  ] + \frac { 4 } {3} r^3 x^3 \e X_i^3  + O (\delta_{x, i}).  \nn
\end{align}	
Under condition \eqref{G2}, we have
$$
\delta_{x, i}  \leq R_{x, i} \leq x^3 ( \e | X_i |^3 + \e | Y_i |^3 ) + r_{x, i} \leq c_1 + 1/64.
$$
As a result for small constant $c_1$,
\begin{align*}
	1/2 \leq  \e [ e^{ G_i  }   ]  \leq 3/2.
\end{align*}
Thus it follows by substituting \eqref{mo-1th} and \eqref{mo-2th} into \eqref{in-1} that
\begin{align}
	& \e \left\{\left(V_n^2-1\right)^2 e^{\sum_{i=1}^{n} G_i  }\right\} \\
	& \leq 2 \prod_{i = 1}^n \e [ e^{ G_i  } ] \Big(  x^{ - 4} R_{x} +  ( x L_{3, n} + x^{ - 2} R_{x} )^3  \Big)  \nn \\
	& \leq  A x^{ - 2} R_x \exp \Big\{ ( 2 r^2 - \omega r^2  ) x^2 \e X_i^2 - 2 \omega r^2 x^3 \e [ X_i Y_i^2  ] + \frac { 4 } {3} r^3 x^3 \e X_i^3  +  A \delta_{x, i}   \Big\}  \nn,
\end{align}
hence the desired result \eqref{error1generalversion} is derived. \eqref{error1} is a special case of \eqref{error1generalversion}, with $ r = \lambda_1 $ and $ \omega = 1 $.

Next we prove \eqref{error2}. We have by conditional Cauchy-Schwarz inequality that
\begin{align}
	&\sum_{i=1}^n \e \Big\{\Big|W_i\Big(Z_i^2+2Z_i
	\sum_{i\ne j} Z_j   \Big)\Big|e^{\lambda_1\sum_{i=1}^nW_i}\Big\}  \label{err-2} \\
	\le& \sum_{i=1}^n \e \left|W_iZ_i^2e^{\lambda_1\sum_{i=1}^nW_i}\right|+2\sum_{i=1}^n \e \Big[\Big|Z_iW_i\sum_{j\ne i}Z_j\Big|e^{\lambda_1\sum_{i=1}^nW_j }\Big] \nonumber \\
	\le& Q_1+Q_2,\nonumber
\end{align}
where
\begin{align*}
	Q_1&= \prod_{i = 1}^n \e [ e^{ \lambda_1 W_i } ] \sum_{i=1}^n\frac{ \e \left|W_iZ_i^2e^{\lambda_1W_i}\right|}{ \e e^{\lambda_1W_i}}, \\
	Q_2&=2\sum_{i=1}^n \e \left(|Z_iW_i|e^{\lambda_1W_i}\right)\Big\{ \e \Big[ \Big(\sum_{j\ne i}Z_j\Big)^2e^{\lambda_1\sum_{j\ne i}W_j} \Big] \Big\}^{\frac 1 2} \left[ \e e^{\lambda_1\sum_{j\ne i}W_j}\right]^{\frac 1 2}.
\end{align*}

Recalling \eqref{0-mo} and \eqref{e4}, we have
\begin{align}
	\prod_{i = 1}^n \e [ e^{ \lambda_1 W_i } ]  \leq  \exp \{ m ( \lambda_1 ) \} e^{ A R_{x} } \label{part1}
\end{align}
and when $ x  $ satisfies \eqref{G2} with small constant $ c_1 $,
\begin{align}
	1/2 \leq   \e [ e^{ \lambda_1 W_i } ] \leq 3 /2 .  \label{part2}
\end{align}
Further, through an analogous proof to \eqref{mo-1th}, we obtain
\begin{align*}
	&  x^2  \e [ | W_i Z_i | e^{\lambda_1 W_i  }  | ] \\
	& \leq   C \e \big[ (  x^3  | X_i  Y_i^2 | + x^4Y_i^4 + x^3 | X_i | \e Y_i^2  +  x^4 Y_i^2 \e Y_i^2 )  e^{ \lambda_1 W_i  }\big] \\
	& \leq C \e \big[ (  x^3  | X_i  Y_i^2 | + x^4 Y_i^4 + x | X_i |   +  x^2 Y_i^2  )  e^{ \lambda_1 W_i  } \IF ( x  X_i   > 1 ) \big]  \\
	& \hspace{0.5cm}  + C \e \big[ (  x^3  | X_i  Y_i^2 | + x^3 |Y_i|^3 + x^3 | X_i | \e Y_i^2  +  x^4 Y_i^2 \e Y_i^2 )  e^{ \lambda_1 W_i  } \IF ( x  X_i   \leq  1 ) \big] \\
	& \leq   C \e \big[ e^{ \frac {10 \lambda_1 } {9} x X_i  - \frac { 9 \lambda_1} { 10  }  x^2 Y_i^2 } \IF ( x  X_i  > 1  ) \big] + C x^3 ( \e | X_i |^3 + \e | Y_i |^3 ),
\end{align*}
where in the last inequality we used the fact that $ ab^2 \leq a^3 + b^3 $	for $  a> 0, b> 0$ and $ x^4 ( \e [Y_i^2] ) ^2 \leq 2 \delta_{x, i }  $ by Lemma \ref{lemma-0}.
By \eqref{cauchy}. we obtain for $ \lambda_{1} \leq 3/4 $,
\begin{align*}
	\e \big[ e^{ \frac {10 \lambda_1 } {9} x X_i  - \frac { 9 \lambda_1} { 10  }  x^2 Y_i^2 } \IF ( x | X_i  | > 1  ) \big]  \leq C \e \Big[ e^{ \min \big \{  \frac { X_i^2  }  { Y_i^2 + c_0 \e Y_i^2  } , 2 x X_i \big\} } \IF ( x  X_i  > 1  )  \Big]  .
\end{align*}	
Hence
\begin{align}
	\e [ | W_i Z_i | e^{ \lambda_1 W_i  }   ] \leq x^{- 2} R_{x, i} + x ( \e | X_i |^3 + \e | Y_i |^3 ). \label{part3}
\end{align}
Similarly, we can obtain
\begin{align}
	\e \big( | W_i Z_i^2  | e^{ \lambda_1 W_i  } \big) \leq C  x^{ - 4 } R_{x, i },  \label{part4}
\end{align}
which together with \eqref{part1} and  \eqref{part2} gives
\begin{align*}
	Q_1  \leq   C_1 \exp \{ m(\lambda_1) \} x^{ - 4 } R_{x} e^{ C_2 R_{x} }.
\end{align*}
As for $ Q_2 $, noting the bound \eqref{part1}--\eqref{part3} and \eqref{error1}, we have
\begin{align*}
	Q_2  & \leq  C  ( x^{ - 2} R_x  + x L_{3, n} ) \big( \exp \{ m ( \lambda_1 ) \} x^{ - 2} R_x e^{ A R_{x} }\big)^{ 1/2  } \big( \exp \{ m ( \lambda_1 ) \}  e^{ A R_x }  \big)^{ 1/2 } \\
	&  \leq C ( x^{ - 2 } R_x  + x L_{3, n} ) \exp \{ m ( \lambda_1 )  \} x^{^{ - 1 } }R_x^{ 1 /2  } e^{ A R_{x} }\\
	& \leq C \exp {m (\lambda_1)} x^{- 2} R_x e^{ A R_x }.
\end{align*}
Here we apply the fact $ x^4 L_{3, n}  \leq 2 \delta_x \leq 2 R_x $ in Lemma \ref{lemma-0} to derive the last inequality. Thus \eqref{error2} is proved. The proof of Lemma \ref{lemma-3} is completed.

\subsection{Proof of Lemma \ref{lemma-rxi}} \label{sec:lemma-rxi}
Assume $A=1$ without loss of generality, as the proof for general $A> 0$ is similar. By Cauchy-Schwarz inequality, it is obvious that $ S_n^2 \le n V_n^2 $ and hence $ \e [ e^{ \frac { S_n^2  } { V_n^2 + c_0 B_n^2  } }   ] < \infty $. We have
\begin{align}
	&  \e \Big[ e^{ \f {S_n^2} {V_n^2 + c_0 B_n^2}}  \IF  ( b |S_n| > 1 )   \Big]   \label{e-1} \\
	& \hspace{0.5cm} = \int_{0}^{n} e^t \, \P \Big( \f { |S_n| } { \sqrt{ V_n^2  + c_0 B_n^2 } } > \sqrt{t}, b |S_n| > 1 \Big)  d t \nn \\
	& \hspace{0.5cm} = \int_{0}^{n} e^t \, \Big[ \P \Big( \f { S_n } { \sqrt{ V_n^2  + c_0 B_n^2 } } > \sqrt{t}, b |S_n| > 1 \Big)  \nn   \\
	& \hspace{2.5cm}+ \P \Big( \f { -S_n } { \sqrt{ V_n^2  + c_0 B_n^2 } } > \sqrt{t}, b |S_n| > 1 \Big) \Big]  d t. \nn
\end{align}
The integral in $ 0< t <1 $ is bounded by $ e b^p \e | S_n |^p $  and hence is dominated by the right hand side of \eqref{rxi} by Rothenthal's inequality.
Thus it remains to consider the integral for $ t>1 $. Let $ \nu =   \f{ c B_n}{\sqrt{t}} $, where $ c $ is a large positive number. Denote $ \hatX_i = X_i \IF ( |X_i| \le \nu ) $, then
\begin{align}
	&  \P \Big( \f { S_n } { \sqrt{ V_n^2  + c_0 B_n^2 } } > \sqrt{t}, b |S_n| > 1 \Big)  \label{e-4} \\
	& \le \P \Big( \f {\hatS_n } { \sqrt{ V_n^2 + c_0 B_n^2 } } > \f{\sqrt{t}}{2} , b |S_n| >1  \Big)  \nn \\
	& \hspace{2cm} + \P \Big( \f {\sum_{i=1}^{n} X_i \IF ( |X_i| > \nu ) } { \sqrt{ V_n^2 + c_0 B_n^2 } } > \f{\sqrt{t}}{2} , b |S_n| >1  \Big) \nn \\
	& \le \P  \Big( \f{ \hatS_n } {\sqrt {c_0} B_n } > \f{\sqrt{t}}{2}, b |S_n| > 1  \Big) + \P  \Big(  \sum_{i=1}^{n} \IF (  |X_i| > \nu )> \f t4 , b |S_n| > 1  \Big) \nn \\
	& \le b^p e^{ - \f { c_0^{1/4} t }{ 2 } } \e  \Big[ |S_n|^p \,  \exp \big\{ c_0^{ - \f 14 } \sqrt{t} \f{\hatS_n} { B_n } \big\}  \Big] + b^p e^{-2t} \e  \Big [ |S_n|^p e^{ 8\sum_{i=1}^{n} \IF ( |X_i| > \nu )  }  \Big] \nn \\
	& :=  b^p e^{ - \f { c_0^{1/4} t }{ 2 } } E_1 + b^p e^{-2t} E_2 , \nn
\end{align}
where the second inequality is obtained by Cauchy inequality. Now we turn to estimating $  E_1 $ and $ E_2 $. Let $ A_{n, p} = \f{1}{2 (4e^8)^{1/p}} \max \{ B_n, (\sum_{i=1}^{n} \e |X_i|^p)^{1/p} \} $. Denote $ \bar{X_i} = X_i \IF  ( |X_i| \le A_{n, p} ) $, $\bar{S}_n = \sum_{i = 1}^n \bar{X}_i $ and $ S_n^{ (i) } = S_n - X_i $. Observing that
\begin{align}
	\IF ( |S_n| > x ) & \le \IF ( | \bar{S}_n | > x ) + \sum_{i=1}^{n} \IF  ( |S_n| > x, | X_i | > A_{n, p} )  \\
	& \le  \IF  ( | \bar{S}_n | > x ) + \sum_{i=1}^{n} \IF  ( |X_i| > \f 12 x)  \nn \\
	&  \hspace{1cm} + \sum_{i=1}^{n} \IF  ( |S_n^{(i)}| > \f 12 x, | X_i | >  A_{n, p} ) , \nn
\end{align}
we have
\begin{align}
	E_2 & = \int_{0}^{\infty} p x^{p-1} \e [ e^{8 \sum_{i=1}^n \IF  ( |X_i| > \nu ) } \IF ( |S_n| > x )  ]  d x   \label{E2} \\
	& \le \int_{0}^{\infty} p x^{p-1} \e [ e^{8 \sum_{i=1}^n \IF ( |X_i| > \nu ) } \IF ( |\bar{S}_n| > x ) ] d  x \nn \\
	& \quad  + \sum_{i=1}^n  \int_{0}^{\infty} p x^{p-1} \e [ e^{8 \sum_{i=1}^n \IF ( |X_i| > \nu ) } \IF  ( |X_i| > \f 12 x )  ] d  x \nn \\
	& \quad + \sum_{i=1}^n  \int_{0}^{\infty} p x^{p-1} \e [ e^{8 \sum_{i=1}^n \IF ( |X_i| > \nu ) } \IF ( S_n^{(i)} > \f 12 x, |X_i | > A_{n,p}  ) ] d x \nn \\
	& \le Q_1 + \sum_{i=1}^n   e^8 \e (e^{8 \sum_{j \neq i} \IF ( |X_j| > \nu ) }) \int_{0}^{\infty} p x^{p-1} \P (|X_i| > \f 12 x )   d x \nn \\
	& \quad + \sum_{i=1}^n e^8 ({A_{n,p}})^{-p} \e |X_i|^p \int_{0}^{\infty} p x^{p-1} \e [ e^{8 \sum_{j \neq i} \IF ( |X_j| > \nu ) } \IF ( S_n^{(i)} > \f 12 x ) ]  d  x  \nn \\
	& \le Q_1 + e^8 2^p \exp \{ \f t4 \}   \sum_{i=1}^n \e |X_i|^p + e^8 2^p A_{n,p}^{-p} (\sum_{i=1}^{n} \e |X_i|^p ) E_2  \nn \\
	& \le Q_1 + e^8 2^p \exp \{ \f t4 \}   \sum_{i=1}^n \e |X_i|^p + \f 14 E_2 , \nn
\end{align}
where
\begin{equation*}
	Q_1 = \int_{0}^{\infty} p x^{p-1} \e [ e^{8 \sum_{i=1}^n \IF ( |X_i| > \nu ) }  \IF ( |\bar{S}_n| > x )  ] d  x .
\end{equation*}
The second to last inequality is based on the fact that $ |x + u|^p $ is a convex function with respect to $x$ and that there exists some large constant $c$ such that for $ \nu =   \f{  c B_n}{\sqrt{t}} $,
\begin{align*}
	\e e^{8 \sum_{i=1}^{n} \IF ( |X_i| > \nu  ) } &\le \textstyle \prod_{i=1}^n [ 1 + e^8 \P ( |X_i| > \nu ) ] \nn \\
	& \le \exp \{  \textstyle  \sum_{i=1}^n e^8 \P ( |X_i| > \nu  ) \} \nn \\
	& \le \exp \{ \f{e^8}{\nu^2} B_n^2  \} \le \exp \{ \f t4 \}.
\end{align*}
As for $ Q_1 $,   it holds that
\begin{equation}
	Q_1 \le \int_{0}^{\infty} p x^{p-1} e^{- \f x A_{n,p}} \prod_{i=1}^{n} \e e^{ 8 \IF ( |X_i | > \nu  ) + \f {1}{A_{n,p}} \bar{X}_i}  d x . \label{e-3}
\end{equation}
Furthermore, as $ |\bar{X}_i| / A_{n,p} \leq 1  $,  we obtain from Taylor expansion that
\begin{align}
	\lefteqn{ \prod_{i=1}^{n} \e e^{ 8 \IF ( |X_i | > \nu  ) + \f {1} {A_{n,p}} \bar{X}_i} } \label{e-2} \\
	& \le \prod_{i=1}^{n} \e \left( 1 + e^8 \IF ( |X_i| > \nu ) \right) \left( 1 + \f {\bar{X}_i }{A_{n,p}} +   \f{e^2 \bar{X}_i^2}{ 2 A_{n,p}^2 } \right) \nn \\
	& \le C \exp \Big\{   \sum_{i=1}^{n} \P( |X_i| > \nu ) +  \f{1}{A_{n,p}}   \e |X_i| \IF  ( |X_i| > \nu ) \nn \\
	& \hspace{4cm} + \f {  1}{ A_{n,p}^2 } \e X_i^2 \IF ( |X_i| > \nu ) \Big\} \nn \\
	& \le C \exp \Big\{ \f{  B_n^2}{\nu^2} +  \f { B_n^2}{A_{n,p} \nu} \sum_{i=1}^n \e |X_i|^p  + \f{1}{2 A_{n,p}^2 } B_n^2 \Big\}  \le A \exp \big\{ \f t4 \big\}, \nn
\end{align}
which yields
\begin{align}
	Q_1 &\le A\, (A_{n,p})^p \,\exp\{ \f t4 \} \int_{0}^{\infty} p x^{p-1} e^{-x}  d  x \label{Q1} \\
	&\le A_1 \exp\{ \f t4 \} \Big( \sum_{i=1}^n \e |X_i|^p + (\e S_n^2)^{\f p2} \Big). \nn
\end{align}
Consequently, it follows from \eqref{E2} and \eqref{Q1} that
\begin{equation}
	E_2 \le A_2 \exp \{ \f t4 \}  \Big[ \sum_{i=1}^{n} \e |X_i|^p + (\e S_n^2)^{\f p2}   \Big]. \label{E-2}
\end{equation}

Next we deal with $ E_1 $. Recalling that $ | \hat{X}_i  |\leq \nu = C B_n / \nu $, we have for large constant $c_0$
\begin{align}
	\e \exp \Big\{ c_0^{- \f 14 } \sqrt{t} \f {\hatS_n}{B_n} \Big\}
	& \le \exp \bigg\{ \sum_{i=1}^n \Big[ \f{\sqrt{t}}{ c_0^{1/4} B_n } \e |X_i| \IF ( |X_i| > \nu )  \\
	& \hspace{1.5cm}+ \exp \{\f{C}{c_0^{1/4}}\} \f{t}{2 c_0^{1/2} B_n^2} \e X_i^2   \Big]  \bigg\} \nn \\
	& \le \exp \bigg\{  \f{\sqrt{t} B_n^2 }{ c_0^{\f14} B_n \nu } + \exp  \Big\{\f{C}{c_0^{1/4}} \Big\} \f{t}{2 c_0^{1/2}}  \bigg\}  \nn \\
	& \le \exp \big\{\f t4 \big\} \nn  .
\end{align}

Denote $ \hat{S}_n^{(i)} = \hat{S}_n - \hat{X}_i $ and
\begin{equation*}
	Q_2 =  \int_0^{\infty} p x^{p-1} e^ {- x/a} \e \Big[\exp \Big\{ c_0^{- \f 14 } \sqrt{t} \f {\hatS_n}{B_n} + \f {\bar{S}_n}{A_{n,p}}  \Big\}  \Big] d x .
\end{equation*}
By a similar procedure to \eqref{E2}--\eqref{E-2},   we have for large $c_0$
\begin{align*}
	E_1 & = \int_0^{\infty} p x^{p-1} \e \Big[\exp \Big\{ c_0^{- \f 14 } \sqrt{t} \f {\hatS_n}{B_n}  \Big\} \IF ( |S_n| > x )  ] dx \nn \\
	& \leq Q_2 + C_1 \sum_{i=1}^n \e \Big[ \exp \Big(  c_0^{- \f 14 } \sqrt{t} \f {\hatS_n}{B_n} \Big) \Big] \cdot \e |X_i|^p \nn \\
	& \quad + 2^p e^{c/c_0^{\f{1}{4}}} \sum_{i=1}^n \f{\e |X_i|^p}{A_{n,p}^{p}} \cdot  \e \Big[ |S_n^{(i)}|^p \exp \Big(  c_0^{- \f 14 } \sqrt{t} \f {\hatS_n^{(i)}}{B_n} \Big)  \Big] \nn  \\
	& \le Q_2 + C_2 \exp\{\f t4\} \sum_{i=1}^n \e |X_i|^p + \f 14 E_1.
\end{align*}
Moreover, it follows from Taylor expansion that
\begin{align}
	Q_2
	& \le C (A_{n,p})^p \prod_{i=1}^n  \e \Bigg[ \bigg( 1 + \f {\sqrt{t}}{ c_0^{1/4} B_n } \hatX_i +   \f{ t e^ { c / c_0^{1/4} } }{2 c_0^{1/2} B_n^2 }  \hatX_i^2 \bigg)\bigg( 1 + \f{ e \bar{X}_i}{A_{n,p}} \bigg) \Bigg] \nn \\
	& \le Ce^{t/4}(A_{n,p})^p \le C e^{t/4}\Big[ \sum_{i=1}^{n} \e |X_i|^p + (\e S_n^2)^{p/2}\Big]. \nn
\end{align}
Consequently,
\begin{equation*}
	E_1 \le  C\exp\{ \f t4 \} \Big[ \sum_{i=1}^{n} \e |X_i|^p + (\e S_n^2)^{p/2} \Big] ,
\end{equation*}
which together with \eqref{E-2}, \eqref{e-4} and \eqref{e-1} gives the result \eqref{rxi}. The proof is completed.    	
%\bibliographystyle{imsart-nameyear} % Style BST file (imsart-number.bst or imsart-nameyear.bst)
%\bibliography{general_CTMD}       % Bibliography file (usually '*.bib')

\end{document}